\documentclass[12pt]{article}

\newenvironment{wileykeywords}{\textsf{Keywords:}\hspace{\stretch{1}}}{\hspace{\stretch{1}}\rule{1ex}{1ex}}

\usepackage{amsmath,amssymb,amsthm}
\newtheorem{thm}{Theorem}
\usepackage{graphicx}
\usepackage{color}
\usepackage{dcolumn}
\usepackage{bm}
\usepackage[numbers,super,comma,sort&compress]{natbib}
\usepackage{hyperref}
\usepackage{stmaryrd}
\usepackage{multirow}

\usepackage{calc}
\usepackage[percent]{overpic}
\usepackage{cellspace, makecell}

\usepackage{float}

\usepackage{subcaption}

\usepackage{enumitem}

\usepackage{tikz}
\usetikzlibrary{quotes,angles}
\usetikzlibrary{positioning,fit,calc}
\tikzset{block/.style={draw,thick,text width=2cm,minimum height=1cm,align=center},
	line/.style={-latex}
}
\usetikzlibrary{calc,trees,positioning,arrows,chains,shapes.geometric,%
	decorations.pathreplacing,decorations.pathmorphing,shapes,%
	matrix,shapes.symbols}
\tikzset{
	>=stealth',
	punktchain/.style={
		rectangle, 
		draw=black, very thick,
		text width=4.2cm, 
		text centered, 
		on chain},
	line/.style={draw, thick, <-}                    , %
	    every join/.style={->, thick,shorten >=1pt}  , %
	    tuborg/.style={decorate}                     , %
	    tubnode/.style={midway, right=2pt} 
		}

\newcommand{\imagespacing}{\vspace{3pt}}

\definecolor{background-color}{gray}{0.98}

\title{Efficient mesh refinement for the Poisson-Boltzmann equation with boundary elements}
\author{Vicente Ramm\thanks{Departamento de Ingenier\'ia Mec\'anica, Universidad T\'ecnica Federico Santa Mar\'ia}, Jehanzeb H. Chaudhry\thanks{Department of Mathematics and Statistics, University of New Mexico}, Christopher D. Cooper\thanks{Departmento de Ingenier\'ia Mec\'anica and Centro Cient\'ifico Tecnol\'ogico de Valpara\'iso (CCTVal), Universidad T\'ecnica Federico Santa Mar\'ia, Valpara\'iso, Chile}}

\begin{document}

\maketitle

\begin{abstract}
The Poisson-Boltzmann equation is a widely used model to study the electrostatics in molecular solvation. Its numerical solution using a boundary integral formulation requires a mesh on the molecular surface only, yielding accurate representations of the solute, which is usually a complicated geometry. Here, we utilize adjoint-based analyses to form two goal-oriented error estimates that allows us to determine the contribution of each discretization element (panel) to the numerical error in the solvation free energy. This information is useful to identify high-error panels to then refine them adaptively to find optimal surface meshes. We present results for spheres and real molecular geometries, and see that elements with large error tend to be in regions where there is a high electrostatic potential. We also find that even though both estimates predict different total errors, they have similar performance as part of an adaptive mesh refinement scheme. Our test cases suggest that the adaptive mesh refinement scheme is very effective, as we are able to reduce the error one order of magnitude by increasing the mesh size less than 20\%. This result sets the basis towards efficient automatic mesh refinement schemes that produce optimal meshes for solvation energy calculations. 
\end{abstract}

\begin{wileykeywords}
Poisson-Boltzmann, Implicit solvent, Goal-oriented adjoint based error estimation, Boundary element method, Adaptive mesh refinement.
\end{wileykeywords}

\clearpage


\begin{figure}[h]
\centering
\colorbox{background-color}{
\fbox{
\begin{minipage}{1.0\textwidth}
\includegraphics[width=50mm,height=50mm]{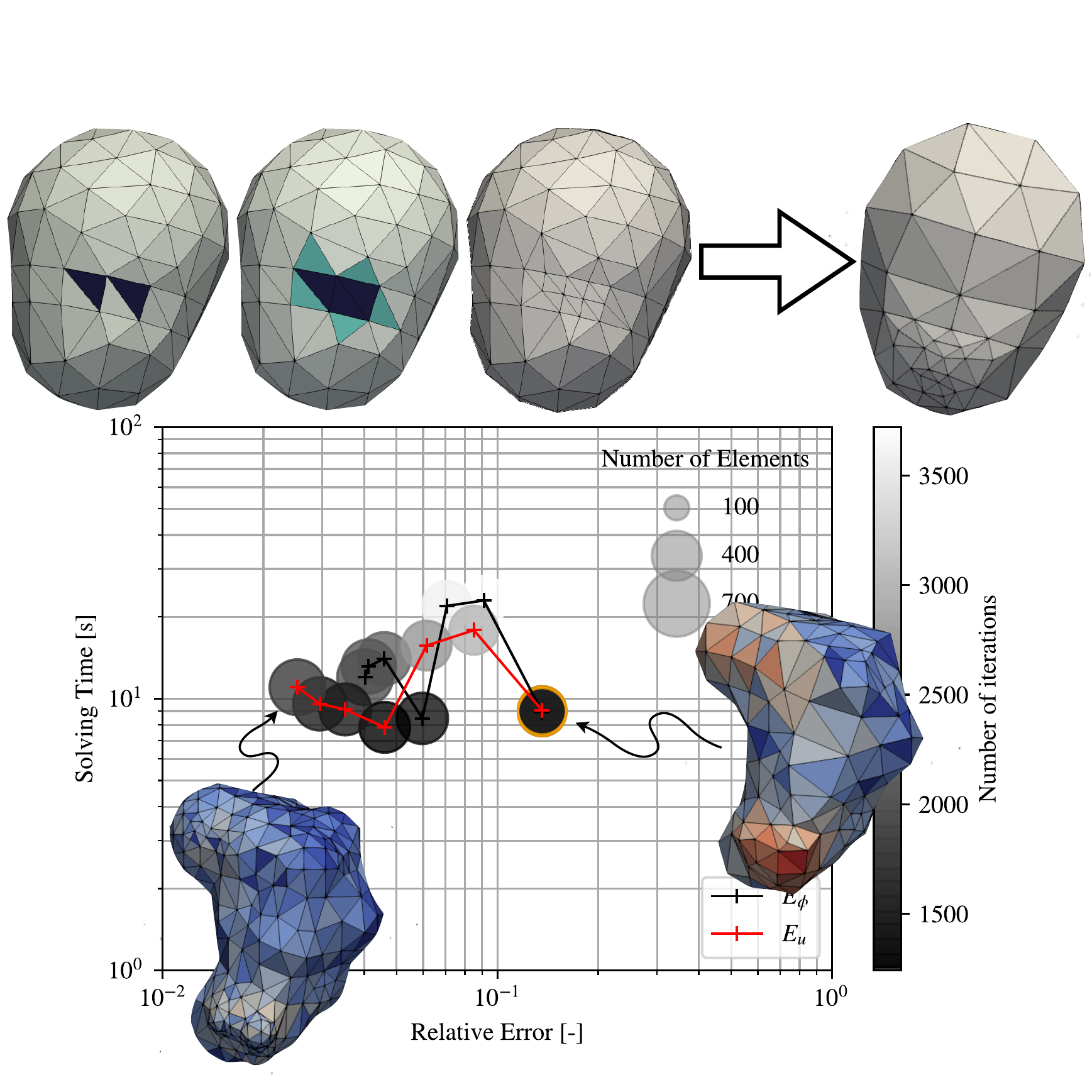} 
\\
The Poisson-Boltzmann equation is widely used to compute solvation energies of molecules. It considers the solute as a cavity region in an infinite dielectric, interfaced by the molecular surface. The boundary element method offers an efficient numerical solution, as it discretizes the interface only. Here, we present an \textit{a posteriori} error estimation method that detects high error elements, with which we generate a highly effective adaptive mesh refinement technique. 
\end{minipage}
}}
\end{figure}

  \makeatletter
  \renewcommand\@biblabel[1]{#1.}
  \makeatother

\bibliographystyle{apsrev}

\renewcommand{\baselinestretch}{1.5}
	
\normalsize

\clearpage

\section*{\sffamily \Large INTRODUCTION} 
In biological settings, biomolecules are found immersed in water with salt, and an appropriate treatment of the solvent is key to have accurate models. 
A popular treatment considers the solvent implicitly, where rather than accounting for each water molecule in a discrete way, they are represented as a continuum material~\cite{RouxSimonson1999,Baker2004}.
The implicit solvent model represents molecules in solution with continuum electrostatic theory, where a solute (region $\Omega_m$) is surrounded by an infinite dielectric (region $\Omega_w$), as sketched by figure \ref{fig:Solvation-impl}.
The solute is a region without water or salt, where the dielectric constant is low ($\epsilon_m=2-4$), and contains partial charges that are represented as point-delta functions.
Outside the solute, we use the permittivity of water ($\epsilon_w=80$), and consider the presence of salt ($\kappa$).
These two regions are interfaced by a molecular surface or interface ($\Gamma$), where several definitions are possible \cite{HarrisBoschitcshFenley2013}: solvent-excluded, solvent-accessible, Gaussian, or van der Waals surface.
In this work, we use the solvent-excluded surface\cite{Connolly83}, which is the result of rolling a spherical probe of the size of a water molecule ($1.4$ \AA radius), and tracking the contact points between the probe and the atoms of the solute (with their corresponding van der Waals radii).

This implicit solvent approximation dramatically reduces the number of degrees of freedom compared to molecular dynamics, and yields a partial differential equation (PDE) based model.
If we consider that, at equilibrium, the mobile ions from the salt in the solvent arrange according to Boltzmann's distribution, electrostatic theory results in the so-called Poisson-Boltzmann equation (PBE).
A common quantity of interest is the solvation free energy, which is the work required to bring the solute molecule from vacuum into its solvated state. The focus of this article is then twofold:
\begin{enumerate}[nosep,label=\roman*]
	\item Derive accurate a posteriori estimates for  the error in a boundary element method approximation of  the solvation free energy.
	\item Design an adaptive mesh refinement algorithm to efficiently arrive at accurate solutions.
\end{enumerate}

The PBE has been solved numerically with a wide variety of techniques, such as finite difference\cite{HonigNicolls1995,BakerETal2001,ChenETal2011}, finite element\cite{HolstETal2000,BakerETal2000,BondETal2010}, and boundary element (BEM)\cite{LuETal2006,Bajaj2011,GengKrasny2013,CooperBardhanBarba2014} methods.
In particular, BEM uses a surface integral formulation\cite{Shaw1985,YoonLenhoff1990,JufferETal1991}, where only the solute-solvent interface is discretized, and the potential goes to zero at infinity by construction. 
Then, the molecular surface is accurately represented, making BEM favorable for high precision simulations\cite{CooperBardhanBarba2014}.
However, BEM generates dense matrices that need fast methods to access large problems, such as fast multipole methods\cite{GreengardRokhlin1987,GreengardHuang2002}, treecodes\cite{DuanKrasny2001,LiJohnstonKrasny2009}, or hierarchical matrices\cite{Bebendorf2000,HoGreengard2012}.

Numerical approximations to the PBE often have large error. Hence, for the reliable use of the PBE in science and engineering, this error needs to be quantified. The tool used to accomplish this task in this article is adjoint based error estimation.  
In this technique, one solves an adjoint problem whose solution provides the
residual weighting to produce the error in the goal functional. The resulting estimate
 also allows to identify the  contributions to the error due to an individual mesh element (panel) and hence aid in forming adaptive algorithms.

Adjoint based analysis has been used for the error estimation of a variety of numerical methods and differential equations~\cite{giles_suli, becker_rannacher_2003, Estep1995, eehj_book_96}, for example, finite element methods~\cite{Estep:Larson:00,Chaudhry2019,BangerthRannacher,AO2000}, finite volume methods~\cite{barth04}, numerous time-integration schemes~\cite{CEG+2015,collins:2015,Logg2004,EJL04,Chaudhry2017}, and parallel-in-time and domain decomposition methods\cite{JESVS15,Chaudhry2019b}. 
A posteriori analysis of the finite element method for the PBE has been considered previously~\cite{ABC+2011,chaudhry_2018}, however, this work is the first such analysis for the PBE with BEM. Residual based a posteriori analysis for BEM has been studied earlier~\cite{CoHe14,NoKa15}, however, they focused on error in some global norm of the solution, whereas here we focus on quantifying the error in a goal functional or quantity of interest. Moreover, if the aim of the computation is to minimize the error in the goal, then forming per-element contributions of the error in the goal is a crucial ingredient in designing an adaptive algorithm. We achieve this aim by combining information from the numerical solution of the PBE with adjoint solutions to classify mesh elements that contribute the most towards the error in the goal. These elements are then refined to decrease the error in numerical approximations of the goal functional.

In the next section, we present the methodology, including the equations governing the implicit solvent model, the boundary integral formulation for the PBE, and finally derive two error estimates using adjoint based error analysis. 
The error estimates are also used to form per-element error indicator to devise a local mesh refinement strategy. 
In the Results and Discussion section we test the performance of the error estimates and the mesh refinement strategy on a variety of molecular setups. 
Finally, the last section presents the conclusions and outlook for future work.

\section*{\sffamily \Large METHODOLOGY}
\section*{\sffamily \Large The Poisson-Boltzmann equation and Solvation Free Energy}

BEM is limited to linear problems, however, the linearized version of the PBE is sufficiently accurate for most protein problems, where charges and potentials are not high~\cite{DecherchiETal2015}.
The linearized PBE is
\begin{equation}
-\nabla\cdot(\epsilon(\mathbf{r})\nabla u(\mathbf{r})) + \bar{\kappa}^2(\mathbf{r})u(\mathbf{r}) =  \sum_{k=1}^{N_q} q_k\delta(|\mathbf{r}-\mathbf{r}_k|)  \label{eq:pbe}
\end{equation}
where $q_k$ is the partial charge of atom $k$ (located at $\mathbf{r}_k$), $\epsilon(\mathbf{r})$ is the permittivity, and $\bar{\kappa}(\mathbf{r}) = \kappa\sqrt{\epsilon_w}$ is the modified Debye-H\"uckel parameter, with $\kappa$  the inverse of the Debye length.
Here, we use $\kappa=0.125$ \AA$^2$ in the solvent region, which corresponds to NaCl dissolved in water at a concentration of 150mM, while we set $\kappa = 0$ in the solute region to indicate ion exclusion there.
The dielectric constant 
also takes two values, $\epsilon_w$ or $\epsilon_m$, depending on the region ($\Omega_w$ or $\Omega_m$).
On the interface ($\Gamma$), the potential and electric displacement are continuous:
\begin{align}\label{eq:interface_condition}
\llbracket u\rrbracket_\Gamma=0, \qquad \left \llbracket \epsilon \frac{\partial u}{\partial\mathbf{n}} \right\rrbracket_\Gamma = 0,
\end{align}
where $\llbracket \psi \rrbracket_\Gamma= \lim_{\alpha\to 0^{+}} (\psi(x+\alpha \mathbf{n}) - \psi(x-\alpha \mathbf{n} ))$ is the jump in $\psi$ across the interface, for $x \in \Gamma$,  and $\mathbf{n}$ a normal vector pointing out of $\Omega_m$.

The electrostatic component of the solvation free energy can be computed as
\begin{equation}\label{eq:Gsolv}
\Delta G_{solv} = \frac{1}{2}\int_{\Omega} \rho(\mathbf{r})u_r(\mathbf{r})d\mathbf{r} = \frac{1}{2}\sum_{k=1}^{N_q} q_ku_r(\mathbf{r}_k) 
\end{equation}
where $u_r$ is the so called reaction potential, due to the polarization of the solvent.
Physically, the reaction potential $u_r$ is the difference in electrostatic potential between the isolated and dissolved states (see figure \ref{fig:ur}). 
In the isolated state, the electrostatic potential in the solvent is exactly zero, whereas in the molecule it is the Coulomb potential from the partial charges: 
\begin{equation}\label{eq:u_c}
u_c(\mathbf{r}) = \frac{1}{\epsilon_m}\sum_{k=1}^{N_q} \frac{q_k}{4\pi|\mathbf{r}-\mathbf{r}_k|},
\end{equation}
Then, when we take the difference in electrostatic potential, we only need to subtract out $u_c$ in the solute region, giving
\begin{equation}\label{eq:u_r_cases}
u_r(\mathbf{r}) = 
\begin{cases}
u - u_c &\text{ if $\mathbf{r}$ in } \Omega_m\\
u &\text{  if $\mathbf{r}$ in } \Omega_w.
\end{cases}
\end{equation}
Hence, the equation, boundary condition, and jump conditions at the interface $\Gamma$ for $u_r$ are
\begin{align}\label{eq:ur_jump}
-\nabla\cdot \left(\epsilon(\mathbf{r}) \nabla u_r(\mathbf{r})\right) + \overline{\kappa}^2(\mathbf{r})u_r(\mathbf{r}) & = 0 \nonumber \\
u_r(\infty) &= 0\nonumber\\
\left \llbracket  u_r \right \rrbracket_\Gamma &= u_c\left(\mathbf{r}_\Gamma\right)  \nonumber\\
\left \llbracket \epsilon(\mathbf{r}) \frac{\partial u_r}{\partial \mathbf{n}} \right  \rrbracket_\Gamma &= \epsilon_m \frac{\partial u_c}{\partial \mathbf{n}}\left(\mathbf{r}_\Gamma\right)
\end{align}
%
In practice,  there is no need to calculate $u_r$ in $\Omega_w$ for $\Delta G_{solv}$ in equation \eqref{eq:Gsolv}.

\section*{\sffamily \Large A boundary integral formulation of the Poisson-Boltzmann equation}

Equation \eqref{eq:pbe} is effectively a coupled system of constant coefficient partial differential equations, where the Poisson-Boltzmann equation governs in $\Omega_w$ and the Poisson equation with point-charge sources in $\Omega_m$.
%
Using Green's second identity, the electrostatic potential anywhere can be computed as
\begin{align}\label{eq:bie_domain}
u(\mathbf{r}) &=-\oint_\Gamma u^-(\mathbf{r'})\frac{\partial}{\partial\mathbf{n}}\left(\frac{1}{4\pi|\mathbf{r}-\mathbf{r}'|}\right)d\mathbf{r}' + \oint_\Gamma \frac{\partial u^-(\mathbf{r}')}{\partial\mathbf{n}}\frac{1}{4\pi|\mathbf{r}-\mathbf{r}'|}d\mathbf{r}' + \frac{1}{\epsilon_m}\sum_{k=1}^{N_q} \frac{q_k}{4\pi|\mathbf{r}-\mathbf{r}_k|} \text{ in $\Omega_m$}, \nonumber \\
u(\mathbf{r}) &= \oint_\Gamma u^+(\mathbf{r'})\frac{\partial}{\partial\mathbf{n}}\left(\frac{e^{-\kappa|\mathbf{r}-\mathbf{r}'|}}{4\pi|\mathbf{r}-\mathbf{r}'|}\right)d\mathbf{r}' + \oint_\Gamma \frac{\partial u^+(\mathbf{r}')}{\partial\mathbf{n}}\frac{e^{-\kappa|\mathbf{r}-\mathbf{r}'|}}{4\pi|\mathbf{r}-\mathbf{r}'|}d\mathbf{r}' \text{ in $\Omega_w$}, 
\end{align}
where the superscripts $``-"$ and $``+"$ indicate that the potential is being evaluated on the internal and external side of $\Gamma$, respectively.         
Evaluating $\mathbf{r}$ on $\Gamma$ and enforcing the interface conditions in equation \eqref{eq:interface_condition}, this becomes~\cite{YoonLenhoff1990}
\begin{align}
\frac{u^-(\mathbf{r})}{2} + \oint_\Gamma u^-(\mathbf{r'})\frac{\partial}{\partial\mathbf{n}}\left(\frac{1}{4\pi|\mathbf{r}-\mathbf{r}'|}\right)d\mathbf{r}' - \oint_\Gamma \frac{\partial u^-(\mathbf{r}')}{\partial\mathbf{n}}\frac{1}{4\pi|\mathbf{r}-\mathbf{r}'|}d\mathbf{r}' &= \frac{1}{\epsilon_m}\sum_{k=1}^{N_q} \frac{q_k}{4\pi|\mathbf{r}-\mathbf{r}_k|}, \nonumber \\
\frac{u^-(\mathbf{r})}{2} - \oint_\Gamma u^-(\mathbf{r'})\frac{\partial}{\partial\mathbf{n}}\left(\frac{e^{-\kappa|\mathbf{r}-\mathbf{r}'|}}{4\pi|\mathbf{r}-\mathbf{r}'|}\right)d\mathbf{r}' +\frac{\epsilon_m}{\epsilon_w} \oint_\Gamma \frac{\partial u^-(\mathbf{r}')}{\partial\mathbf{n}}\frac{e^{-\kappa|\mathbf{r}-\mathbf{r}'|}}{4\pi|\mathbf{r}-\mathbf{r}'|}d\mathbf{r}'  &=0. \label{eq:bie}
\end{align}
where integrals are now principal value integrals. 

To compute $u_r$ in $\Omega_m$, we use equation \eqref{eq:bie_domain} and subtract out $u_c$
\begin{equation}\label{eq:u_reac}
u_r(\mathbf{r}) = u^- - u_c = -\oint_\Gamma u^-(\mathbf{r'})\frac{\partial}{\partial\mathbf{n}}\left(\frac{1}{4\pi|\mathbf{r}-\mathbf{r}'|}\right)d\mathbf{r}' + \oint_\Gamma \frac{\partial u^-(\mathbf{r}')}{\partial\mathbf{n}}\frac{1}{4\pi|\mathbf{r}-\mathbf{r}'|}d\mathbf{r}'. 
\end{equation}

We solve the system in equation \eqref{eq:bie} numerically using a boundary element method (BEM), available in the \texttt{bempp} library~\cite{SmigajETal2015}.
The BEM uses a triangulated surface to generate a finite dimensional representation of $u^-$, which we call $U^-$, on panel $p$, of the form
\begin{equation}
U^-(\mathbf{r}_p) = \sum_{l=1}^{N_l} \Psi(\mathbf{r}_p,\mathbf{r}_l) U^-(\mathbf{r}_l)
\end{equation}
where $\Psi$ is the shape function, which in this work may belong to the space of piecewise constant (for equation \eqref{eq:bie}) or linear (to solve the adjoint in equation \eqref{eq:adjoint}) functions.
The \texttt{bempp} library uses a Galerkin discretization to arrive to a linear system for $U^-$ and $\frac{\partial U^-}{\partial\mathbf{n}}$, which is solved using GMRES.
Then, we replace these results in equation \eqref{eq:u_reac} to compute $U_r$ (the numerical approximation of $u_r$), and then  the numerical approximation of solvation energy $\Delta \widehat{G}_{solv}$ as
\begin{equation}\label{eq:Gsolv_hat}
\Delta \widehat{G}_{solv} = \sum_{k=1}^{N_q} q_kU_r(\mathbf{r}_k)
\end{equation}

In this BEM implementation, the surface $\Gamma$ is discretized in $N_p$ flat triangular panels using \texttt{msms}~\cite{SannerOlsonSpehner1995} or \texttt{Nanoshaper}~\cite{DecherchiRocchia2013}, we assume a piecewise constant ansatz, and compute integrals with Gaussian quadrature rules, to obtain the numerical approximations $U^-$ and $\frac{\partial U^-}{\partial \mathbf{n}}$ on the interface.
Also, we set the GMRES tolerance to 10$^{-8}$ in all tests.

\section*{\sffamily \Large Goal-oriented error estimation}
In this work, we use an adjoint-based error estimation method, which allows us to approximate the contribution of each element to the error in a goal or quantity of interest (QoI).
In general, let us consider the QoI as
\begin{equation}\label{eq:qoi}
\text{QoI} = \int_\Omega \psi(\mathbf{r})u_r(\mathbf{r})d\mathbf{r},
\end{equation}
where $\psi(\mathbf{r})$ is a weight function chosen to specify the QoI. Here, the QoI is the solvation free energy ($\Delta G_{solv}$). Comparing equation \eqref{eq:qoi} with the expression for
the solvation free energy in equation \eqref{eq:Gsolv} we see that $\psi(\mathbf{r})$ is the charge distribution,
\begin{equation}\label{eq:psi_rho}
	\psi(\mathbf{r}) = \rho(\mathbf{r}).
\end{equation}

Note that we ignored the factor $1/2$ present in equation \eqref{eq:Gsolv} when defining $\psi$. It is trivial to account for this constant factor when forming error estimates, and the error estimates we form later do indeed account for this.

\subsection*{\sffamily \large The adjoint operator}

Given a differential operator $\mathcal{D}$, the adjoint operator, $\mathcal{D}^\ast$, is defined as,
\begin{equation}
\int_\Omega \mathcal{D} w(\mathbf{r})v(\mathbf{r})d\mathbf{r} = \int_\Omega w(\mathbf{r})\mathcal{D}^*v(\mathbf{r})d\mathbf{r},
\end{equation}
where $w$ and $v$ are functions for $\mathbf{r} \in \Omega$. 
We can write $\psi(\mathbf{r})$ in terms of an adjoint function ($\phi$) as
\begin{equation}
\mathcal{D}^\ast\phi(\mathbf{r}) = \psi(\mathbf{r}).
\end{equation}
%
The corresponding differential operator on the adjoint is \cite{ABC+2011}
\begin{align}
\label{eq:adjoint}
\mathcal{D}^\ast\phi(\mathbf{r}) &= -\nabla\cdot\left(\epsilon(\mathbf{r})\nabla \phi(\mathbf{r})\right) + \overline{\kappa}^2 \phi(\mathbf{r}) = \psi(\mathbf{r})  \quad \text{ in } \Omega,
\end{align}
which implies the following jump conditions at the interface
\begin{align}
\left \llbracket\phi(\mathbf{r})\right  \rrbracket_{\Gamma} &= 0, \nonumber  \\
\left \llbracket \epsilon(\mathbf{r})\frac{\partial\phi}{\partial\mathbf{n}}(\mathbf{r})\right  \rrbracket_{\Gamma} &= 0. \label{eq:phi_jump}
\end{align}
%

%
%

\subsection*{\sffamily \large Exact Error Representations}

Let $e_r = u_r - U_r$ denote the error in the reaction potential. Our aim is to compute the error in the numerical approximation to the solvation energy,
\begin{equation}
	\frac{1}{2}\int_{\Omega} e_r(\mathbf{r}) \psi(\mathbf{r}) d\mathbf{r} = \frac{1}{2}\int_{\Omega} (u_r(\mathbf{r}) - U_r(\mathbf{r})) \psi(\mathbf{r}) d\mathbf{r}. \label{eq:error_approx}
\end{equation}
Here, we present two alternatives to build error estimates. The proofs of the theorems are given in the Appendix.

\begin{thm}
The error in the approximation to the solvation free energy is,
\begin{align}
	\label{eq:err_rep_I}
	\frac{1}{2}\int_\Omega e_r(\mathbf{r}) \psi(\mathbf{r}) d\mathbf{r} = E_\phi + R_\phi	
\end{align}
where 
\begin{align*}
E_\phi &= \frac{\epsilon_{m}}{2} \oint_{\Gamma} \left(  \frac{\partial \phi^-}{\partial \mathbf{n}} (\mathbf{r})u_c(\mathbf{r}) - \phi^-(\mathbf{r})\frac{\partial u_c}{\partial \mathbf{n}}(\mathbf{r}) \right)d\mathbf{r} + \frac{\epsilon_{m}}{2}\oint_{\Gamma}\left( \frac{\partial \phi^-}{\partial \mathbf{n}}(\mathbf{r}) U_r^-(\mathbf{r}) - \phi^-(\mathbf{r}) \frac{\partial U_r^-}{\partial \mathbf{n}}(\mathbf{r}) \right)d\mathbf{r},\nonumber\\
R_\phi &= \frac{1}{2} \int_{\Omega_w \cup \Omega_m} \phi(\mathbf{r}) (-\nabla\cdot(\epsilon(\mathbf{r}) \nabla u_r(\mathbf{r}))+ \overline{\kappa}^{2}(\mathbf{r})u_r(\mathbf{r})) d\mathbf{r}
 +\frac{\epsilon_{m}}{2}\int_{\Omega_{m}} \phi(\mathbf{r}) \nabla^{2}U_r(\mathbf{r})d\mathbf{r}.
\end{align*}

\label{th:Ephi}

\end{thm}

Alternatively, we derive a different error representation.

\begin{thm}
The error in the approximation to the solvation free energy is,
\begin{align}
	\label{eq:err_rep_II}
	\frac{1}{2}\int_\Omega e_r(\mathbf{r}) \psi(\mathbf{r}) d\mathbf{r}  = E_u + R_u	
\end{align}
where 
\begin{align*}
E_u &= \frac{\epsilon_{m}}{2} \oint_{\Gamma} \left(  \frac{\partial \phi^-}{\partial \mathbf{n}} (\mathbf{r})u_c(\mathbf{r}) - \phi^-(\mathbf{r})\frac{\partial u_c}{\partial \mathbf{n}}(\mathbf{r}) \right)d\mathbf{r} - \frac{\epsilon_{m}}{2}\oint_{\Gamma}\left( u_c(\mathbf{r}) \frac{\partial U^-}{\partial \mathbf{n}}(\mathbf{r})-\frac{\partial u_c}{\partial \mathbf{n}}(\mathbf{r}) U^-(\mathbf{r})  \right)d\mathbf{r},\nonumber\\
R_u &= \frac{1}{2} \int_{\Omega_w \cup \Omega_m} \phi(\mathbf{r}) (-\nabla\cdot(\epsilon(\mathbf{r}) \nabla u_r(\mathbf{r}))+ \overline{\kappa}^{2}(\mathbf{r})u_r(\mathbf{r})) d\mathbf{r}.
\end{align*}

\label{th:Eu}
\end{thm}

\subsection*{\sffamily \large Error Estimates}

The error representations in  equations \eqref{eq:err_rep_I} and \eqref{eq:err_rep_II} contain surface integrals (represented by $E_\phi$ and $E_u$) and volume integrals (represented by $R_\phi$ and $R_u$). In the context of BEM, evaluating volume integrals is computationally expensive. Hence, we propose the following error estimates,
\begin{align}
\label{eq:err_est_from_err_rep}
	\int_\Omega (u_r(\mathbf{r}) - U_r(\mathbf{r})) \psi(\mathbf{r}) d\mathbf{r}  &\approx E_\phi, \text{ and} \nonumber\\
	\int_\Omega (u_r(\mathbf{r}) - U_r(\mathbf{r})) \psi(\mathbf{r}) d\mathbf{r}  &\approx E_u,
\end{align}
where we neglect the contributions from the volumetric terms $R_\phi$ and $R_u$. Our numerical results will be useful to determine if these approximations are appropriate.


\subsubsection*{\sffamily \normalsize Element-wise error estimation}


%
Our aim is to find the contribution of each discretization element to the error in the numerical solution of $\Delta G_{solv}$. To this end, we decompose the surface integrals in $E_\phi$ and $E_u$ as follows,
\begin{equation}
\begin{aligned}
|E_\phi| &= \left|  \frac{\epsilon_{m}}{2} \oint_{\Gamma} \left(  \frac{\partial \phi^-}{\partial \mathbf{n}} (\mathbf{r})u_c(\mathbf{r}) - \phi^-(\mathbf{r})\frac{\partial u_c}{\partial \mathbf{n}}(\mathbf{r}) \right)d\mathbf{r} + \frac{\epsilon_{m}}{2}\oint_{\Gamma}\left( \frac{\partial \phi^-}{\partial \mathbf{n}}(\mathbf{r}) U_r^-(\mathbf{r}) - \phi^-(\mathbf{r}) \frac{\partial U_r^-}{\partial \mathbf{n}}(\mathbf{r}) \right)d\mathbf{r} \right|, \\
&\leq \sum_i^{N_p} \left| \frac{\epsilon_{m}}{2} \oint_{\Gamma_i} \left(  \frac{\partial \phi^-}{\partial \mathbf{n}} (\mathbf{r})u_c(\mathbf{r}) - \phi^-(\mathbf{r})\frac{\partial u_c}{\partial \mathbf{n}}(\mathbf{r}) \right)d\mathbf{r} + \frac{\epsilon_{m}}{2}\oint_{\Gamma_i}\left( \frac{\partial \phi^-}{\partial \mathbf{n}}(\mathbf{r}) U_r^-(\mathbf{r}) - \phi^-(\mathbf{r}) \frac{\partial U_r^-}{\partial \mathbf{n}}(\mathbf{r}) \right)d\mathbf{r} \right|,  \\
 &=\sum_i^{N_p} E_\phi^i,
 \end{aligned}
 \label{eq:error_element_I}
 \end{equation}
 and
 \begin{equation}
 \begin{aligned}
|E_u| &= \left|  \frac{\epsilon_{m}}{2} \oint_{\Gamma} \left(  \frac{\partial \phi^-}{\partial \mathbf{n}} (\mathbf{r})u_c(\mathbf{r}) - \phi^-(\mathbf{r})\frac{\partial u_c}{\partial \mathbf{n}}(\mathbf{r}) \right)d\mathbf{r} - \frac{\epsilon_{m}}{2}\oint_{\Gamma}\left( \frac{\partial U_r^-}{\partial \mathbf{n}}(\mathbf{r}) u_c(\mathbf{r}) - U_r^-(\mathbf{r}) \frac{\partial u_c}{\partial \mathbf{n}}(\mathbf{r}) \right)d\mathbf{r} \right|, \\
 &\leq \sum_i^{N_p} \left|  \frac{\epsilon_{m}}{2} \oint_{\Gamma_i} \left(  \frac{\partial \phi^-}{\partial \mathbf{n}} (\mathbf{r})u_c(\mathbf{r}) - \phi^-(\mathbf{r})\frac{\partial u_c}{\partial \mathbf{n}}(\mathbf{r}) \right)d\mathbf{r} - \frac{\epsilon_{m}}{2}\oint_{\Gamma_i}\left( \frac{\partial U_r^-}{\partial \mathbf{n}}(\mathbf{r}) u_c(\mathbf{r}) - U_r^-(\mathbf{r}) \frac{\partial u_c}{\partial \mathbf{n}}(\mathbf{r}) \right)d\mathbf{r} \right|, \\ 
 &=\sum_i^{N_p} E_u^i,              
\end{aligned}
\label{eq:error_element_II}
 \end{equation}
where $\Gamma_i$ corresponds to panel $i$ in the discretization of $\Gamma$, and $E_\phi^i$ or $E^i_u$ represent the contribution of element $i$ to the error.


\subsubsection*{\sffamily \normalsize Numerical calculation of $\phi$ and $U_r$}
The error estimates in equations \eqref{eq:error_element_I} and \eqref{eq:error_element_II} need the numerical approximation $U_r$ and the adjoint $\phi$, and their normal derivatives on $\Gamma$.
We compute $U_r$ using equations \eqref{eq:bie} and \eqref{eq:u_reac} assuming  piecewise constant boundary elements.
Calculating $\phi$, which is in principle exact, requires more work.

Starting from the mesh where $U_r$ is solved, we subdivide each triangular panel into four sub-triangles by placing new vertices in each edge center. 
We repeat this process iteratively to obtain an arbitrarily finer mesh conserving the shape of $\Gamma$.
Equation \eqref{eq:adjoint} showed us that $\phi$ equals $u$, hence, we can also compute $\phi$ on $\Gamma$ using equation \eqref{eq:bie}, this time, on the finer mesh and using a piecewise linear ansatz. 
Finally, the local errors, $E^i_\phi$ or $E^i_u$, are computed for each panel on the coarse mesh.
Figure \ref{fig:error_estimation_flowchart} is a summary of the algorithm to compute $E^i_\phi$ and $E^i_u$.

\subsubsection*{\sffamily \normalsize The effectivity ratio $\gamma_{eff}$}
The effectivity ratio, $\gamma_{eff}$, is an indicator of the quality of the error estimation in equation \eqref{eq:err_est_from_err_rep}:
\begin{equation}\label{eq:gamma}
\gamma_{eff} = \frac{E}{\Delta G_{solv} - \Delta\widehat{G}_{solv}}.
\end{equation}
Here, $\Delta G_{solv}$ is the exact value of the solvation free energy, and $E$ can be $E_\phi$ or $E_u$. A $\gamma_{eff}$ close to 1 indicates that the error estimate is accurate. Note that 
$\Delta G_{solv}$ is not available for realistic molecular geometries. 
We compute $\Delta G_{solv}$ using Richardson extrapolation with three consecutive uniform mesh surface-conforming refinements (see figure \ref{fig:refinement_flowchart}), where every boundary element was divided into 4 subtriangles. 
Then, each refinement contains four times more elements than the previous mesh.
Using those three meshes, and knowing that error scales with the average area for a piecewise constant BEM, Richardson extrapolation finds an approximation for a infinitely refined mesh~\cite{Roache1998,CooperBardhanBarba2014}.

\section*{\sffamily \Large Local mesh refinement}

Having the contribution of each triangular element to the error, we can decide if an element should be refined. 
In this work, we sort the elements in descending order according to $E^i_\phi$ or $E_u^i$, and refine those that contribute most to the error, all the way until they add up to 10\% of the total error ($|E_{\phi}|$ or $|E_u|$).
Once the high-error triangles are identified, we perform a barycentric refinement into four subtriangles by adding a vertex on the midpoint of each edge. 
If a triangle that should not be refined shares an edge with one that was refined, this triangle is split into two, by adding an edge between the newly created vertex and the one opposite to it.
However, if a triangle that should not be refined shares two edges with triangles that were refined, we use the same barycentric refinement technique on it.
We can perform this procedure iteratively to obtain finer meshes.
Figure \ref{fig:refinement} is a sketch of the local mesh refinement.

The {\it flat} mesh refinement described in figure \ref{fig:refinement} preserves the geometry of the original mesh.
However, the molecular surface is smooth, and the newly created vertices should adapt to the molecular geometry to represent it more accurately.
To do so, we extended the procedure from figure \ref{fig:refinement} to {\it conform} to the molecular geometry by using an highly refined mesh in the background, and rather than adding the vertex in the edge midpoint, we add the closest vertex of the background mesh.
This way, we make sure the vertex is on the molecular surface.
This process, however, may generate elongated elements, which affect the matrix condition and quality of the solution. 
Then, we use the {\it ImproveSurfMesh} script from GAMer~\cite{YuETal2008}, which improves the mesh of the molecular surface preserving its geometry. 
We call this extended scheme {\it surface-conforming refinement}, and is summarized in figure \ref{fig:refinement_flowchart}.




\section*{\sffamily \Large RESULTS AND DISCUSSION}
In this section, we start by assessing the accuracy of the estimators $E_\phi$ and $E_u$ on methanol, to then use them in an adaptive mesh refinement scheme on a spherical molecule. 
We then apply this method on realistic molecular geometries, in particular, methanol and arginine. First, we use methanol to find optimal meshes for the calculation of $\phi$, and then adaptively refine a mesh on arginine to analyze the efectiveness of the algorithm.

\section*{\sffamily \Large Accuracy of the estimators $E_\phi$ and $E_u$.}

We first assessed the performance of the error estimators $E_\phi$ and $E_u$ using the effectivity ratio $\gamma_{eff}$. 
The results in Table \ref{tab:Results_on_gamma_eff} show $\gamma_{eff}$ for methanol using a mesh with 0.5 and 1 elements per \AA$^2$.
The adjoint $\phi$ was computed on a finer mesh with $N_\phi$ elements, obtained by flat-refining all elements of the original mesh into four subtriangles recursively (see figure \ref{fig:refinement_flowchart}).

We expect that a good error estimator would yield $\gamma_{eff}\approx 1$. 
In table \ref{tab:Results_on_gamma_eff}, we see that as $\phi$ is computed accurately (higher $N_\phi$), $\gamma_{eff}$ does approach 1 for $E_u$, but not for $E_\phi$.
This is an indication that the approximation of setting $R_\phi=0$ in theorem \ref{th:Ephi} is not accurate, and the volumetric integral term has a large contribution. 
Regardless, we continued studying the performance of $E_\phi$ as a per-element error indicator for adaptive mesh refinement.

The detailed distributions of the per-element error estimates are presented in figures \ref{fig:methanol_d05} and \ref{fig:methanol_d1}.
It is interesting to note that the error maps for $E_u$ and $E_\phi$ look almost identical, indicating that they estimate the per-element error distribution similarly, and they recognize the same high-error elements, despite a poor $\gamma_{eff}$ for $E_\phi$.

\section*{\sffamily \Large Mesh refinement on a spherical molecule}
Even though a spherical cavity with internal charges may not be relevant physically, this test becomes useful to understand how the error is distributed on the mesh, as we can control the charge locations, and its effect on the solvation energy, since it has an analytical solution.
Here, we analyzed the error distribution for two configurations, according to figure \ref{fig:Sphere_drawings}: a single off-centered point charge located at half the radius, and a charge and dipole (created by two point charges) placed opposite to each other, $0.62$ into the radius (figure.
For each case, we performed 20 recursive iterations of the surface-conforming adaptive mesh refinement scheme in figure \ref{fig:refinement_flowchart}.
The adjoint $\phi$ was computed on a uniformly flat refined mesh, where every triangle was divided once into four subtriangles.

Figures \ref{fig:dG_offcenter} and \ref{fig:dG_charge-dipole} show $\Delta\widehat{G}_{solv}$ and error with respect to an analytical solution ($\Delta G_{solv}=-52.462648$ kcal/mol and $\Delta G_{solv}=-65.467255$ kcal/mol for the Off-centered and Charge-dipole respectively)~\cite{Kirkwood1934}, computed with the meshes generated from our adaptive mesh refinement scheme. 
The red crosses and black triangles correspond to results with meshes generated using $E_u$ and $E_\phi$, respectively, and they have similar behavior approaching the analytical solution (black segmented line).
This is an indication that in the context of solvation energy, both estimators are equivalent.
The grey dotted lines are the results using a uniformly refined mesh (all elements refined with a surface-conforming method). 
The fact that the red crosses and black triangles are consistently closer to the analytical solution compared to the grey line is evidence of the effectivity of the adaptive mesh refinement.
This is further supported by figures \ref{fig:dG_offcenter} and \ref{fig:dG_charge-dipole}, that evidence a linear relation between the error and number of elements, which is the expected behavior of a piecewise constant BEM~\cite{CooperBardhanBarba2014}.

Figure \ref{table:Sphere_drawings} shows the initial mesh, and the resulting meshes after performing the adaptive refinement scheme iteratively 10, 15, and 20 times, with the corresponding per-element error estimation.
As one could expect, it is evident that refinement is more intense closer to the charge (top of the sphere), and that the bottom of the sphere is more refined for the charge-dipole configuration, compared to the offcentered charge. 
Moreover, the resulting mesh of the charge-dipole case has more elements near the charge than near the dipole, showing that the charge has a higher influence in the error.
Comparing the error distribution on the final mesh between $E_\phi$ and $E_u$, we see that the error is lower and more homogeneously distributed for $E_\phi$.
This is a surprising fact, considering $E_\phi$ has a worse $\gamma_{eff}$ compared to $E_u$.

\section*{\sffamily \Large Adaptive mesh refinement for realistic molecular geometries}

In this section, we study the behavior of our adaptive mesh refinement method on molecular geometries, in particular, for methanol and arginine.
These results are useful to determine the impact of this scheme in real applications.

	\subsection*{\sffamily \large Accuracy of the adjoint}
	
    As as first study of our mesh refinement scheme for realistic molecular geometries, we look at the influence of the accuracy in the calculation of $\phi$ on the resulting meshes.
    Our aim is to determine how fine of a mesh for $\phi$ is required for it to be effective in an adaptive refinement setting.
    Moreover, the calculation of $\phi$ represents the most time consuming part of the algorithm, which may be mitigated by using a coarse mesh.
    Starting from meshes of the molecular surface of methanol with 0.5 and 1 elements per \AA$^2$, we applied the adaptive mesh refinement technique iteratively 6 times, using the surface-conforming scheme.
    We computed $\phi$ in two ways: (a) on the same mesh as $U_r$ ({\it coarse}), and (b) on a mesh with 64 times more elements ({\it fine}), obtained by flat refining all triangles recursively 3 times.
   Figure \ref{fig:methanol_energy} shows $\Delta\widehat{G}_{solv}$ and the error with respect to a Richardson extrapolated value, where the red and black lines use the coarse and fine meshes for $\phi$, respectively.
   The results with the fine and coarse meshes are very similar, which is evidence that the accuracy of the adjoint has a weak effect on the adaptive mesh refinement, and computing $\phi$ on the same mesh as $U_r$ is sufficient. 

    We also present more fine-grained comparisons in figures \ref{table:methanol_adaptive_error} and \ref{table:methanol_adaptive_potential}, where the element-wise error ($E_u$) and electrostatic potential are plotted, respectively, on the original mesh and after 2, 4, and 6 iterations of the adaptive mesh refinement.
    We did not include the equivalent result to figure \ref{table:methanol_adaptive_error} for $E_\phi$ to avoid redundancy, as it was very similar. 
    From these plots, we can see that the resulting meshes with a coarse and fine calculation of $\phi$ are different, even though they yield similar $\Delta \widehat{G}_{solv}$. 
    This happens because even with a poor approximation of $\phi$, the error estimate $E_u$ is capable of detecting high error elements appropriately.
    Moreover, we can see from the results in figure \ref{table:methanol_adaptive_potential} that regions with high electrostatic potential coincide with large error panels (figure \ref{table:methanol_adaptive_error}), and our adaptive mesh refinement technique adds more elements in that area.

	\subsection*{\sffamily \large Mesh refinement on larger structures}

    To estimate the error in larger molecules using a very fine mesh for $\phi$ would be time consuming, however, the results in figure \ref{fig:methanol_energy} indicate that we can compute it on the same mesh as $U_r$.
    Here, we use this fact to perform an adaptive mesh refinement on arginine, starting from meshes with 0.5, 1, and 2 elements per \AA$^2$, aiming towards finding an optimal mesh to compute the solvation energy.

    Figure \ref{fig:arginine_results} shows the convergence of the solvation energy as we apply the adaptive mesh refinement iteratively, using $E_\phi$ and $E_u$.
    We can see that both error estimates generate meshes that are approaching an exact value (obtained with Richardson extrapolation), however, $E_u$ slightly outperforms $E_\phi$.
    The effectiveness of the adaptive mesh refinement becomes evident as, for example, using a mesh refined adaptively from 0.5 elements per \AA$^2$ can reach an error that is lower than the calculation with 2 elements per \AA$^2$, with near half the number of elements.
    
    We can find further evidence of the power of the mesh refinement technique in figure \ref{fig:Arginine-plots}.
    Those plots show the relative error in the $x$ axis and time to solution in the $y$ axis for runs performed with each mesh of the adaptive mesh refinement process, using $E_u$ (red) and $E_\phi$ (black). 
    The size of the markers correspond to the number of elements, and their color to the number of GMRES iterations. 
    Note that the large number of iterations is due to a tight GMRES tolerance (10$^{-8}$), and the fact we use an integral formulation that yields an ill-conditioned matrix, without a preconditioner.
    We can see that as we refine the mesh, the number of elements grows slightly (symbols' sizes remain similar), however, the error drops significantly.
    For example, performing six iterations of the adaptive mesh refinement scheme on the mesh with 0.5 elements per \AA$^2$, the element count only increases from 282 to 328 ($\sim$16\%), whereas the error decreases one order of magnitude. 
    This same behavior is present in the tests starting from 1 and 2 elements per \AA$^2$.
    While the number of elements only increases mildly with adaptive refinement, the computer time increases somewhat, as shown by the rightmost plot in figure \ref{fig:Arginine-plots}. This is because the computer time in these simulations is dominated by the number of iterations rather than the mesh size (see that the simulations that took longer have consistently  lighter symbols). 
    This computer time would decrease if we used better conditioned integral formulations\cite{JufferETal1991} or preconditioners\cite{AltmanBardhanWhiteTidor09}, which control the iteration count.
    Also, these plots show that the meshes generated with $E_u$ (red line) tend to outperform those generated with $E_\phi$ (black line), as the red line is mostly underneath the black line. 
%
%

	%



\section*{\sffamily \Large CONCLUSIONS}
This work presents two adjoint-based goal-oriented error estimates, $E_\phi$ and $E_u$, for the Poisson-Boltzmann equation with BEM, where the quantity of interest is the solvation energy of a solute molecule.
These estimates are written in such a way that we can compute the contribution of each discretization element to the total error, which is useful to detect high-error panels of the mesh.
We saw that even though $E_u$ predicts the error better than $E_\phi$ ($\gamma_{eff}$ closer to 1), both estimates detect the same high-error areas of the mesh.

We used these per-element error indicators to build an adaptive mesh refinement technique, which we tested on a spherical molecule, methanol, and arginine.
Both error estimators had a similar performance in the adaptive mesh refinement, and we found that errors tend to concentrate in areas with high electrostatic potential.
Also, we realized that the accuracy in the calculation of the adjoint had a weak effect on the resulting mesh from the adaptive refinement, and we did not require a finer mesh to resolve it.
We showcase the power of the adaptive mesh refinement in the results for arginine, where the error dropped a factor of 10, by increasing the number of elements only 16\%.

As future work, we plan to use the adaptive mesh refinement technique to generate optimal meshes automatically, aiming towards computing the solvation energy in large-scale molecular simulations efficiently~\cite{MartinezETal2019}.


\section*{\sffamily \large SUPPORTING INFORMATION}
All data and code to reproduce the results of this paper can be obtained from the Github repository \url{https://github.com/RammVI/Error-Formulation}.

\subsection*{\sffamily \large ACKNOWLEDGMENTS}
The authors thank Michael Karkulik (USM) and Ignacio Muga (PUCV) for very helpful conversations. C. Cooper and V. Ramm acknowledge the financial support from ANID (Chile) through FONDECYT Iniciaci\'on 11160768. C. Cooper is also funded by ANID through Basal FB0821.  J. Chaudhry's work is supported by the NSF-DMS (USA) 1720402.


\section*{\sffamily \Large APPENDIX}
\subsection*{\sffamily \large Proof of theorem \ref{th:Ephi}}
\begin{proof}
Using \eqref{eq:adjoint},
\begin{align}
\frac{1}{2}\int_\Omega e_r(\mathbf{r}) \psi(\mathbf{r}) d\mathbf{r} &= \frac{1}{2}\int_\Omega (u_r(\mathbf{r}) - U_r(\mathbf{r})) (-\nabla\cdot \left(\epsilon(\mathbf{r})\nabla \phi(\mathbf{r}) \right) + \overline{\kappa}^2(\mathbf{r}) \phi(\mathbf{r}) ) d\mathbf{r}, \\
& = \frac{1}{2}\int_\Omega u_r(\mathbf{r}) (-\nabla\cdot \left(\epsilon(\mathbf{r})\nabla \phi(\mathbf{r}) \right) + \overline{\kappa}^2(\mathbf{r}) \phi(\mathbf{r}) )d\mathbf{r}, \nonumber\\
& - \frac{1}{2}\int_\Omega U_r(\mathbf{r}) (-\nabla\cdot \left(\epsilon(\mathbf{r})\nabla \phi(\mathbf{r}) \right) + \overline{\kappa}^2(\mathbf{r}) \phi(\mathbf{r}) )d\mathbf{r},
=I - II \label{eq:err_I_II}
\end{align}
where
\begin{align}
I =  \frac{1}{2}\int_\Omega u_r(\mathbf{r}) (-\nabla\cdot \left(\epsilon(\mathbf{r})\nabla \phi(\mathbf{r}) \right)+ \overline{\kappa}^2(\mathbf{r}) \phi(\mathbf{r}) )d\mathbf{r} \label{eq:I}\\
II = \frac{1}{2}\int_\Omega U_r(\mathbf{r}) (-\nabla\cdot \left(\epsilon(\mathbf{r})\nabla \phi(\mathbf{r}) \right)+ \overline{\kappa}^2(\mathbf{r}) \phi(\mathbf{r}) )d\mathbf{r}.
\end{align} \label{eq:II}
First we consider the term $I$.

\begin{align}\label{eq:Gsolv_ex2}
 I &= \frac{1}{2} \int_{\Omega} \left[-\nabla\cdot \left(\epsilon(\mathbf{r})\nabla \phi(\mathbf{r}) \right) + \overline{\kappa}^2 \phi(\mathbf{r})\right]u_r(\mathbf{r})d\mathbf{r},\nonumber\\
&= \frac{1}{2} \left[\int_{\Omega_m} -\epsilon_m\nabla^2 \phi(\mathbf{r})u_r(\mathbf{r})d\mathbf{r} + \int_{\Omega_w} -\epsilon_s\nabla^2 \phi(\mathbf{r}) u_r(\mathbf{r})d\mathbf{r} + \int_{\Omega_w}\overline{\kappa}^2 \phi(\mathbf{r})u_r(\mathbf{r})d\mathbf{r}\right].
\end{align}
We now use separation of variables to write
\begin{align}
I =\frac{1}{2}&\left[ -\int_{\Omega_m} \left\lbrace \epsilon_m\nabla\cdot(\nabla\phi(\mathbf{r}) u_r(\mathbf{r})) - \epsilon_m\nabla\cdot(\phi(\mathbf{r})\nabla u_r(\mathbf{r})) + \phi(\mathbf{r})(\epsilon_m\nabla^2u_r(\mathbf{r})) \right\rbrace d\mathbf{r}\right.\nonumber\\
 &- \int_{\Omega_w} \left\lbrace\epsilon_w\nabla\cdot(\nabla\phi(\mathbf{r}) u_r(\mathbf{r})) - \epsilon_w\nabla\cdot(\phi(\mathbf{r})\nabla u_r(\mathbf{r})) + \phi(\mathbf{r})(\epsilon_w\nabla^2u_r(\mathbf{r})) \right\rbrace d\mathbf{r}\nonumber\\
&+\left. \int_{\Omega_w}\phi(\mathbf{r})\overline{\kappa}^2u_r(\mathbf{r})d\mathbf{r}\right]
\end{align}
and then use the divergence theorem to obtain
\begin{align}
I=\frac{1}{2}&\left[-\oint_\Gamma \left( \epsilon_m \frac{\partial\phi^-}{\partial \mathbf{n}}(\mathbf{r})u_r^-(\mathbf{r}) - \phi^-(\mathbf{r})\epsilon_m\frac{\partial u_r^-}{\partial\mathbf{n}}(\mathbf{r}) - \epsilon_w\frac{\partial\phi^+}{\partial\mathbf{n}}(\mathbf{r})u_r^+(\mathbf{r}) + \phi^+(\mathbf{r})\epsilon_w\frac{\partial u_r^+}{\partial\mathbf{n}}(\mathbf{r}) \right) d\mathbf{r}\right.\nonumber\\
&+\left.\int_{\Omega_w+\Omega_m}\phi(\mathbf{r})\left(-\nabla\cdot(\epsilon(\mathbf{r})\nabla u_r(\mathbf{r})) + \overline{\kappa}^2(\mathbf{r}) u_r(\mathbf{r})\right) d\mathbf{r}\right],
\end{align}
recalling that $\mathbf{n}$ points out of $\Omega_m$.
Then, we can apply the interface conditions on $\phi$ (equation \eqref{eq:phi_jump}) and $u_r$ (equation \eqref{eq:ur_jump}) to write

\begin{align}\label{eq:Gsolv_ex_decomp}
I &= \frac{\epsilon_{m}}{2} \oint_{\Gamma} \left(  \frac{\partial \phi^-}{\partial \mathbf{n}} (\mathbf{r})u_c(\mathbf{r}) - \phi^-(\mathbf{r})\frac{\partial u_c}{\partial \mathbf{n}}(\mathbf{r}) \right)d\mathbf{r} \nonumber\\
        &+ \frac{1}{2} \int_{\Omega_w+\Omega_m} \phi(\mathbf{r}) (-\nabla\cdot(\epsilon(\mathbf{r}) \nabla u_r(\mathbf{r}))+ \overline{\kappa}^{2}(\mathbf{r})u_r(\mathbf{r})) d\mathbf{r}.
\end{align}

Now, considering the term $II$ in equation \eqref{eq:II},

\begin{align}\label{eq:Gsolv_ap}
II &= \frac{1}{2}\int_{\Omega} \left[ \sum_{k}^{N_{q}} q_{k}\delta (\mathbf{r}-\mathbf{r}_{k}) \right] U_{r}(\mathbf{r})d\mathbf{r}, \nonumber\\
&= \frac{1}{2}\int_{\Omega} \left[-\nabla\cdot\left(\epsilon(\mathbf{r})\nabla \phi(\mathbf{r})\right) + \overline{\kappa}^2 \phi(\mathbf{r})\right]U_r(\mathbf{r})d\mathbf{r},\nonumber\\
&= \frac{1}{2}\int_{\Omega_m} -\epsilon_m\nabla^2 \phi(\mathbf{r})U_{r}(\mathbf{r})d\mathbf{r} + \frac{1}{2}\int_{\Omega_w} \left[-\epsilon_s\nabla^2 \phi(\mathbf{r}) + \overline{\kappa}^2 \phi(\mathbf{r})\right]U_r(\mathbf{r})d\mathbf{r}.
\end{align}
Considering that there are no point charges in the outer region, equation \eqref{eq:adjoint} equals zero in $\Omega_w$, and the last integral equation \eqref{eq:Gsolv_ap} cancels out.
Using separation of variables and the divergence theorem, we can write
\begin{align}\label{eq:Gsolv_ap2}
II & = -\frac{\epsilon_{m}}{2}\int_{\Omega_{m}} \nabla^{2} \phi(\mathbf{r}) U_{r}(\mathbf{r})d\mathbf{r}, \nonumber\\
& = -\frac{\epsilon_{m}}{2} \int_{\Omega_{m}} \left[ \nabla \cdot (\nabla \phi(\mathbf{r}) U_{r}(\mathbf{r})) - \nabla\cdot ( \phi(\mathbf{r}) \nabla U_r(\mathbf{r})) + \phi(\mathbf{r}) \nabla^{2}U_r(\mathbf{r}) \right] d\mathbf{r}, \nonumber \\
& = -\frac{\epsilon_{m}}{2} \left[ \oint_{\Gamma}\left( \frac{\partial \phi^-}{\partial \mathbf{n}}(\mathbf{r}) U_r^-(\mathbf{r}) - \phi^-(\mathbf{r}) \frac{\partial U_r^-}{\partial n}(\mathbf{r}) \right)d\mathbf{r} + \int_{\Omega_{m}} \phi(\mathbf{r}) \nabla^{2}U_r(\mathbf{r})d\mathbf{r} \right].
\end{align}

Combining equations \eqref{eq:err_I_II}, \eqref{eq:Gsolv_ex_decomp} and \eqref{eq:Gsolv_ap2} completes the proof.

\end{proof}

\subsection*{\sffamily \large Proof of theorem \ref{th:Eu}}

\begin{proof}
We start from equation \eqref{eq:error_approx}, to decompose the error as
\begin{equation}\label{eq:err_I_II_b}
\frac{1}{2}\int_\Omega (u_r(\mathbf{r}) - U_r(\mathbf{r})) \psi d\mathbf{r} = \frac{1}{2}\int_\Omega u_r(\mathbf{r}) \psi(\mathbf{r}) d\mathbf{r} - \frac{1}{2}\int_\Omega U_r(\mathbf{r}) \psi(\mathbf{r}) d\mathbf{r} = I - II,
\end{equation}
where
\begin{align}
I &=  \frac{1}{2}\int_\Omega u_r(\mathbf{r}) \psi(\mathbf{r}) d\mathbf{r} = \frac{1}{2} \int_{\Omega} \left[-\nabla\cdot\left(\epsilon(\mathbf{r})\nabla \phi(\mathbf{r})\right) + \overline{\kappa}^2 \phi(\mathbf{r})\right]u_r(\mathbf{r})d\mathbf{r}, \label{eq:I_b}\\
II &= \frac{1}{2}\int_\Omega U_r(\mathbf{r}) \psi(\mathbf{r}) d\mathbf{r}.\label{eq:II_b}
\end{align} 

We already derived an expression for term $I$ in equation \eqref{eq:Gsolv_ex_decomp}, and we will move straight to term $II$.
Using equation \eqref{eq:psi_rho} and considering $\psi = \rho = \sum_{k=1}^{N_q} q_k\delta(|\mathbf{r}-\mathbf{r}_k|)$ in equation \eqref{eq:pbe}, we find
\begin{align}
II &=  \frac{1}{2}\int_\Omega U_r(\mathbf{r})  \sum_{k=1}^{N_q} q_k\delta(|\mathbf{r}-\mathbf{r}_k|) d\mathbf{r} = \frac{1}{2}\sum_{k=1}^{N_q} U_r(\mathbf{r}_k) q_k.
\end{align}
$U_r$ can be further expanded with the numerical approximation of equation \eqref{eq:u_reac} to write
\begin{align}
II &= \frac{1}{2}\sum_{k=1}^{N_q} q_k\left[-\oint_\Gamma U^-(\mathbf{r})\frac{\partial}{\partial\mathbf{n}}\left(\frac{1}{4\pi|\mathbf{r}_k-\mathbf{r}|}\right)d\mathbf{r} + \oint_\Gamma \frac{\partial U^-(\mathbf{r})}{\partial\mathbf{n}}\frac{1}{4\pi|\mathbf{r}_k-\mathbf{r}|}d\mathbf{r}\right].
\end{align}
Considering the sum and integral are linear operators, we can swap them to bring $\sum_{k=1}^{N_k} q_k$ into the integral, as 
\begin{align}
II &= \frac{1}{2}\Bigg[-\oint_\Gamma U^-(\mathbf{r})\frac{\partial}{\partial\mathbf{n}}\underbrace{\left(\sum_{k=1}^{N_q} q_k\frac{1}{4\pi|\mathbf{r}_k-\mathbf{r}|}\right)}_{=\epsilon_m u_c}d\mathbf{r} + \oint_\Gamma \frac{\partial U^-(\mathbf{r})}{\partial\mathbf{n}}\underbrace{\sum_{k=1}^{N_q} q_k\frac{1}{4\pi|\mathbf{r}_k-\mathbf{r}'|}}_{=\epsilon_m u_c}d\mathbf{r}\Bigg],
\end{align}
where, as indicated in the previous equation, we have recovered $u_c$ from equation \eqref{eq:u_c}.
Then, we get
\begin{equation}\label{eq:II_final}
II = \frac{\epsilon_m}{2}\left[-\oint_\Gamma U^-(\mathbf{r})\frac{\partial u_c}{\partial\mathbf{n}}(\mathbf{r})d\mathbf{r} + \oint_\Gamma \frac{\partial U^-(\mathbf{r})}{\partial\mathbf{n}}u_c(\mathbf{r})d\mathbf{r}\right].
\end{equation}

Combining equations \eqref{eq:err_I_II_b}, \eqref{eq:Gsolv_ex_decomp} and \eqref{eq:II_final} completes the proof.

\end{proof}


\clearpage




\clearpage
\begin{figure}
\caption{ \label{fig:Solvation-impl}
(Left) Solute molecule surrounded by explicit solvent molecules. (Right) Representation of the implicit solvent model.
The hatched area corresponds to the unbounded implicit solvent ($\Omega_w$, with $\epsilon_w$ and $\kappa$), which has a cavity containing the solute molecule ($\Omega_m$, with $\epsilon_m$).}
\end{figure}

\begin{figure}
    \caption{\label{fig:ur} Thermodynamic cycle of molecular solvation. Initially (pane I), the solute is isolated in vacuum (only point charges) and the electrostatic potential in the solvent is zero. Then, we place the charges inside the solvent, generating a reaction potential (pane II).}
\end{figure}

\begin{figure}
	\caption{\label{fig:error_estimation_flowchart} Process of creating the mesh, solving and calculating the error.}
\end{figure}

\begin{figure}
\caption{\label{fig:refinement}Example of a local mesh refinement procedure. In the left pane, the two elements marked in blue are identified to have a high $E^i_{\phi}$ or $E_u^i$, and will be subdivided into four triangles. In the middle pane we identify neighbor elements: those marked in light blue only share one edge with the high-error triangles, and they will be divided in two, however, there is one triangle in the middle that shares two edges with refined triangles, then, we mark it in blue, and it will be refined in four triangles. This generates a new neighbor triangle to be divided in two, marked in light blue below the blue triangles. Finally, the right pane shows the resulting refined mesh. We can perform this procedure iteratively for higher mesh refinements.}
\end{figure}

\begin{figure}
\caption{\label{fig:refinement_flowchart}Summary of the local mesh refinement procedure.}
\end{figure}

\begin{figure}
	\caption{\label{fig:methanol_d05}Per-element error estimation for methanol with a 0.5 el/\AA$^{2}$ mesh.}
\end{figure}

\begin{figure}
	\caption{\label{fig:methanol_d1}Per-element error estimation for methanol with a 1 el/\AA$^{2}$ mesh.}
\end{figure}

\begin{figure}
    \caption{\label{fig:Sphere_drawings}Offcenter charge (top) and charge-dipole (bottom) distributions for the spherical cavity}
\end{figure}

\begin{figure}
	\caption{\label{fig:dG_offcenter}Solvation energy (left) and error (right) for the off-centered charge configuration. Meshes for the results marked with red crosses were obtained with $E_u$ ($\Delta\widehat{G}_{solv}^u$) and the ones marked with black triangles were obtained with $E_\phi$ ($\Delta\widehat{G}_{solv}^\phi$). The grey dotted lines correspond to uniformly refined meshes (all elements with a surface conforming method, $\Delta\widehat{G}_{solv}^{unif}$). The segmented black line is the true solution ($\Delta G_{solv}$), computed analytically. 
	}
    
\end{figure}

\begin{figure}
		\caption{\label{fig:dG_charge-dipole}Solvation energy (left) and error (right) for the charge-dipole configuration. Meshes for the results marked with red crosses were obtained with $E_u$ ($\Delta\widehat{G}_{solv}^u$) and the ones marked with black triangles were obtained with $E_\phi$ ($\Delta\widehat{G}_{solv}^\phi$). The grey dotted lines correspond to uniformly refined meshes (all elements with a surface conforming method, $\Delta\widehat{G}_{solv}^{unif}$). The segmented black line is the true solution ($\Delta G_{solv}$), computed analytically. }
\end{figure}

\begin{figure}
	\caption{Initial and resulting meshes after 10, 15, and 20 adaptive refinements for the sphere cases. Colors correspond to the per-element error estimations $E_{\phi}$ and $E_{u}$.}
    \label{table:Sphere_drawings}
\end{figure}

\begin{figure}
		\caption{\label{fig:methanol_energy}Solvation energy (top) and error (bottom) for methanol with meshes refined adaptively, computing $\phi$ on a coarse (red line) and fine (black line) mesh, starting from a mesh with $0.5[El/$\AA$^{2}]$ (left) and $1.0[El/$\AA$^{2}]$ (right). Solid lines correspond to results using $E_u$, whereas segmented lines used $E_\phi$.}
\end{figure}

\begin{figure}
	\caption{\label{table:methanol_adaptive_error}Per element error estimation ($E_u^i$) using a {\it fine} and {\it coarse} mesh to obtain $\phi$.}
\end{figure}

\begin{figure}
    \caption{\label{table:methanol_adaptive_potential}Total electrostatic potential on the molecular surface.}
\end{figure}
	
\begin{figure}
    \caption{\label{fig:arginine_results}$\Delta \widehat{G}_{solv}$ (left) and error (right) for arginine using meshes generated with the adaptive mesh refinement technique with a surface-conforming scheme. Results in red use $E_u$ ($\Delta\widehat{G}_{solv}^u$) whereas those in black use $E_\phi$ ($\Delta\widehat{G}_{solv}^\phi$). The black segmented line in the left pane corresponds to a Richardson extrapolated value for the energy, which is also used as the reference for the error calculations.}
\end{figure}

\begin{figure}
		\caption{\label{fig:Arginine-plots}Relative error versus time to solution for each mesh using adaptive refinement. Size and color of the markers indicate the number of mesh elements and number of GMRES iterations, respectively. Results following the red line use $E_u$ whereas those following the black line use $E_\phi$.}
\end{figure}

%


\clearpage
\begin{center}
\includegraphics[width=0.7\textwidth, trim= 4cm 5cm 4cm 5cm, clip]{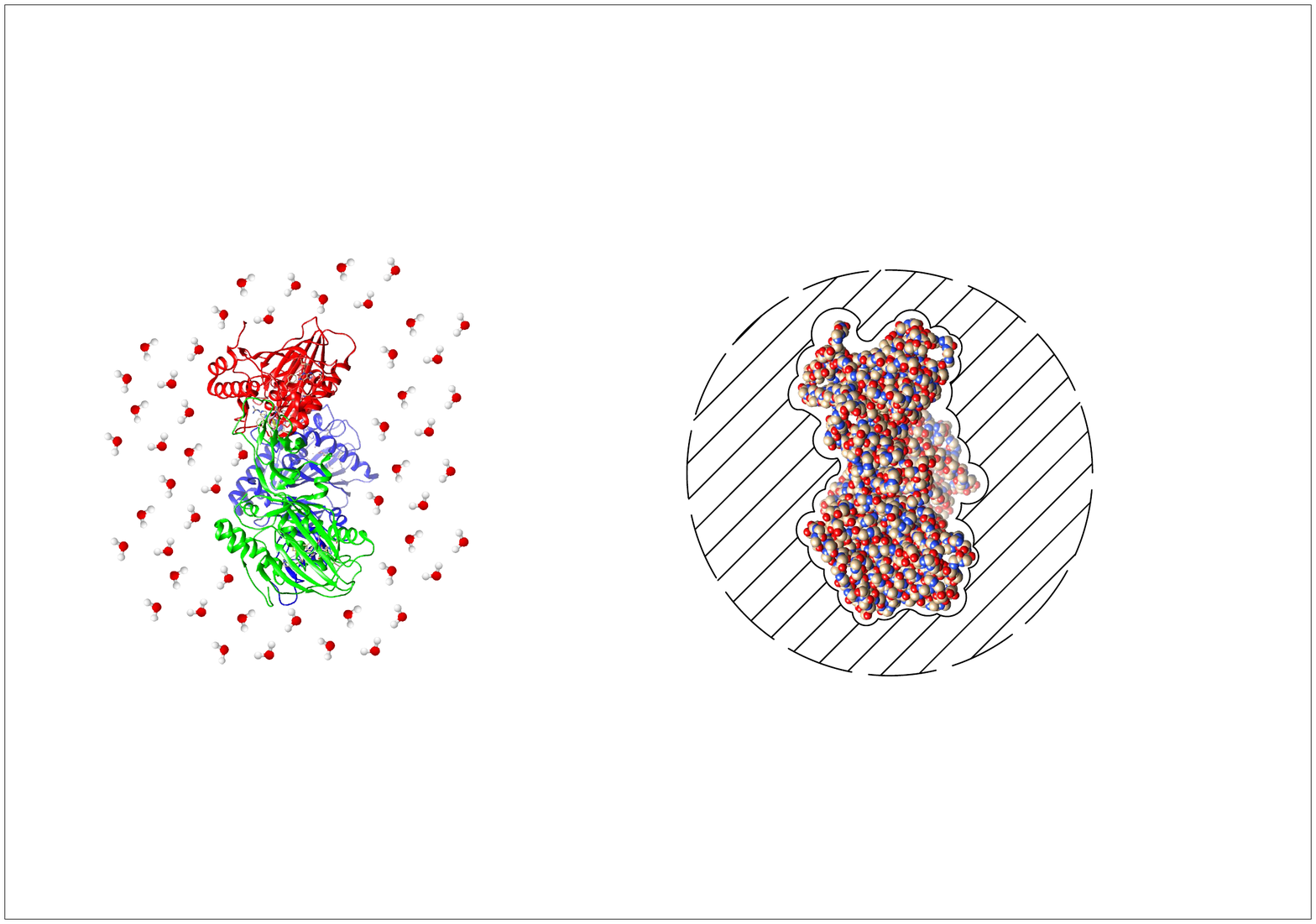}
\end{center}
\vspace{0.25in}
\hspace*{3in}
{\Large
	\begin{minipage}[t]{3in}
		\baselineskip = .5\baselineskip
		Figure 1 \\
		Vicente Ramm, Jehanzeb H. Chaudhry, Christopher D. Cooper \\
		J.\ Comput.\ Chem.
	\end{minipage}
}
\clearpage

\begin{center}
	\begin{overpic}[width=0.6\linewidth, trim= 2cm 5cm 1cm 5cm, clip]{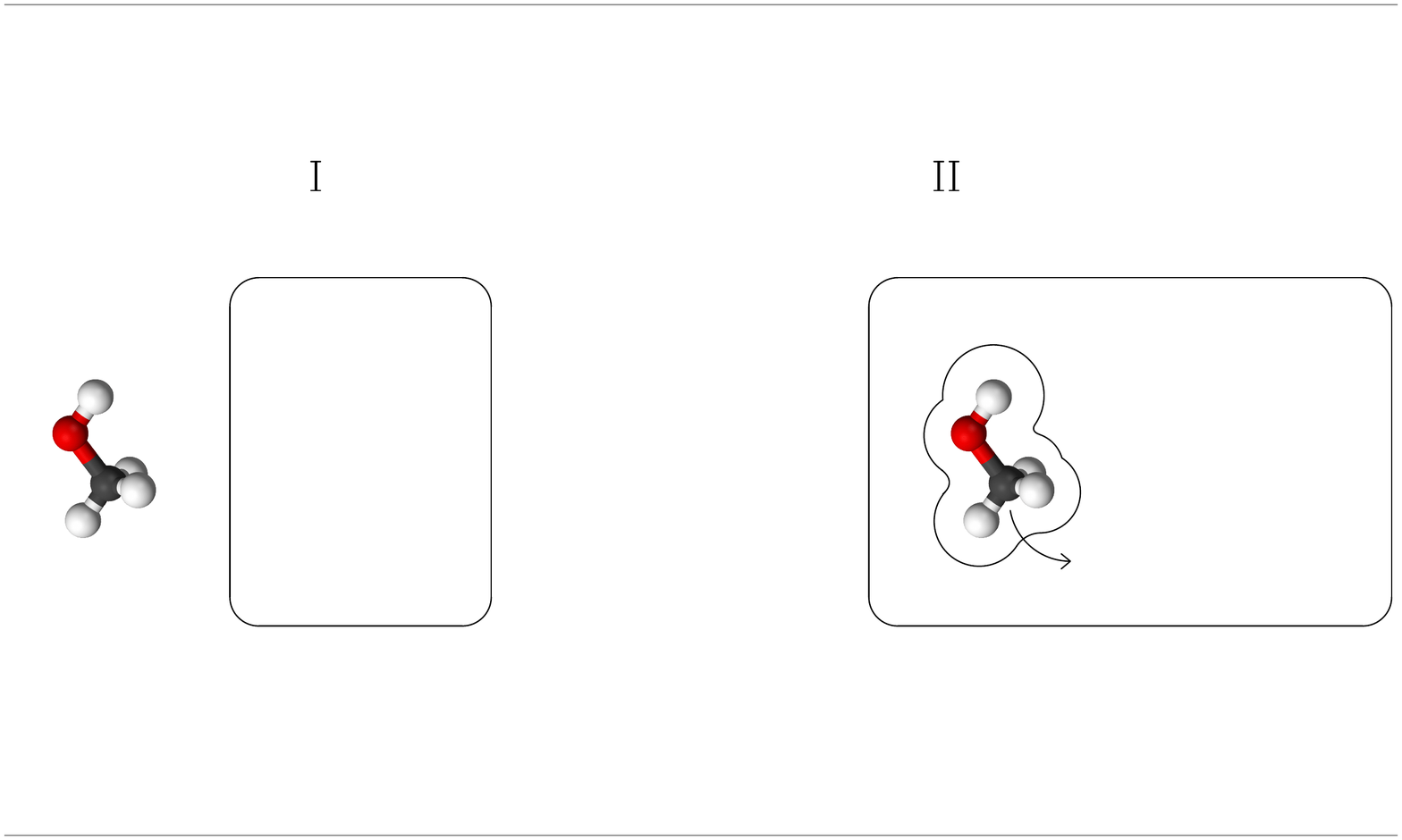} 
		\put (2,27)  {\scriptsize Vacuum} %
		\put (18,27) {\scriptsize Solvent}%
		\put (1,24) {\scriptsize $u=u_{c}$}%
		\put (48,27)  {\scriptsize Vacuum} %
		\put (63,27) {\scriptsize Solvent}%
		\put (75,24) {\scriptsize $u=u_{r}$}%
		\put (76,9) {\scriptsize $u=u_{r}+u_{c}$}%
	\end{overpic}
\end{center}
\vspace{0.25in}
\hspace*{3in}
{\Large
	\begin{minipage}[t]{3in}
		\baselineskip = .5\baselineskip
		Figure 2 \\
		Vicente Ramm, Jehanzeb H. Chaudhry, Christopher D. Cooper \\
		J.\ Comput.\ Chem.
	\end{minipage}
}
\clearpage

\begin{center}
	\begin{tikzpicture}
	[node distance=.4cm ,
	start chain=going below,]
	\node[punktchain, join] (node1) {Generate the mesh}         ;
	\node[punktchain, join] (node2) {Solve $U^{-}$ and $\frac{\partial U^{-}}{\partial n^{-}}$}     ;
	\begin{scope}[start branch=Ephi,]
		\node[punktchain,on chain=going right] (node3) {Calculate $U_{r}^{-}$ by eq. \eqref{eq:u_reac}}; 
	\end{scope}
	\node[punktchain, join ] (node4)[below=1.2cm]{Solve $\phi^{-}$ and $\frac{\partial \phi^{-}}{\partial n^{-}}$ on finer Adjoint mesh}        ;
	\node[punktchain, join ] (node5){Calculate $E_{\phi}^{i}$ or $E^i_u$}            ;
	\draw[dashed,-, thick] (node3.south) |-+(0,-1em)-| (node4.north);
	\draw[dashed,->, thick] (node2)--(node4);
	\end{tikzpicture}
\end{center}
\vspace{0.25in}
\hspace*{3in}
{\Large
	\begin{minipage}[t]{3in}
		\baselineskip = .5\baselineskip
		Figure 3 \\
		Vicente Ramm, Jehanzeb H. Chaudhry, Christopher D. Cooper \\
		J.\ Comput.\ Chem.
	\end{minipage}
}
\clearpage

\begin{center}
\begin{tikzpicture}
    \node[block] (a) {\includegraphics[width=\linewidth]{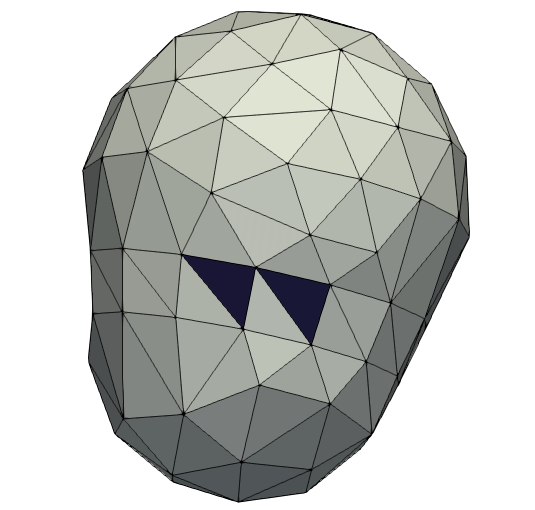}};
    \node[block,right=of a] (b) {
        \includegraphics[width=\linewidth]{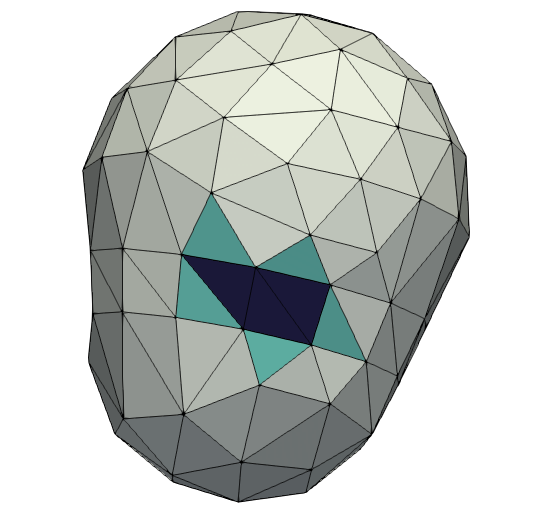}};
    \node[block,right=of b] (c) {
        \includegraphics[width=\linewidth]{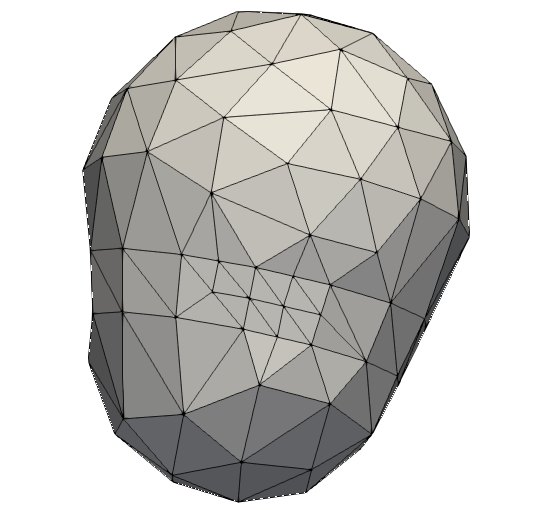}
    };
    \draw[line] (b)-- (a);
    \draw[line] (c)-- (b);
\end{tikzpicture}
\end{center}
\vspace{0.25in}
\hspace*{3in}
{\Large
	\begin{minipage}[t]{3in}
		\baselineskip = .5\baselineskip
		Figure 4 \\
		Vicente Ramm, Jehanzeb H. Chaudhry, Christopher D. Cooper \\
		J.\ Comput.\ Chem.
	\end{minipage}
}
\clearpage

\begin{center}
\begin{tikzpicture}
[node distance=.4cm,
start chain=going below,]
\node[punktchain, join] (node1) {Mark elements}         ;
\node[punktchain, join] {Assign neighbor criteria}     ;
\node[punktchain, join] (node3){Flat refinement}        ;
\node[punktchain, join] {Add closest vertex from background mesh}            ;
\node[punktchain, join] (end){\textit{ImproveSurfMesh}} ;
\draw[tuborg, decoration={brace}] let \p1=(node1.north),
\p2=(node3.south) in ($(3, \y1)$) -- ($(3, \y2)$) node[tubnode]
{Flat refinement}                                        ;
\draw[tuborg, decoration={brace}] let \p1=(node1.north),
\p2=(end.south) in ($(7, \y1)$) -- ($(7, \y2)$) node[tubnode]
{Surface-conforming refinement}                                            ;
\end{tikzpicture}
\end{center}
\vspace{0.25in}
\hspace*{3in}
{\Large
	\begin{minipage}[t]{3in}
		\baselineskip = .5\baselineskip
		Figure 5 \\
		Vicente Ramm, Jehanzeb H. Chaudhry, Christopher D. Cooper \\
		J.\ Comput.\ Chem.
	\end{minipage}
}
\clearpage

\begin{center}
	\begin{tabular}{m{1cm}m{1cm}|m{2.1cm}|m{2.1cm}|m{2.1cm}|m{2.1cm}|m{2.1cm}|m{2.1cm}}   \cline{3-7}
		& & \multicolumn{5}{c|}{$N_{\phi}$} & \\
		 \cline{2-7}
		& \multicolumn{1}{|c|}{View} & \multicolumn{1}{c|}{48} & \multicolumn{1}{c|}{192} & \multicolumn{1}{c|}{768} &  \multicolumn{1}{c|}{3072} & \multicolumn{1}{c|}{12288} &  \\ \cline{1-7}
		\multicolumn{1}{|c}{\multirow{2}{*}{ \begin{minipage}{.9cm} \vspace{36pt}\centering
					$E_{\phi}$
		\end{minipage} }} &
		\multicolumn{1}{|c|}{Front} &
		\imagespacing
		\includegraphics[width=\linewidth]{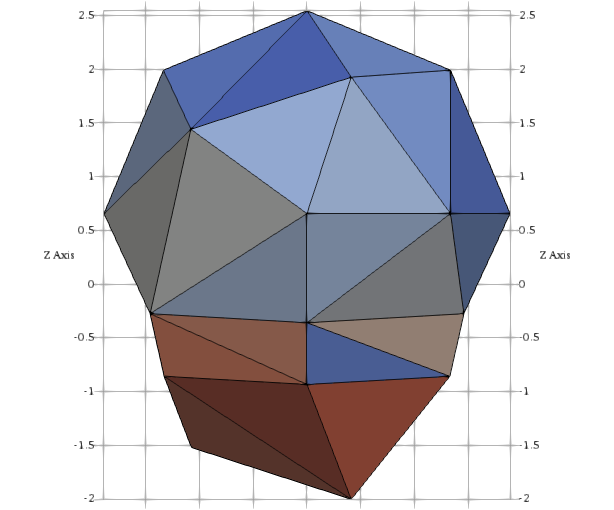} &	
		\imagespacing
		\includegraphics[width=\linewidth]{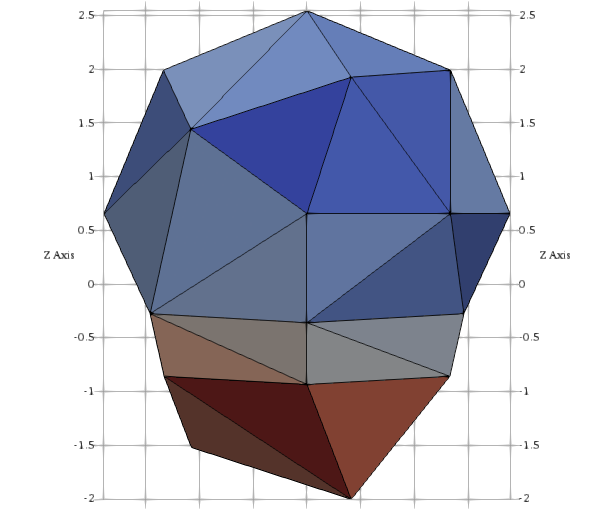} &	
		\imagespacing
		\includegraphics[width=\linewidth]{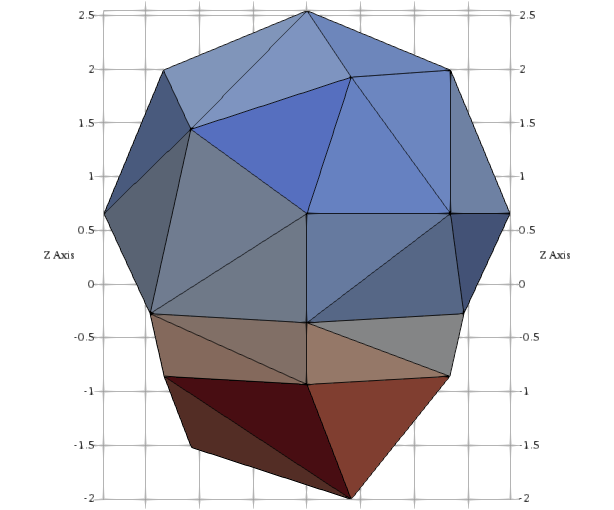} &	
		\imagespacing
		\includegraphics[width=\linewidth]{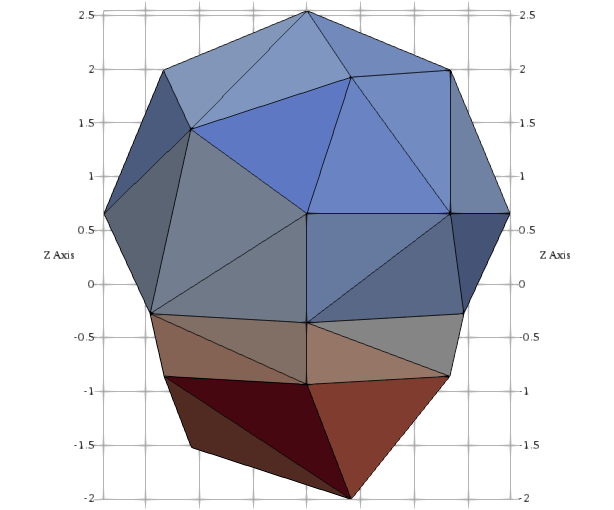} &	
		\imagespacing
		\includegraphics[width=\linewidth]{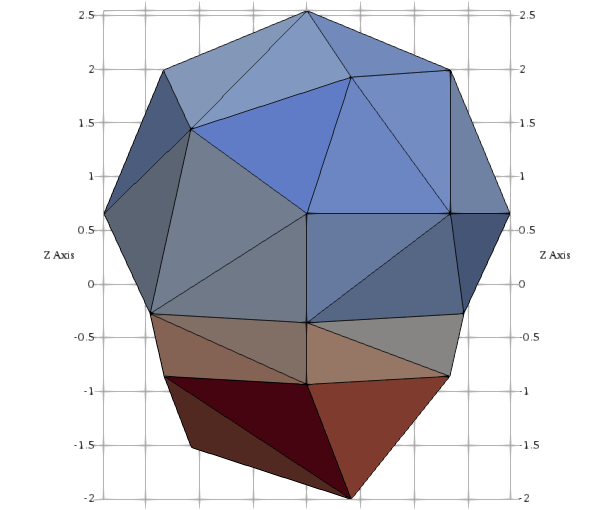} &	 
		\multirow{4}{\linewidth}{\vspace{2pt}
			\begin{minipage}{2.1cm}
				\includegraphics[height=2.4\linewidth]{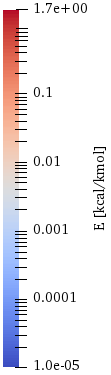}
			\end{minipage} 
		} \\ \cline{2-7}
		\multicolumn{1}{|c}{} & \multicolumn{1}{|c|}{Back} & 
		\imagespacing \includegraphics[width=\linewidth]{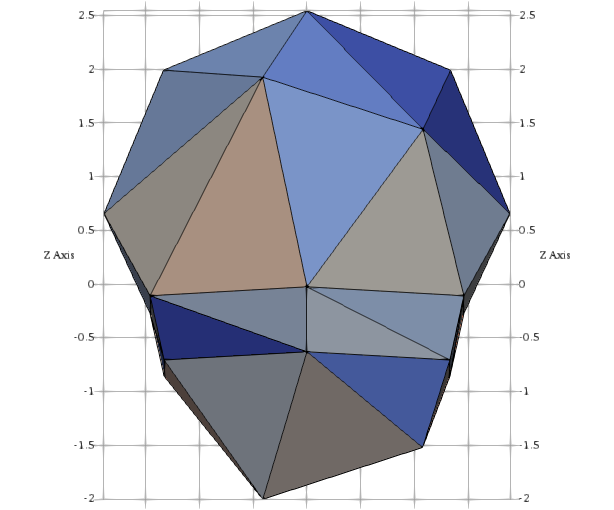} &	
		\imagespacing
		\includegraphics[width=\linewidth]{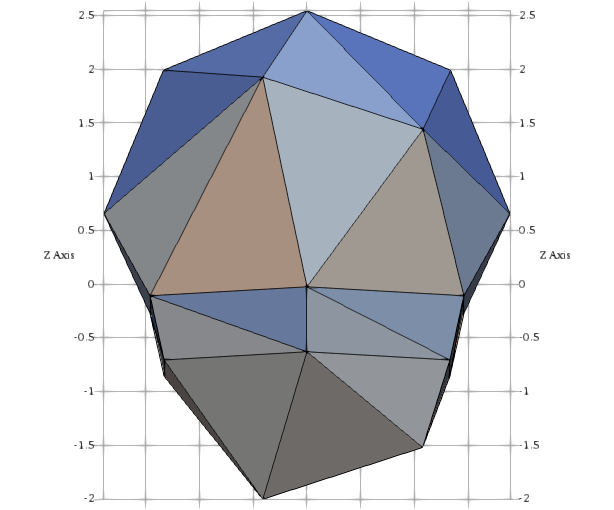} &	
		\imagespacing
		\includegraphics[width=\linewidth]{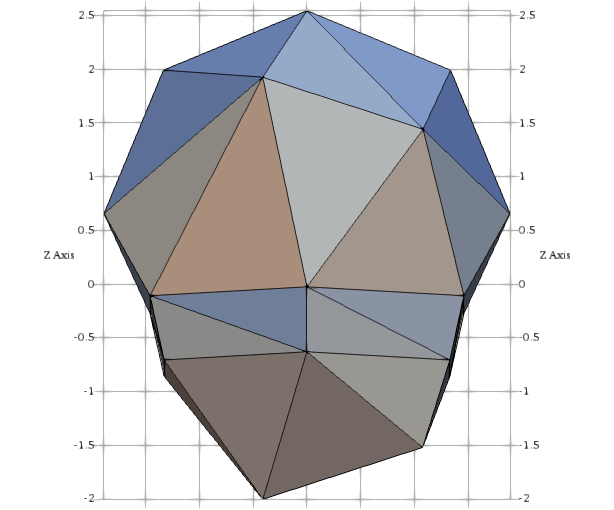} &	
		\imagespacing
		\includegraphics[width=\linewidth]{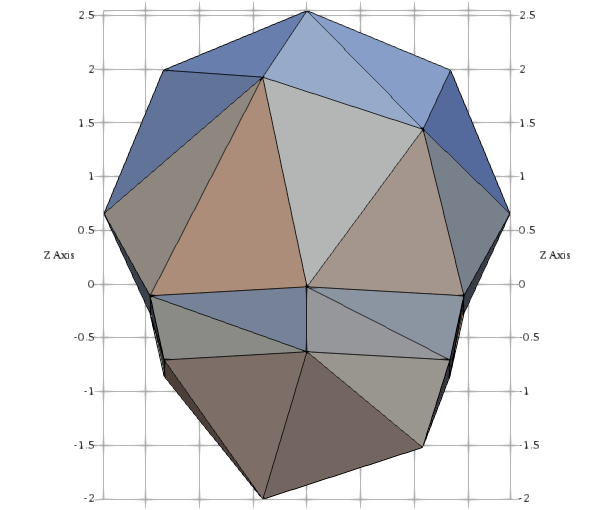} &	
		\imagespacing
		\includegraphics[width=\linewidth]{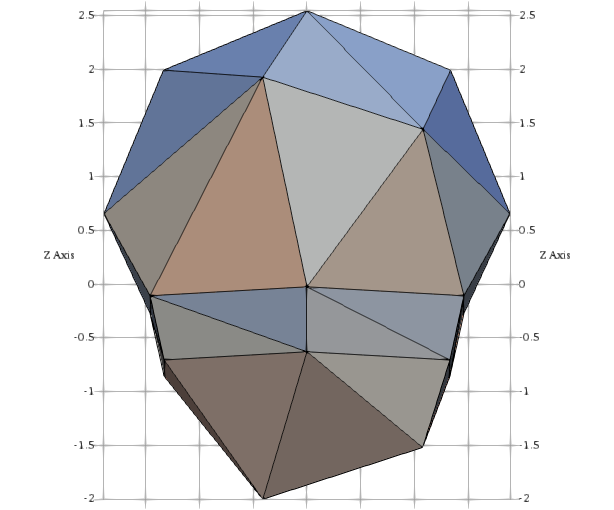} &	\\ \cline{1-7}
		\multicolumn{1}{|c}{\multirow{2}{*}{ \begin{minipage}{.9cm} \vspace{36pt}\centering
					$E_{u}$
		\end{minipage} }} & \multicolumn{1}{|c|}{Front} &
		\imagespacing \includegraphics[width=\linewidth]{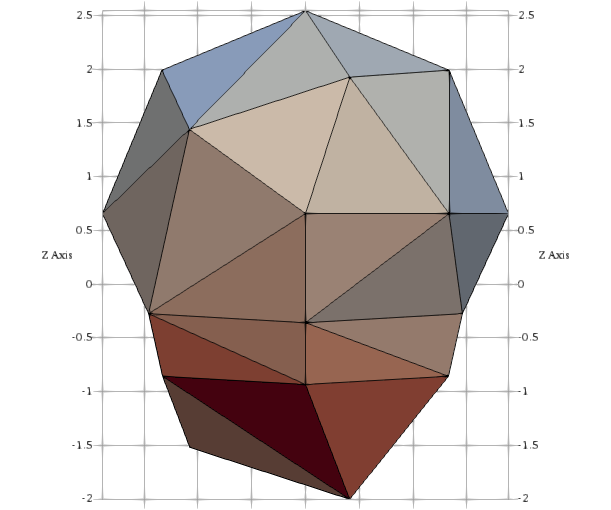}   &
		\imagespacing
		\includegraphics[width=\linewidth]{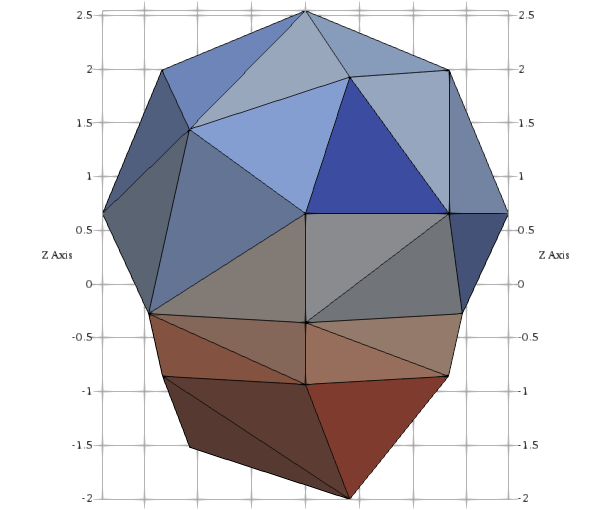}   &
		\imagespacing \includegraphics[width=\linewidth]{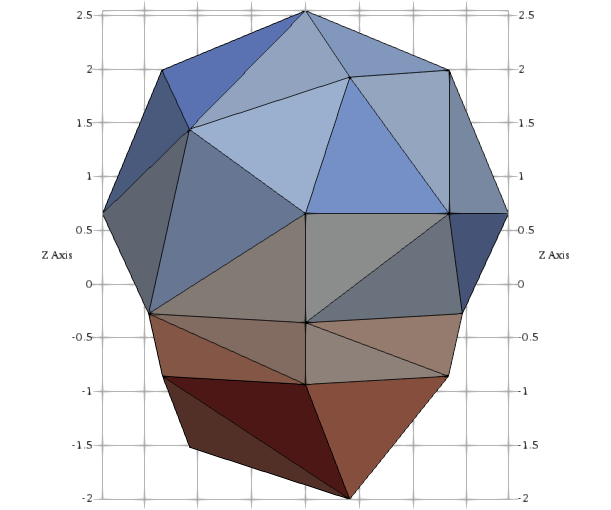}   &
		\imagespacing \includegraphics[width=\linewidth]{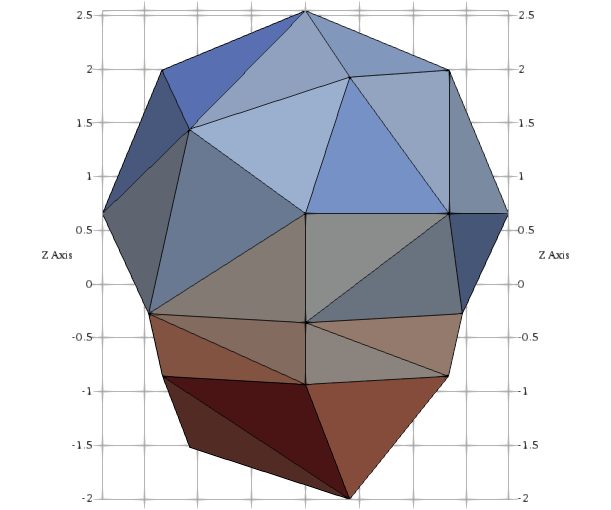}   &
		\imagespacing \includegraphics[width=\linewidth]{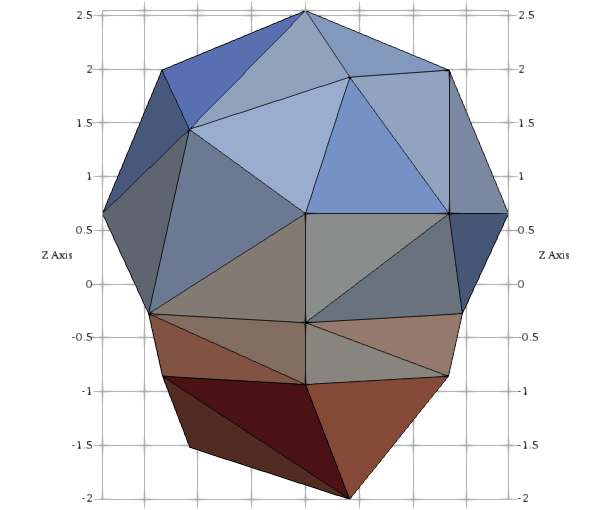} & 
		\\ \cline{2-7}
		\multicolumn{1}{|c}{}& \multicolumn{1}{|c|}{Back}  &
		\imagespacing \includegraphics[width=\linewidth]{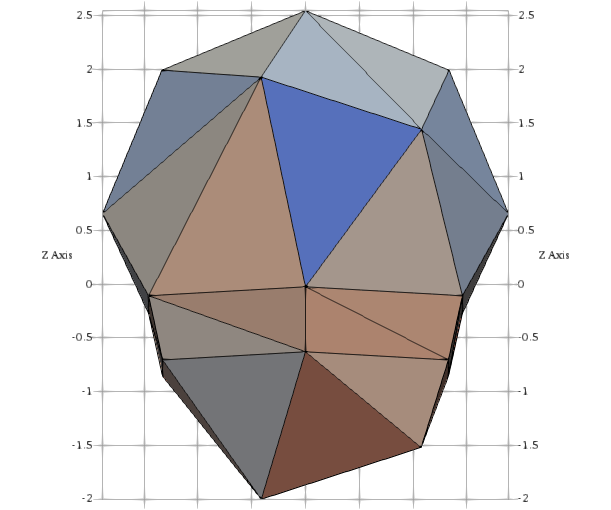}   &
		\imagespacing
		\includegraphics[width=\linewidth]{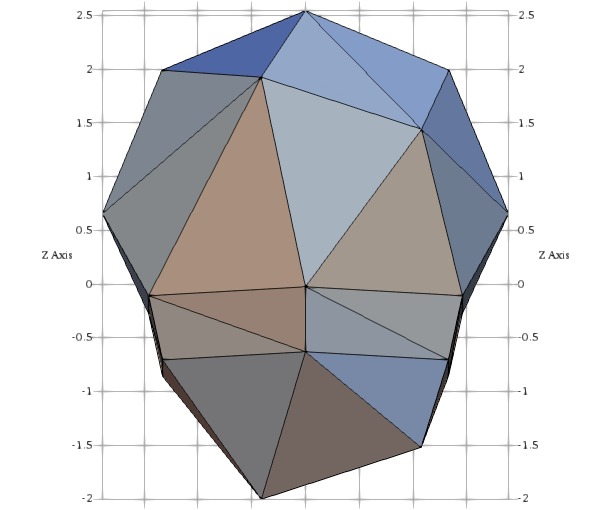}   &
		\imagespacing \includegraphics[width=\linewidth]{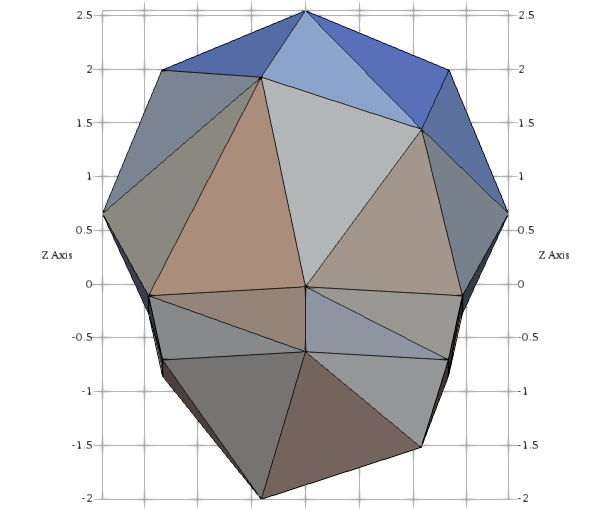}   &
		\imagespacing \includegraphics[width=\linewidth]{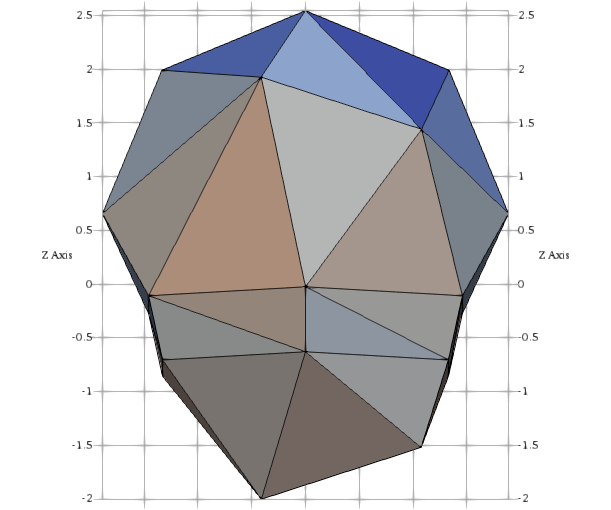}   &
		\imagespacing
		\includegraphics[width=\linewidth]{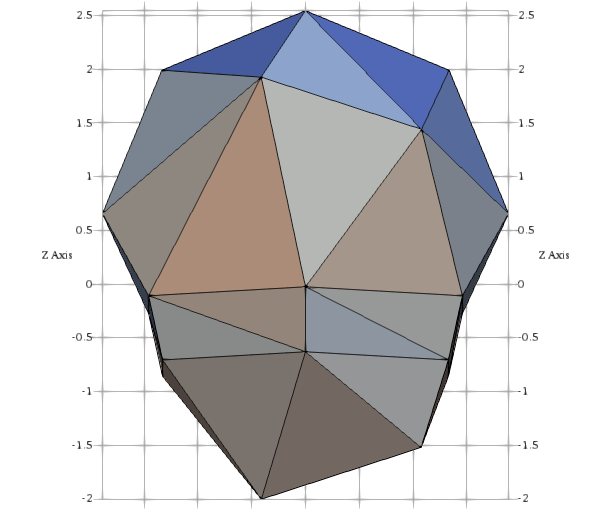}   &    \\ \cline{1-7}
	\end{tabular}
\end{center}
\vspace{0.25in}
\hspace*{3in}
{\Large
	\begin{minipage}[t]{3in}
		\baselineskip = .5\baselineskip
		Figure 6 \\
		Vicente Ramm, Jehanzeb H. Chaudhry, Christopher D. Cooper \\
		J.\ Comput.\ Chem.
	\end{minipage}
}
\clearpage

\begin{center}
	\begin{tabular}{m{1cm}m{1cm}|m{2.1cm}|m{2.1cm}|m{2.1cm}|m{2.1cm}|m{2.1cm}|m{2.1cm}} \cline{3-7}
		& & \multicolumn{5}{c|}{$N_{\phi}$} & \\
		\cline{2-7}
		& \multicolumn{1}{|c|}{View} & \multicolumn{1}{c|}{86} & \multicolumn{1}{c|}{344} & \multicolumn{1}{c|}{1376} &  \multicolumn{1}{c|}{5504} & \multicolumn{1}{c|}{22016} &  \\ \cline{1-7}
		\multicolumn{1}{|r}{\multirow{2}{*}{ \begin{minipage}{.9cm} \vspace{36pt}\centering
					$E_{\phi}$
		\end{minipage} }} & \multicolumn{1}{|c|}{Front} &
		\imagespacing
		\includegraphics[width=\linewidth]{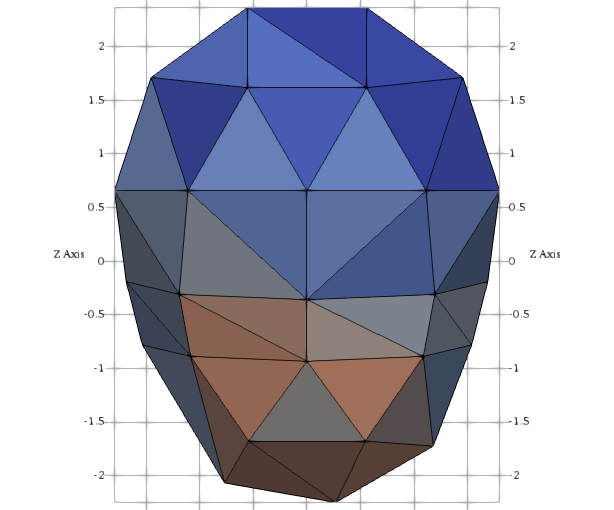} &	
		\imagespacing
		\includegraphics[width=\linewidth]{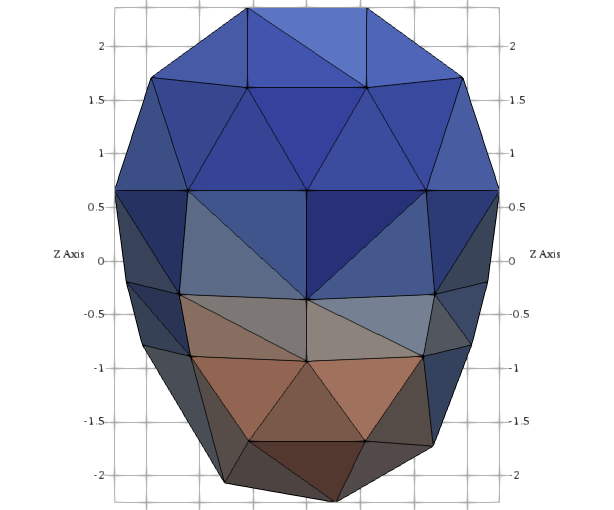} &	
		\imagespacing
		\includegraphics[width=\linewidth]{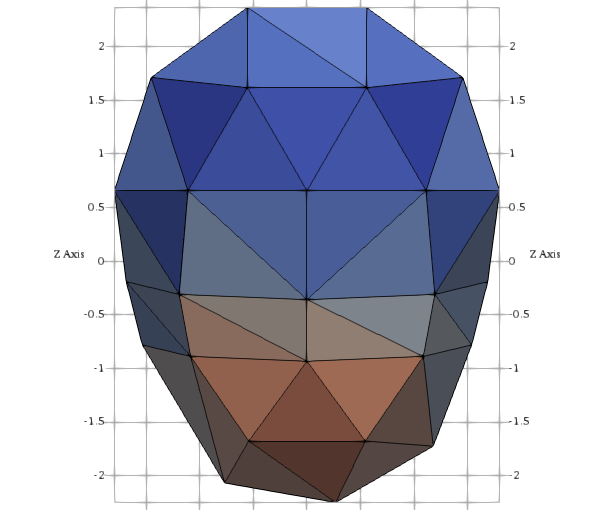} &	
		\imagespacing
		\includegraphics[width=\linewidth]{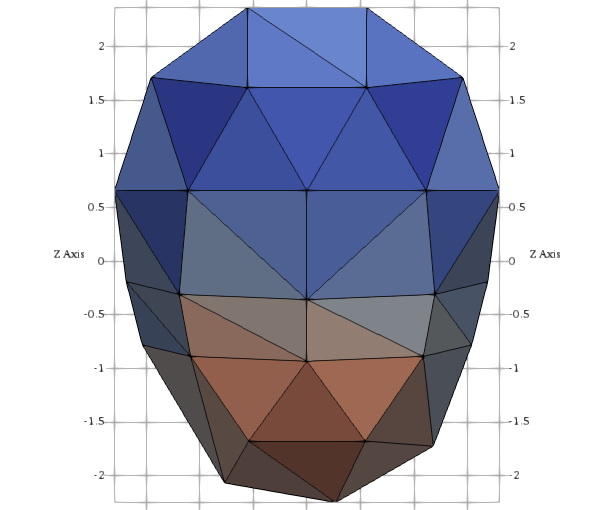} &	
		\imagespacing
		\includegraphics[width=\linewidth]{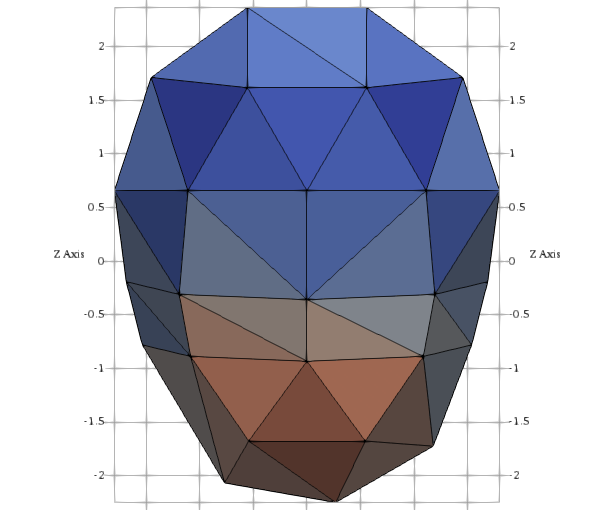} &	\multirow{4}{\linewidth}{\vspace{2pt}
			\begin{minipage}{2.1cm}
				\includegraphics[height=2.4\linewidth]{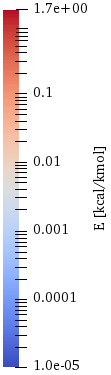} 
			\end{minipage} 
		} \\ \cline{2-7}
		\multicolumn{1}{|c}{} & \multicolumn{1}{|c|}{Back} &
		\imagespacing \includegraphics[width=\linewidth]{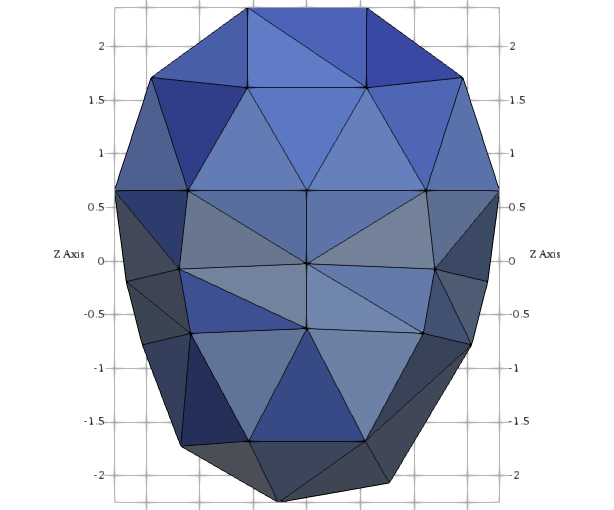} &
		\imagespacing	
		\includegraphics[width=\linewidth]{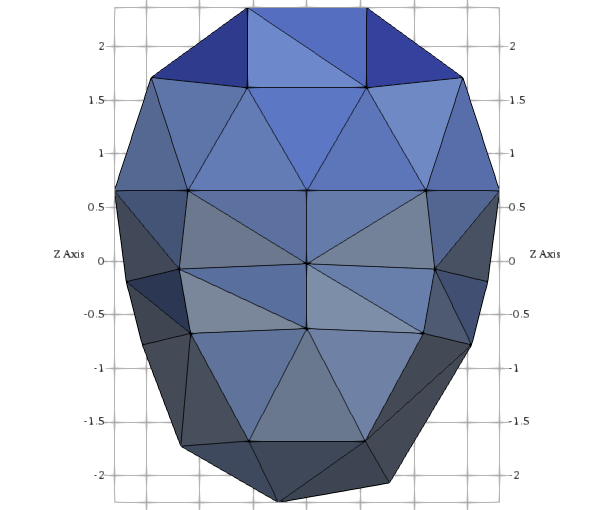} &
		\imagespacing	
		\includegraphics[width=\linewidth]{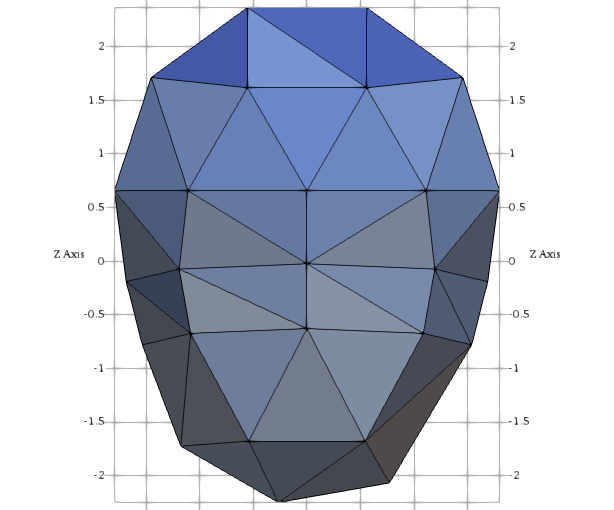} &
		\imagespacing	
		\includegraphics[width=\linewidth]{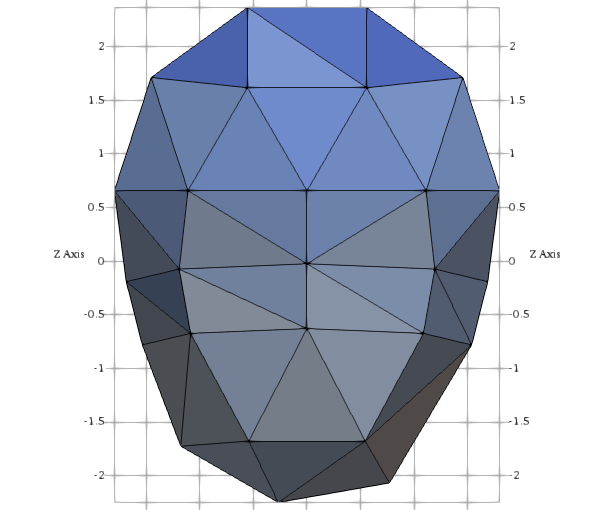} &	
		\imagespacing
		\includegraphics[width=\linewidth]{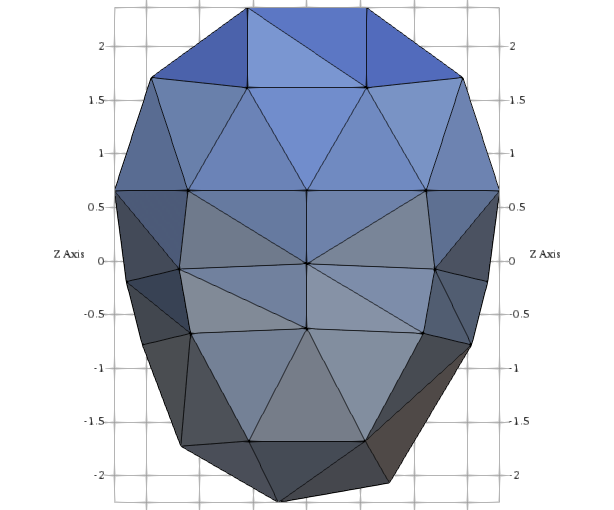} &	\\ \cline{1-7}
		\multicolumn{1}{|r}{\multirow{2}{*}{ \begin{minipage}{.9cm} \vspace{36pt}\centering
					$E_{u}$
		\end{minipage} }}& \multicolumn{1}{|c|}{Front}  &  \includegraphics[width=\linewidth]{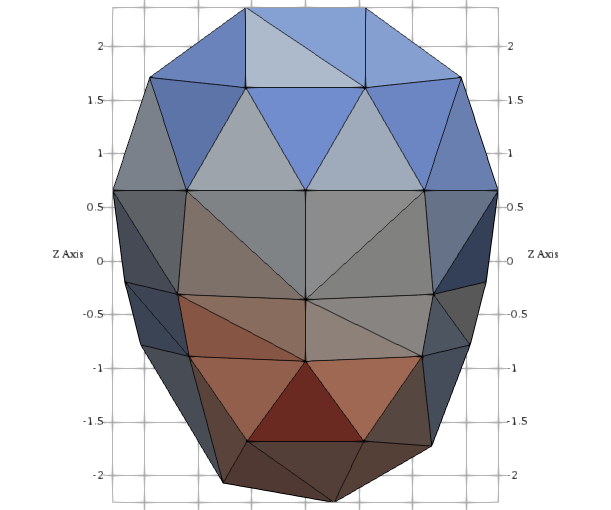}   &
		\imagespacing
		\includegraphics[width=\linewidth]{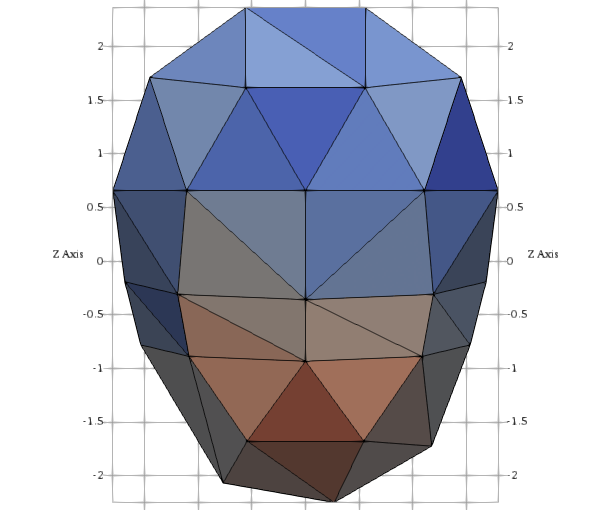}   &
		\imagespacing \includegraphics[width=\linewidth]{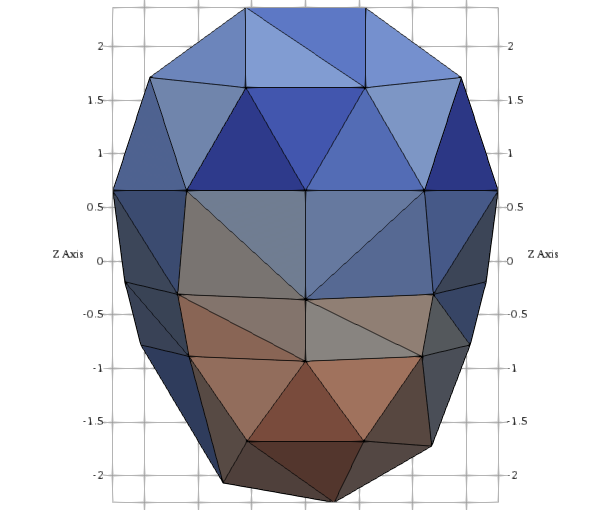}   &
		\imagespacing \includegraphics[width=\linewidth]{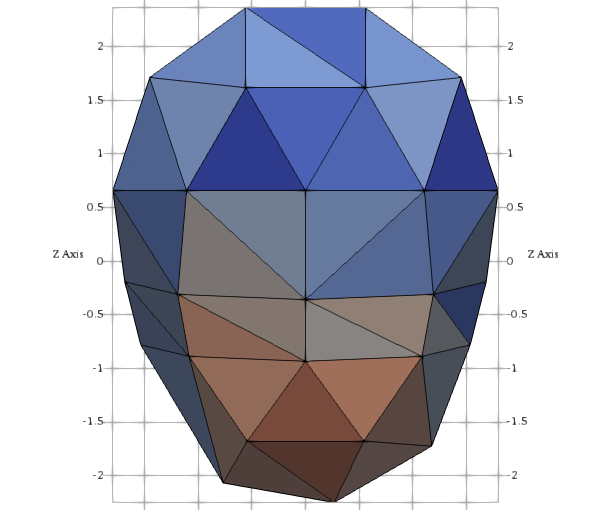}   &
		\imagespacing \includegraphics[width=\linewidth]{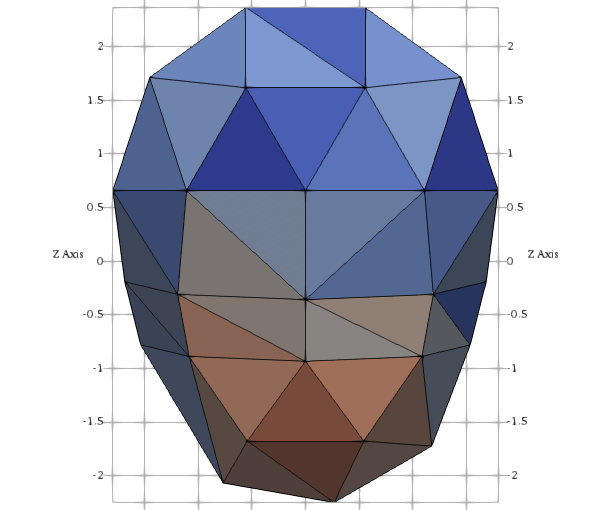}   &
		\\ \cline{2-7}
		\multicolumn{1}{|c}{} & \multicolumn{1}{|c|}{Back}  & \includegraphics[width=\linewidth]{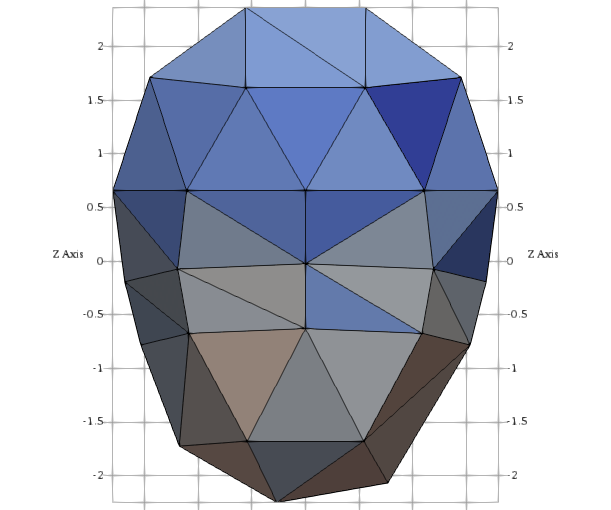}   &
		\imagespacing
		\includegraphics[width=\linewidth]{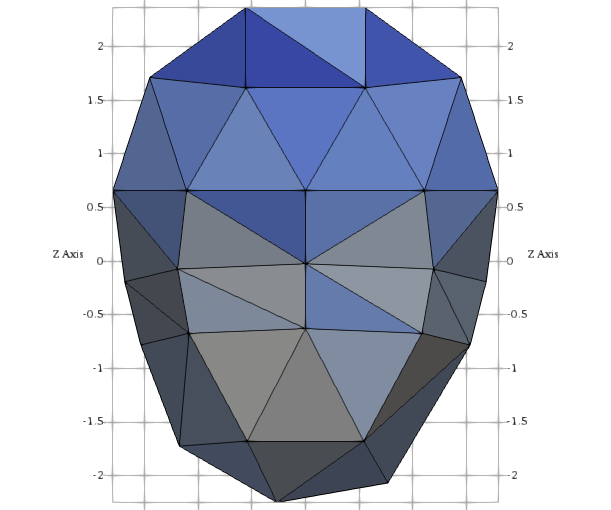}   &
		\imagespacing \includegraphics[width=\linewidth]{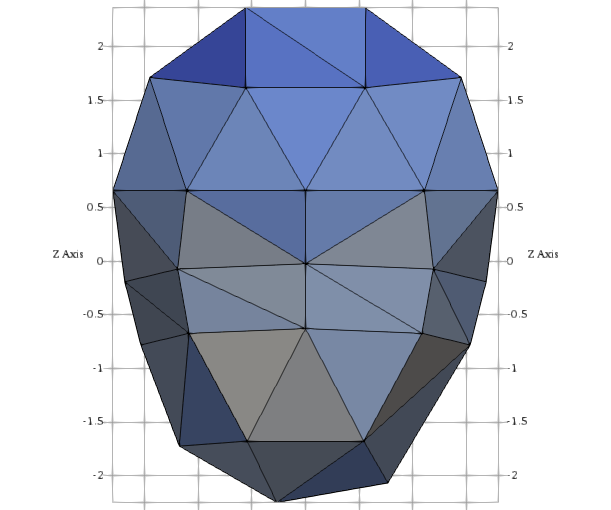}   &
		\imagespacing \includegraphics[width=\linewidth]{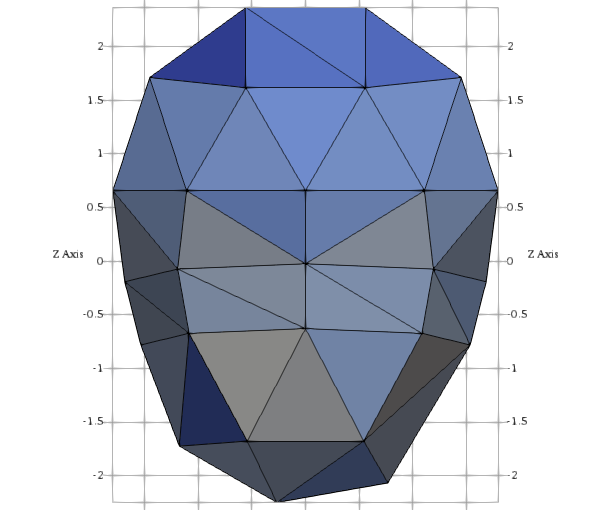}   &
		\imagespacing \includegraphics[width=\linewidth]{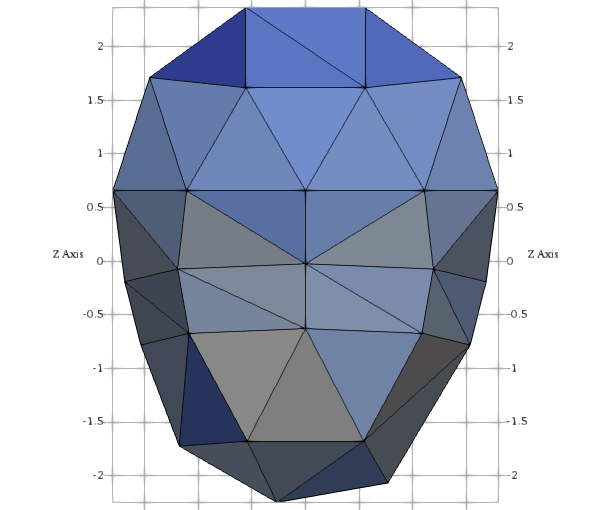}   &   \\ \cline{1-7}
	\end{tabular}
\end{center}
\vspace{0.25in}
\hspace*{3in}
{\Large
	\begin{minipage}[t]{3in}
		\baselineskip = .5\baselineskip
		Figure 7 \\
		Vicente Ramm, Jehanzeb H. Chaudhry, Christopher D. Cooper \\
		J.\ Comput.\ Chem.
	\end{minipage}
}
\clearpage

\begin{center}
			\begin{center}
			\newlength{\radius}
			\setlength{\radius}{2cm}
			\scalebox{0.7}{
				\begin{tikzpicture}
				\draw (0,0) circle [radius=\radius];
				\draw[dotted] (0,0) circle [radius=0.5*\radius];
				\fill (0,0.5*\radius) circle [radius=0.1];
				\node at (\radius*0.1,\radius*0.6) {$+q$};
				\draw[<-]  (\radius*0.707107 , \radius*0.707107 ) -- (\radius , \radius );
				\draw      (\radius , \radius ) -- (\radius + 2cm , \radius );
				\node at   (\radius + 1cm , \radius + 10pt ) {$R= 1 [\AA]$};
				\draw[<-]  (0.5*\radius*0.866025 , 0.5*\radius*0.5 ) -- (1.2*\radius*0.866025 , 1.2*\radius*0.5 );
				\draw      (1.2*\radius*0.866025 , 1.2*\radius*0.5 ) -- (1.2*\radius*0.866025 + 2cm , 1.2*\radius*0.5 );
				\node at   (1.2*\radius*0.866025 + 1cm , 1.2*\radius*0.5 + 10pt ) {$R= 0.5 [\AA]$};
				\draw[->]  (-\radius , -\radius) -- (-\radius        , -\radius+ 0.5cm) ;
				\draw[->]  (-\radius , -\radius) -- (-\radius + 0.5cm, -\radius) ;
				\node at (-\radius +2pt   , -\radius+ 0.6cm) {$z$} ;
				\node at (-\radius + 0.6cm, -\radius +2pt ) {$y$};
				\end{tikzpicture}}
			\end{center}

			\begin{center}
				\setlength{\radius}{2cm}
				\scalebox{0.7}{
					\begin{tikzpicture}
					\draw (0,0) circle [radius=\radius];
					\draw[dotted] (0,0) circle [radius=0.62*\radius];
					\fill (0,0.62*\radius) circle [radius=0.1];
					\node at (\radius*0.1,\radius*0.7) {$+q$};
					\draw (0.62*\radius*-0.08715574274,0.62*\radius*-0.99619469809) circle [radius=0.1];
					\node at (-\radius*0.2,-\radius*0.5) {$-q$};
					\fill (0.62*\radius* 0.08715574274,0.62*\radius*-0.99619469809) circle [radius=0.09];
					\node at (\radius*0.2,-\radius*0.5) {$+q$};
					\draw      (0,0) -- (0.9*\radius*-0.08715574274,0.9*\radius*-0.99619469809);
					\draw      (0,0) -- (0.9*\radius* 0.08715574274,0.9*\radius*-0.99619469809);
					\draw      (0,0) -- (0                         , 0.9*\radius );
					\draw[<-]  (\radius*0.707107 , \radius*0.707107 ) -- (\radius , \radius );
					\draw      (\radius , \radius ) -- (\radius + 2cm , \radius );
					\node at   (\radius + 1cm , \radius + 10pt ) {$R= 1 [\AA]$};
					\draw[<-]  (0.62*\radius*0.866025 , 0.62*\radius*0.5 ) -- (1.2*\radius*0.866025 , 1.2*\radius*0.5 );
					\draw      (1.2*\radius*0.866025 , 1.2*\radius*0.5 ) -- (1.2*\radius*0.866025 + 2cm , 1.2*\radius*0.5 );
					\node at   (1.2*\radius*0.866025 + 1cm , 1.2*\radius*0.5 + 10pt ) {$R= 0.62 [\AA]$};
					\draw[->]  (-\radius , -\radius) -- (-\radius        , -\radius+ 0.5cm) ;
					\draw[->]  (-\radius , -\radius) -- (-\radius + 0.5cm, -\radius) ;
					\node at (-\radius +2pt   , -\radius+ 0.6cm) {$z$} ;
					\node at (-\radius + 0.6cm, -\radius +2pt ) {$y$};
					\coordinate (A) at (0.4*\radius* 0.707107     ,0.4*\radius* -0.5     ) ;
					\coordinate (B) at (0.4*\radius*-0.707107     ,0.4*\radius* -0.5     ) ;
					\coordinate (T) at (0.4*\radius* 0.707107 + 0.5cm    ,0.4*\radius* -0.3 ) ;
					\draw[->] (B) arc (225:265:0.4*\radius);
					\draw[->] (A) arc (315:275:0.4*\radius);
					\node at  (T) {$\theta=10^{\circ}$};
					\end{tikzpicture}}
			\end{center}
\end{center}
\vspace{0.25in}
\hspace*{3in}
{\Large
	\begin{minipage}[t]{3in}
		\baselineskip = .5\baselineskip
		Figure 8 \\
		Vicente Ramm, Jehanzeb H. Chaudhry, Christopher D. Cooper \\
		J.\ Comput.\ Chem.
	\end{minipage}
}
\clearpage

\begin{center}

    \includegraphics[width=0.45\linewidth]{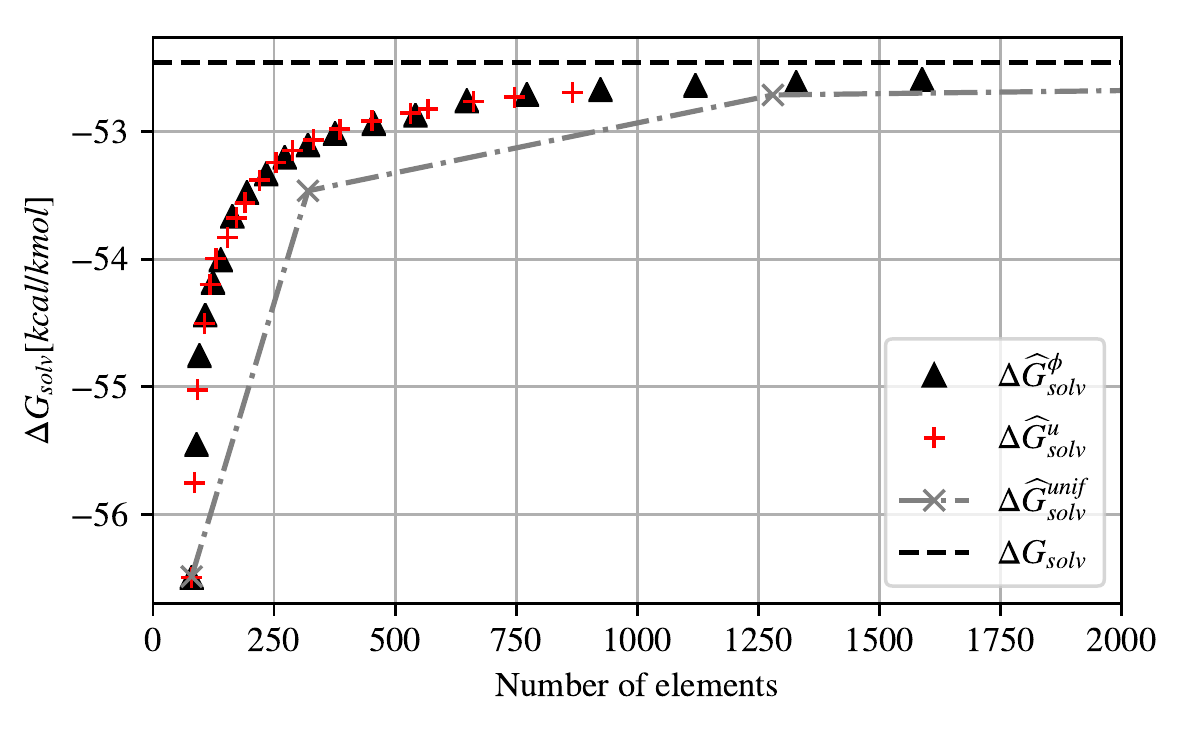}
    \includegraphics[width=0.45\linewidth]{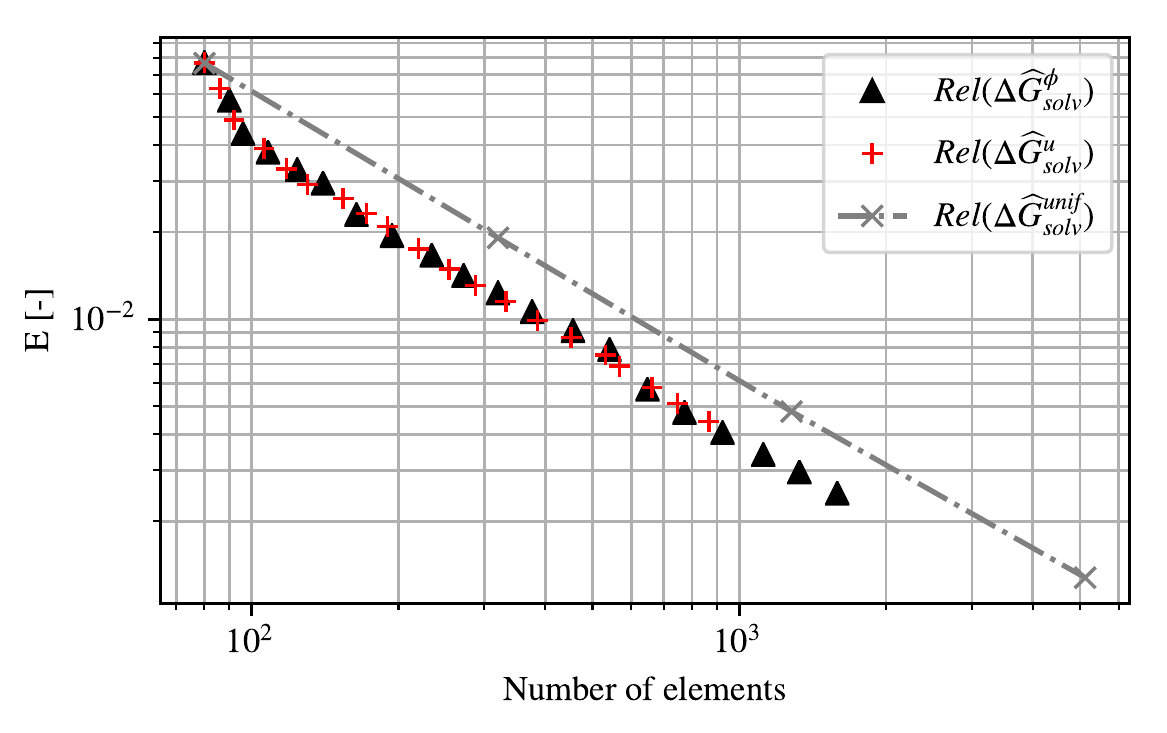}
\end{center}
\vspace{0.25in}
\hspace*{3in}
{\Large
	\begin{minipage}[t]{3in}
		\baselineskip = .5\baselineskip
		Figure 9 \\
		Vicente Ramm, Jehanzeb H. Chaudhry, Christopher D. Cooper \\
		J.\ Comput.\ Chem.
	\end{minipage}
}
\clearpage

\begin{center}
    \includegraphics[width=0.45\linewidth]{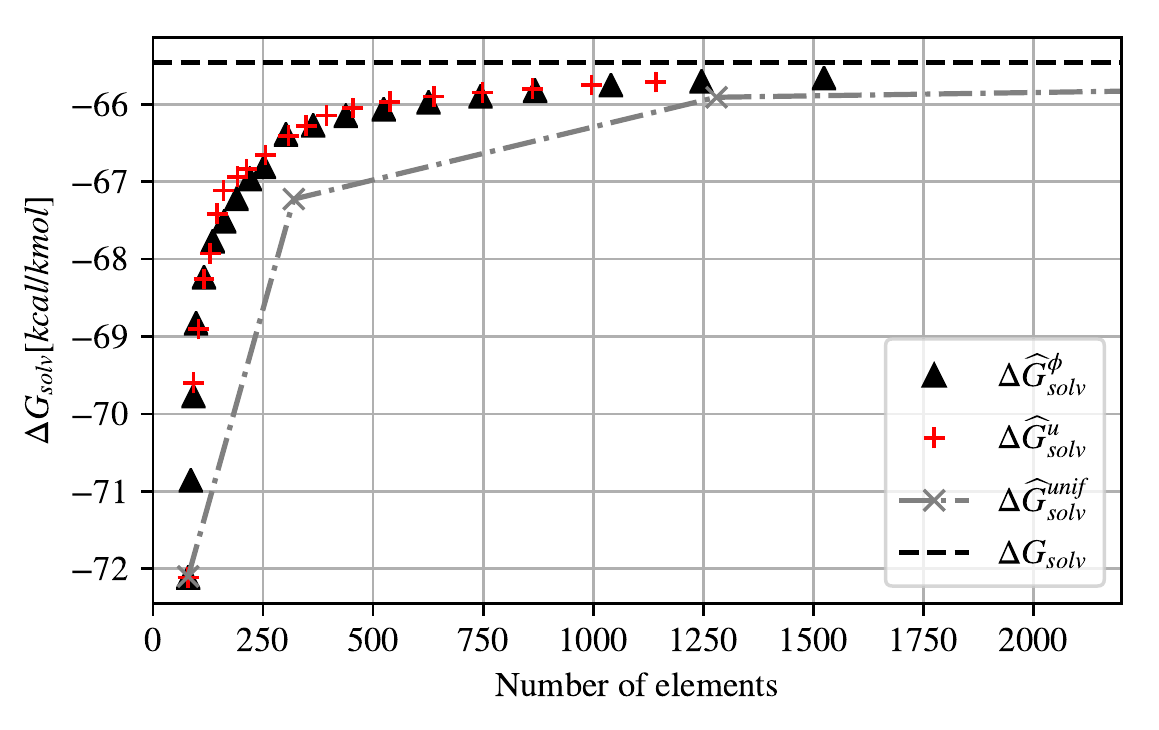}
    \includegraphics[width=0.45\linewidth]{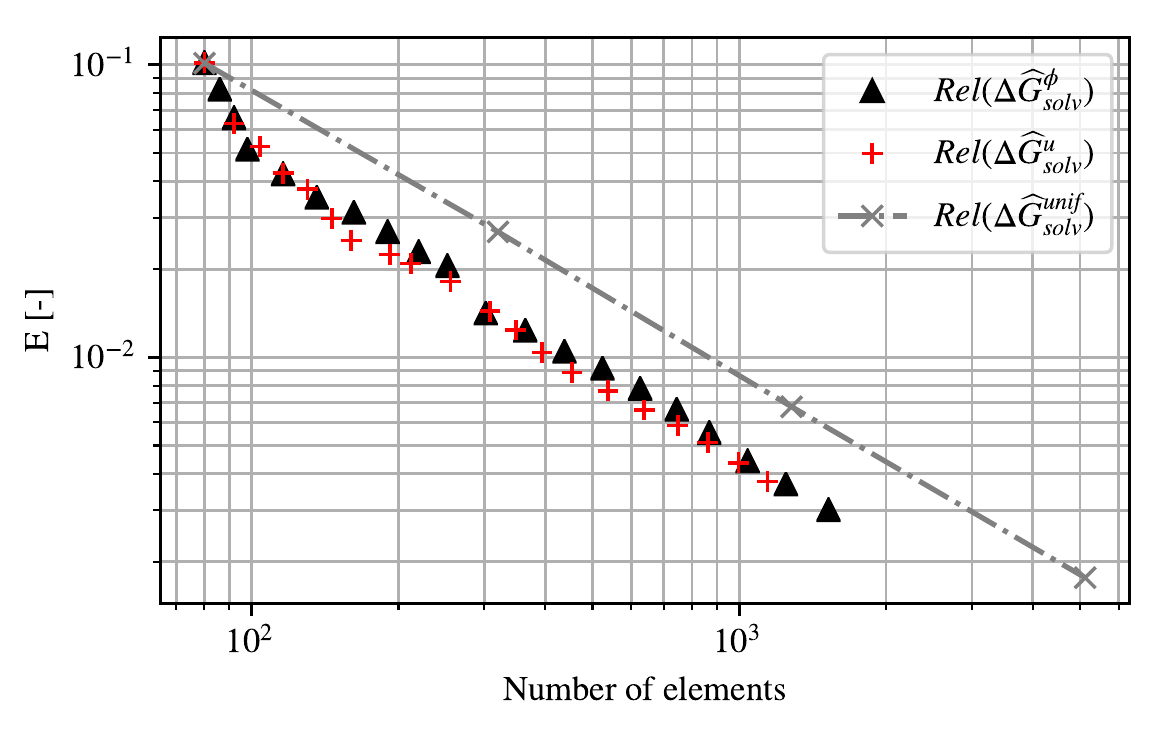}
\end{center}
\vspace{0.25in}
\hspace*{3in}
{\Large
	\begin{minipage}[t]{3in}
		\baselineskip = .5\baselineskip
		Figure 10 \\
		Vicente Ramm, Jehanzeb H. Chaudhry, Christopher D. Cooper \\
		J.\ Comput.\ Chem.
	\end{minipage}
}
\clearpage

\begin{center}
	\begin{tabular}{lm{.4cm}m{2.2cm}|m{2cm}|m{2cm}|m{2cm}|m{2cm}|m{4cm}}
		\cline{4-7}
		& &  & \multicolumn{4}{c|}{Iteration}  & \\
		\cline{3-7}
		& & \multicolumn{1}{|c|}{Estimator} & \multicolumn{1}{c|}{0} & \multicolumn{1}{c|}{10} & \multicolumn{1}{c|}{15} & \multicolumn{1}{c|}{20} &  \\  \cline{2-7}
		\multicolumn{1}{l|}{} & \multirow{2}{.4cm}{\rotatebox{90}{Charge-dipole}} & \multicolumn{1}{|c|}{$E_{u}$}    &
		\vspace{2pt}		 \includegraphics[width=\linewidth]{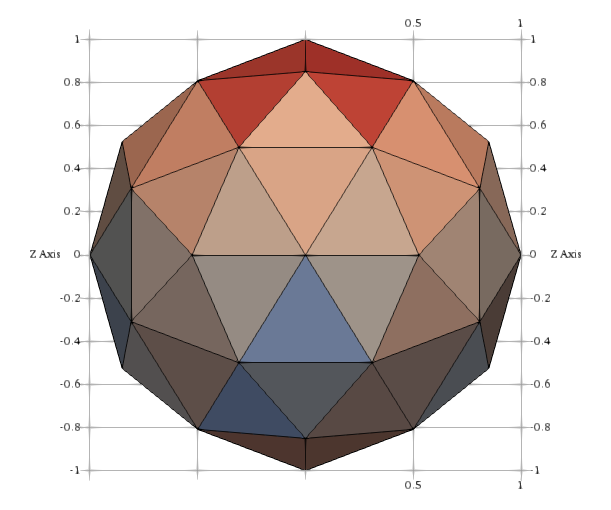}   &
		\vspace{2pt}	 \includegraphics[width=\linewidth]{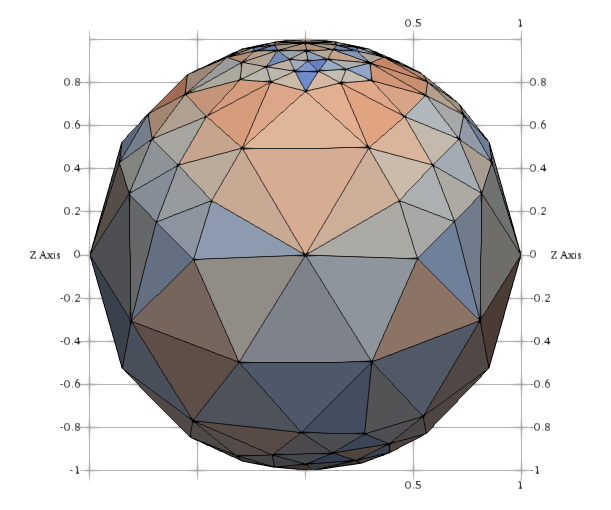}  &
		\vspace{2pt}	 \includegraphics[width=\linewidth]{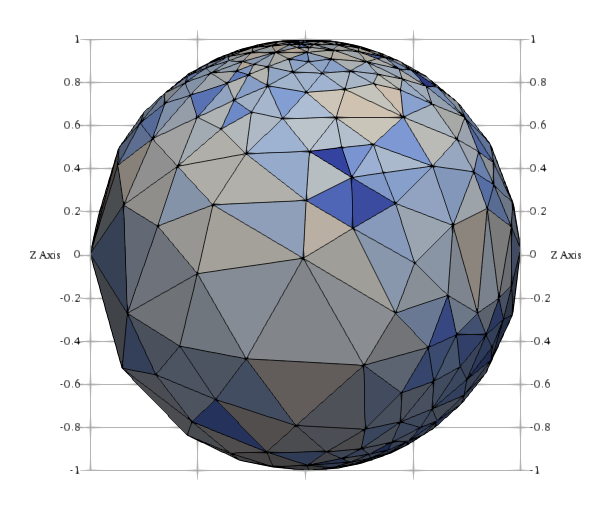}  &
		\vspace{2pt}	 \includegraphics[width=\linewidth]{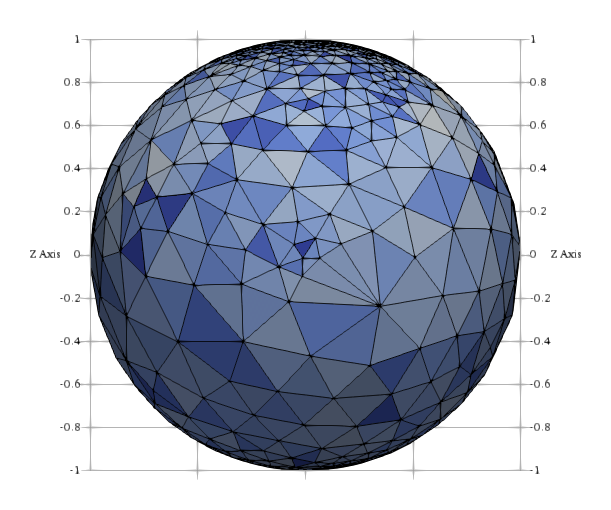}  & \multirow{4}{\linewidth}{\vspace{2pt}
			\begin{minipage}{2.1cm}
				\vspace{21pt}
				\includegraphics[height=5cm]{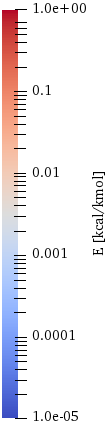}
			\end{minipage} 
		} \\ \cline{3-7}
		\multicolumn{1}{l|}{} & &\multicolumn{1}{|c|}{$E_{\phi}$}    &
		\vspace{2pt}	 \includegraphics[width=\linewidth]{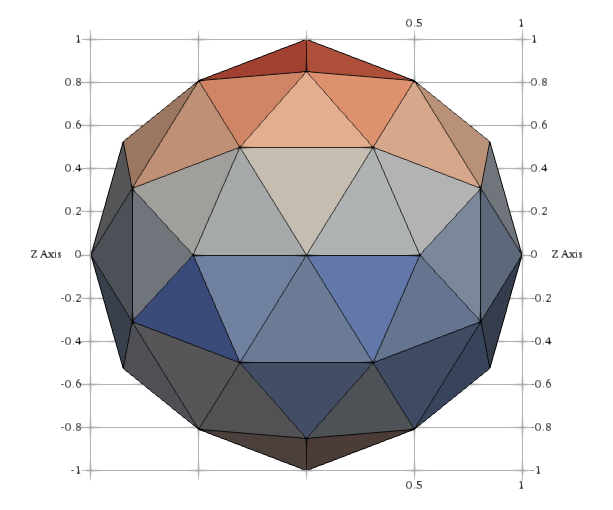}   &
		\vspace{2pt}	 \includegraphics[width=\linewidth]{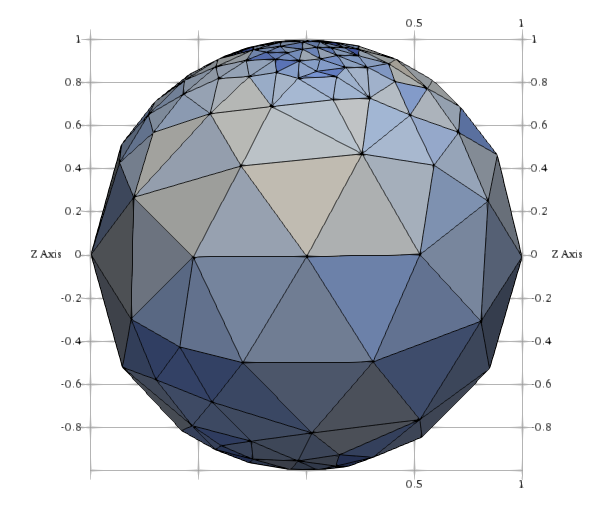}  &
		\vspace{2pt}	 \includegraphics[width=\linewidth]{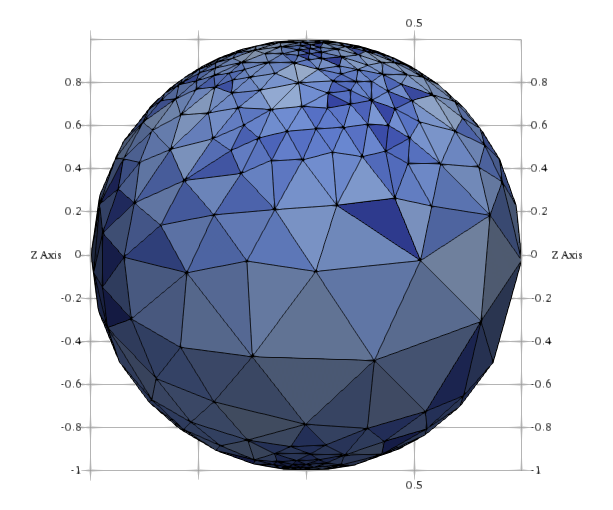}  &
		\vspace{2pt}	 \includegraphics[width=\linewidth]{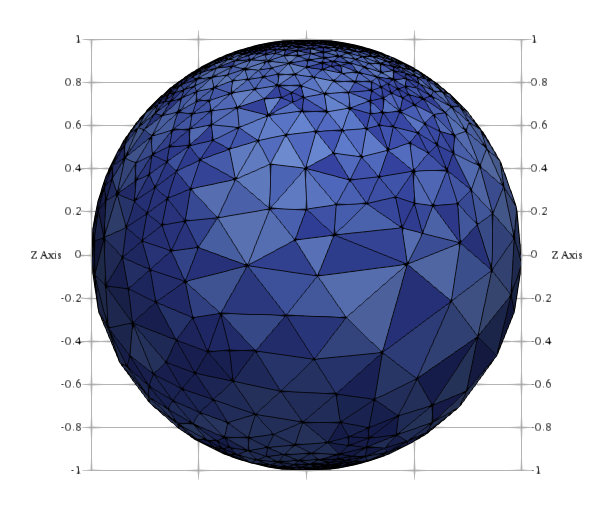}  &
		\\ \cline{2-7}
		\multicolumn{1}{l|}{} & \multirow{2}{\linewidth}{\rotatebox{90}{Off-centered}} & \multicolumn{1}{|c|}{$E_{u}$}    &
		\vspace{2pt}	 \includegraphics[width=\linewidth]{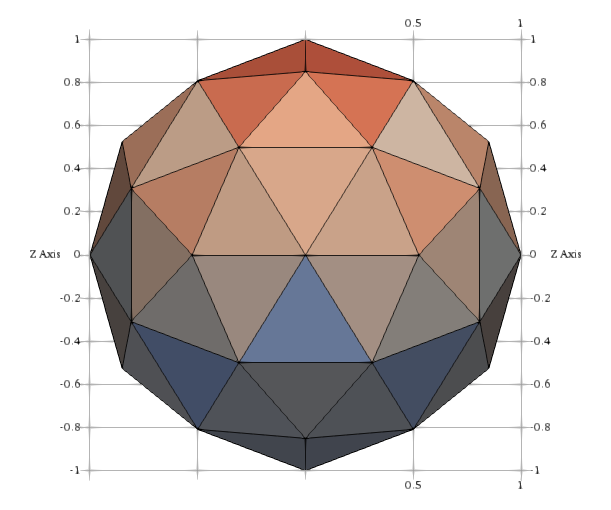}   &
		\vspace{2pt}	 \includegraphics[width=\linewidth]{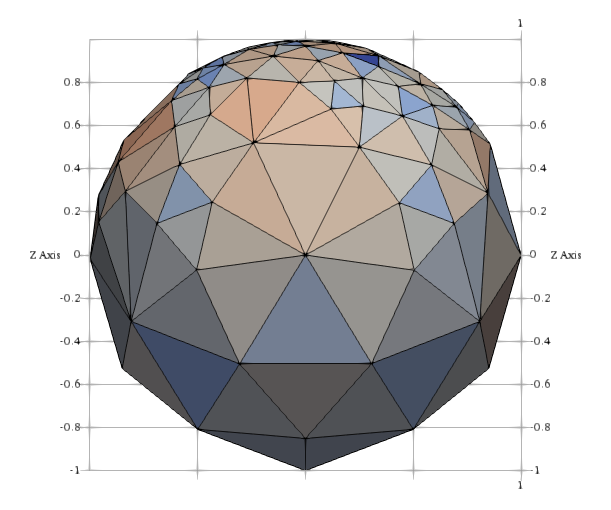}  &
		\vspace{2pt}	 \includegraphics[width=\linewidth]{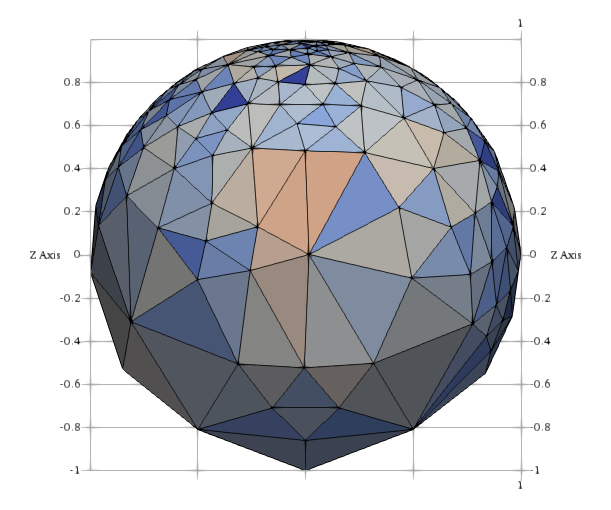}  &
		\vspace{2pt}	 \includegraphics[width=\linewidth]{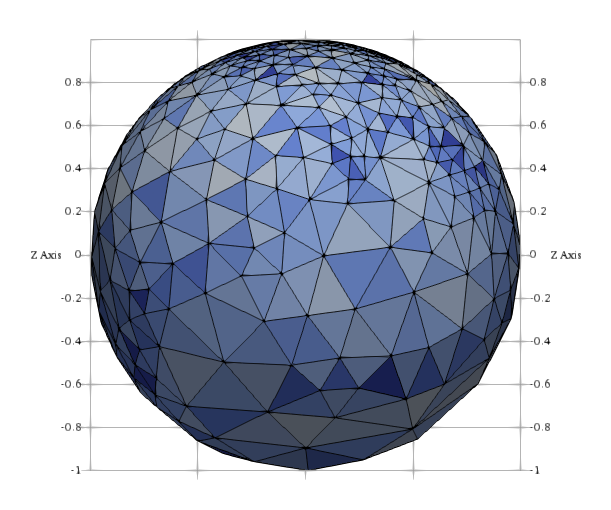}  & 
		\\ \cline{3-7}
		\multicolumn{1}{l|}{} & & \multicolumn{1}{|c|}{$E_{\phi}$}       &
		\vspace{2pt}	 \includegraphics[width=\linewidth]{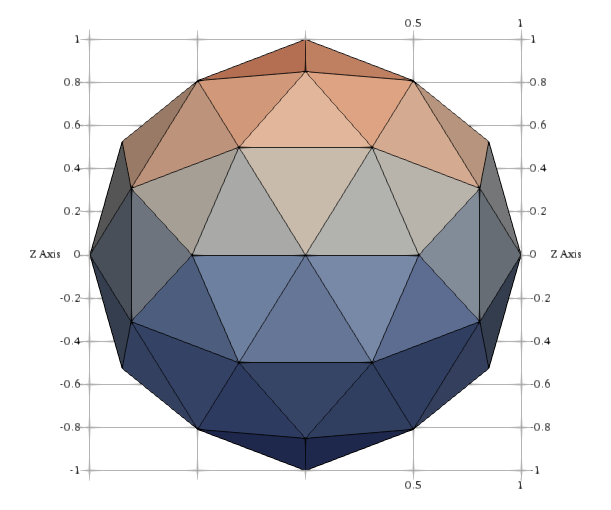}   &
		\vspace{2pt}	 \includegraphics[width=\linewidth]{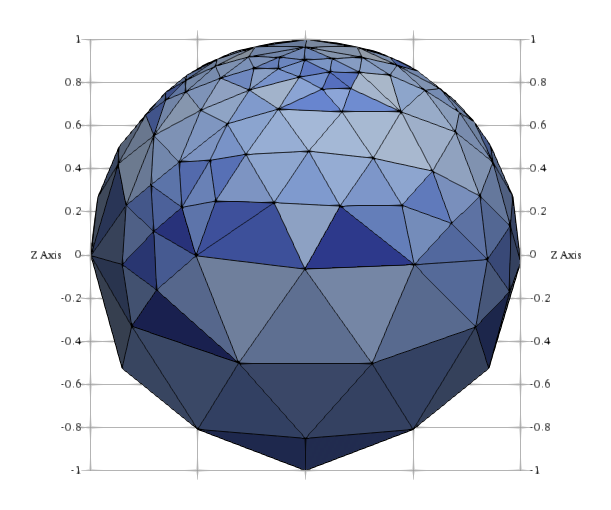}  &
		\vspace{2pt}	 \includegraphics[width=\linewidth]{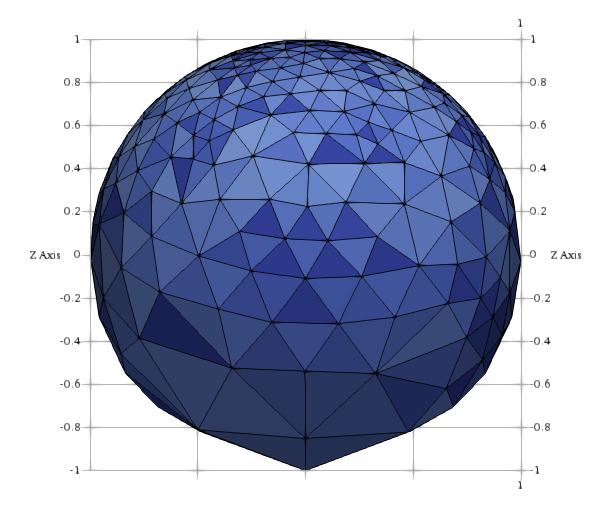}  &
		\vspace{2pt}	 \includegraphics[width=\linewidth]{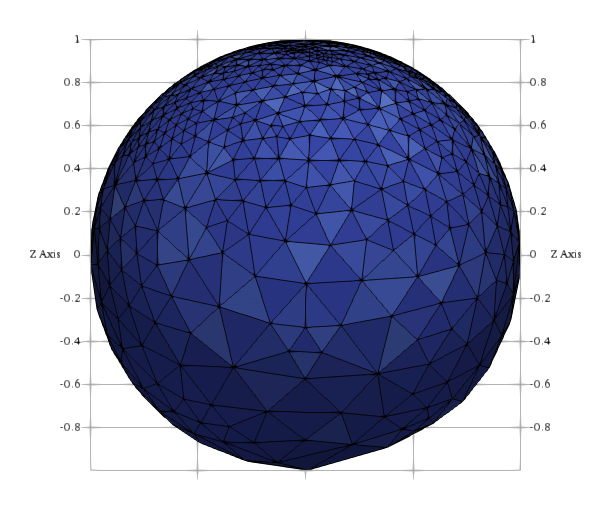}  & 
		\\ \cline{2-7}
	\end{tabular}
\end{center}
\vspace{0.25in}
\hspace*{3in}
{\Large
	\begin{minipage}[t]{3in}
		\baselineskip = .5\baselineskip
		Figure 11 \\
		Vicente Ramm, Jehanzeb H. Chaudhry, Christopher D. Cooper \\
		J.\ Comput.\ Chem.
	\end{minipage}
}
\clearpage

\renewcommand\cellspacetoplimit{2pt}
\renewcommand\cellspacebottomlimit{2pt}
\begin{center}
    \includegraphics[width=0.49\linewidth]{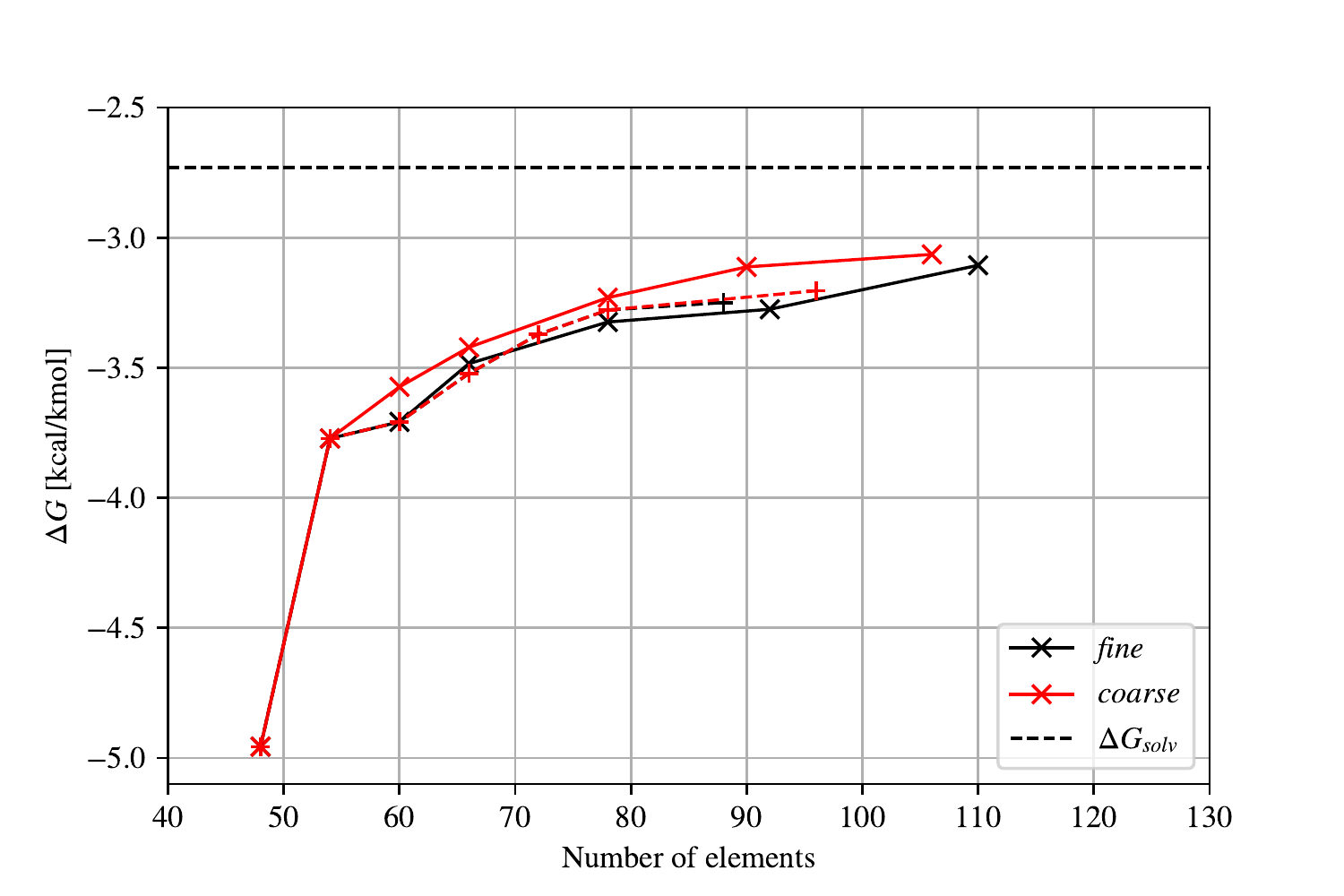}
    \includegraphics[width=0.49\linewidth]{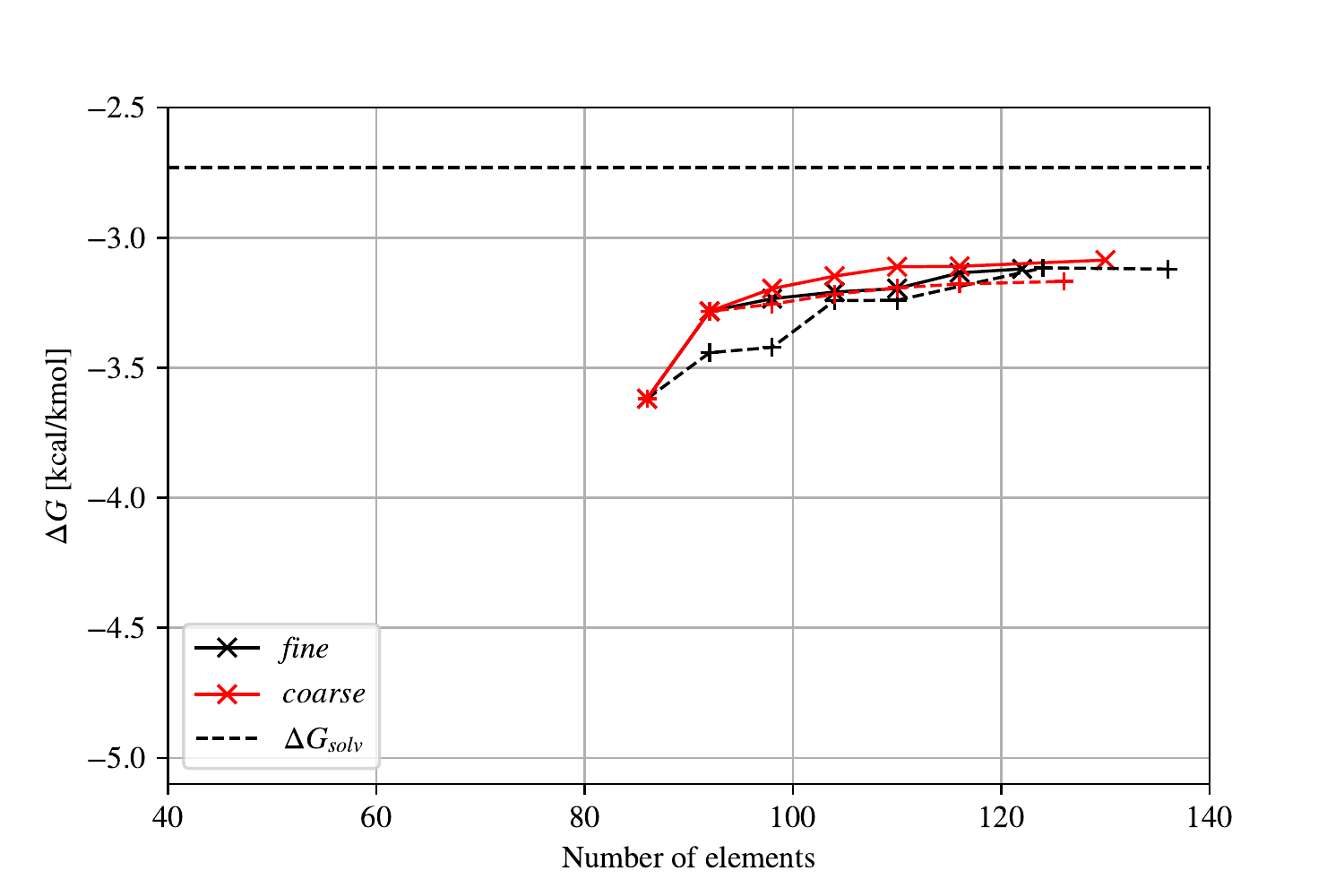}
    \includegraphics[width=0.49\linewidth]{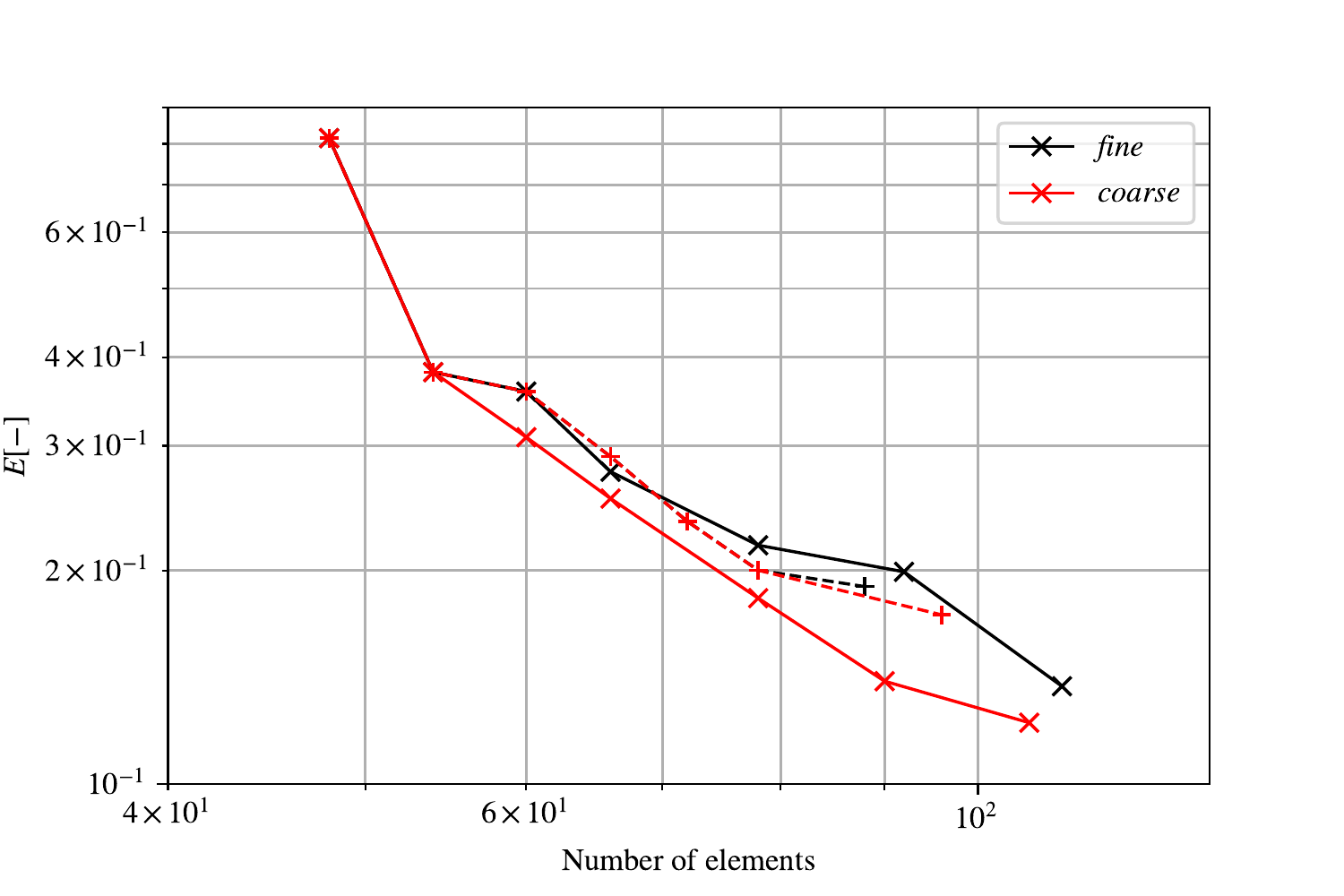}
    \includegraphics[width=0.49\linewidth]{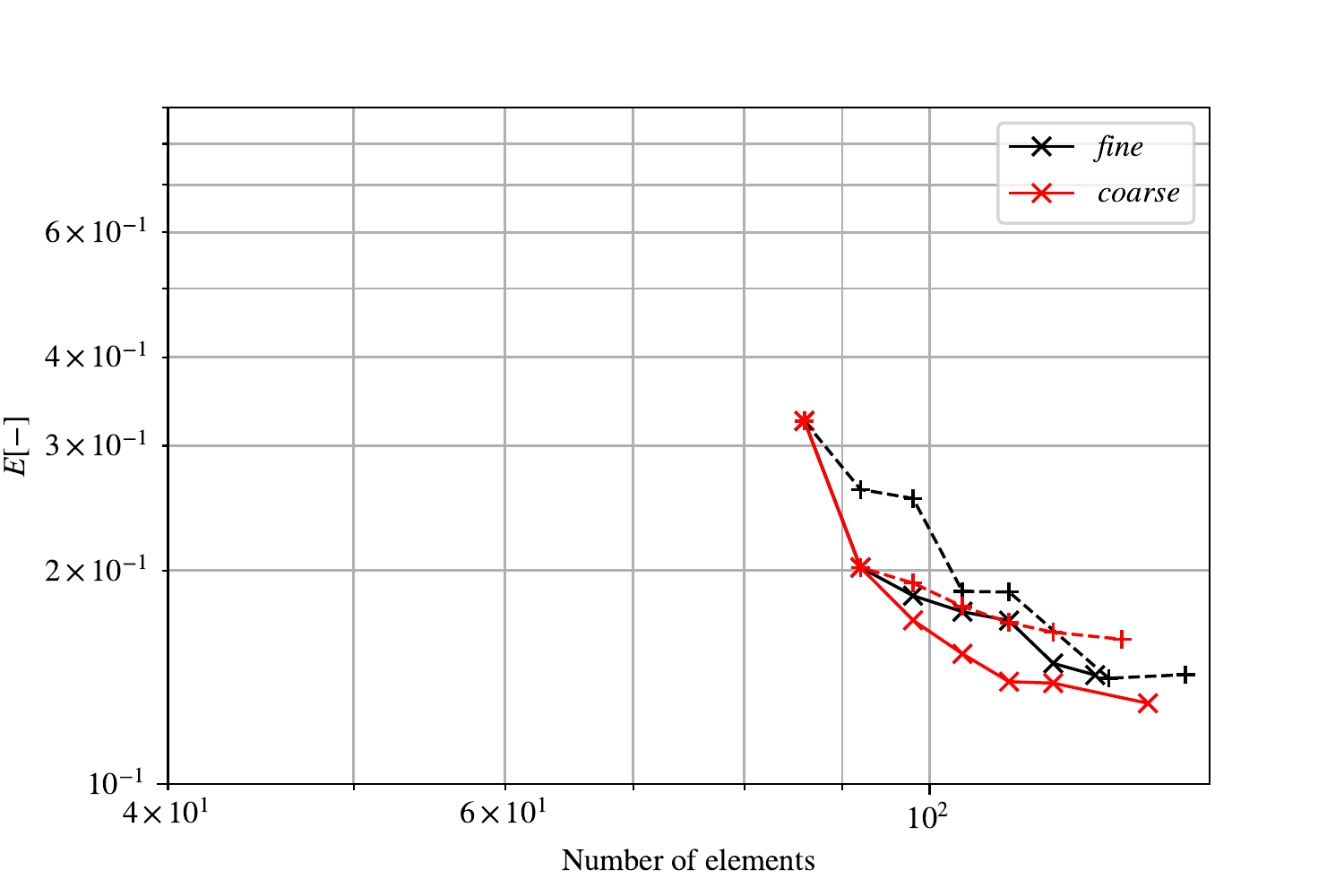}
\end{center}
\vspace{0.25in}
\hspace*{3in}
{\Large
	\begin{minipage}[t]{3in}
		\baselineskip = .5\baselineskip
		Figure 12 \\
		Vicente Ramm, Jehanzeb H. Chaudhry, Christopher D. Cooper \\
		J.\ Comput.\ Chem.
	\end{minipage}
}
\clearpage

\begin{center}
		\begin{tabular}{|c|c|m{1.8cm}|m{1.8cm}|m{1.8cm}|m{1.8cm}|m{1.8cm}}
			\cline{1-6}
			Mesh    & & \multicolumn{4}{c|}{Iteration} &  \\ \cline{3-6}
			\multicolumn{1}{|c|}{Density} & $\phi$ mesh & \multicolumn{1}{c|}{0}  & \multicolumn{1}{c|}{2}  & \multicolumn{1}{c|}{4}  & \multicolumn{1}{c|}{6}  & \\ \cline{1-6}
			0.5     & fine &
			\vspace{2pt}
			\includegraphics[width=\linewidth]{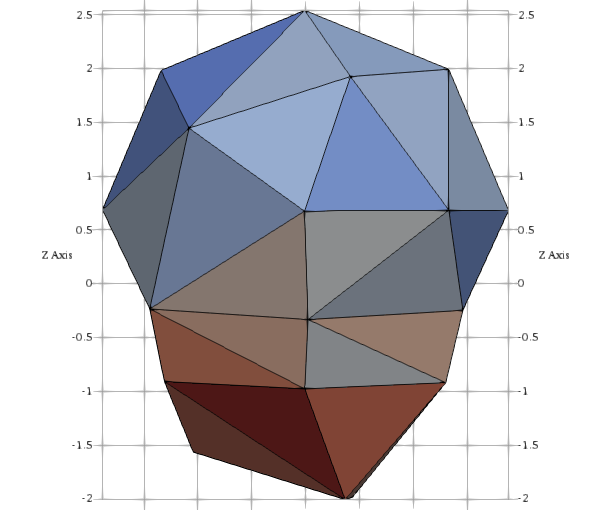}    &
			\vspace{2pt}
			\includegraphics[width=\linewidth]{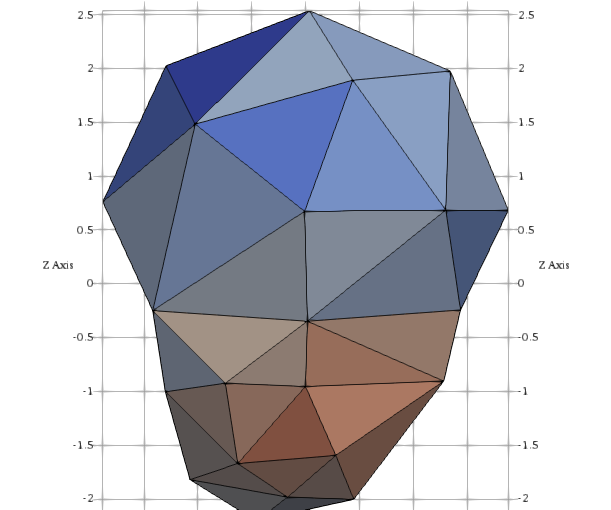}    &
			\vspace{2pt} \includegraphics[width=\linewidth]{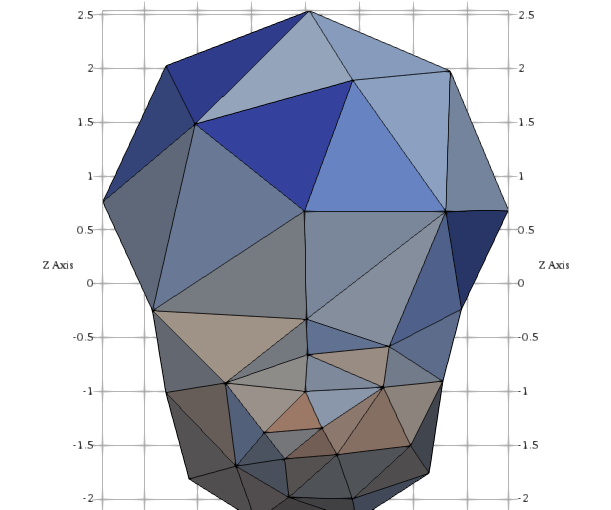}    &
			\vspace{2pt}
			\includegraphics[width=\linewidth]{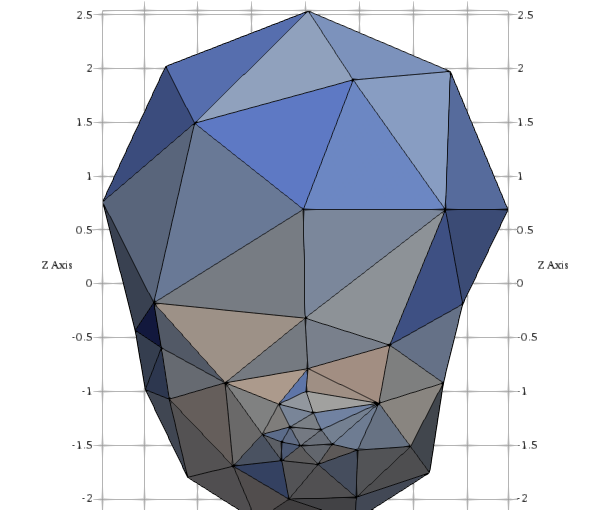}    &
			 \multirow{4}{\linewidth}{\vspace{2pt}
			 	\begin{minipage}{2.1cm}
			 		\includegraphics[width=0.6\linewidth]{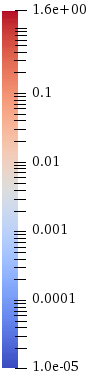}
			 	\end{minipage} 
			    }     
			\\ \cline{1-6}
			0.5     & coarse &
			\vspace{2pt}
			\includegraphics[width=\linewidth]{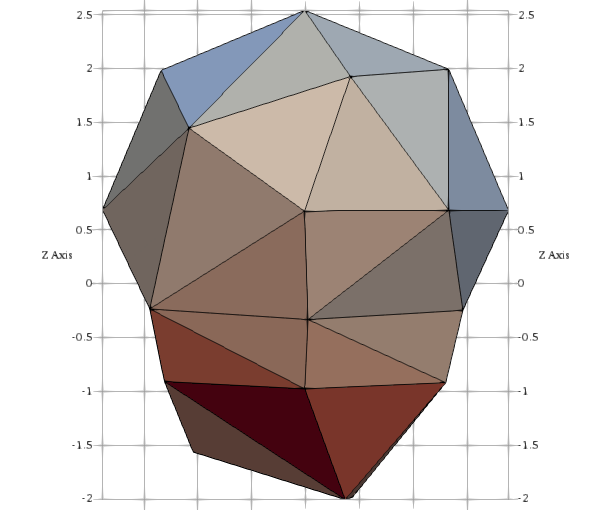}    &
			\vspace{2pt}
			\includegraphics[width=\linewidth]{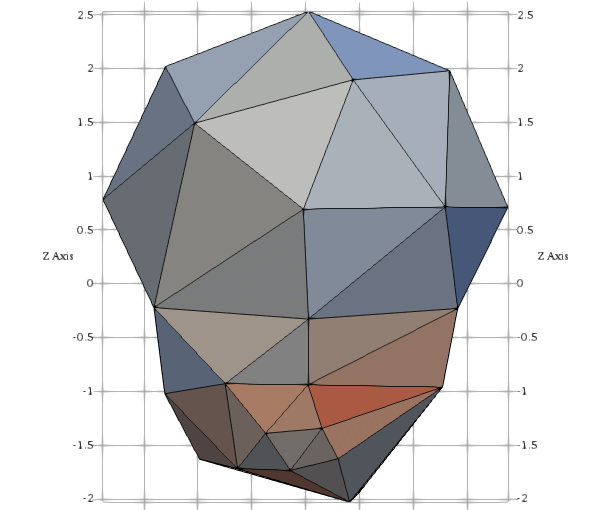}    &
			\vspace{2pt} \includegraphics[width=\linewidth]{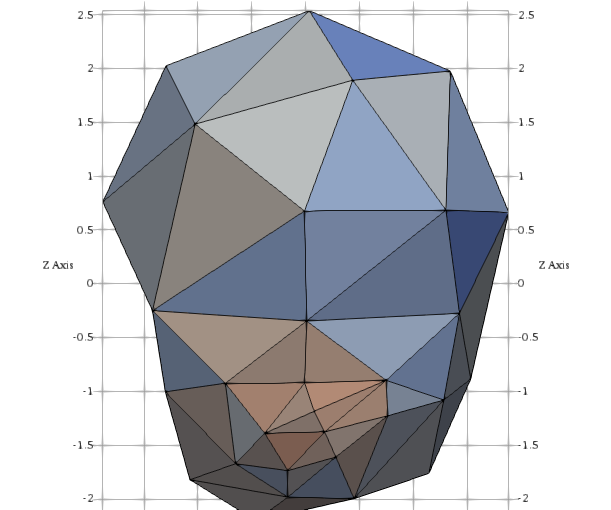}    &
			\vspace{2pt}
			\includegraphics[width=\linewidth]{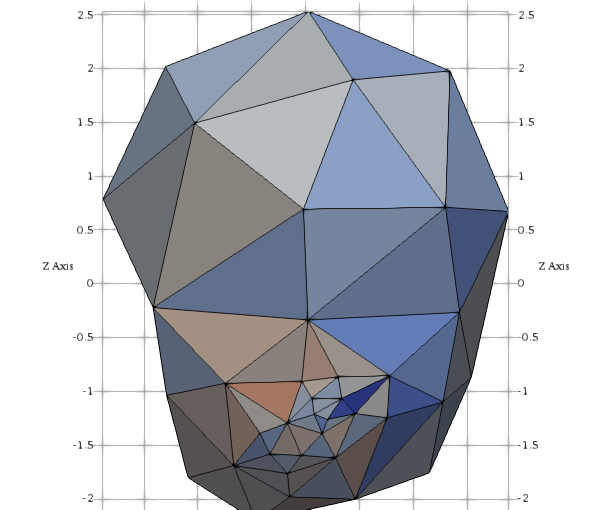}    &
			\\ \cline{1-6}
			1.0     & fine &
			\vspace{2pt}
			\includegraphics[width=\linewidth]{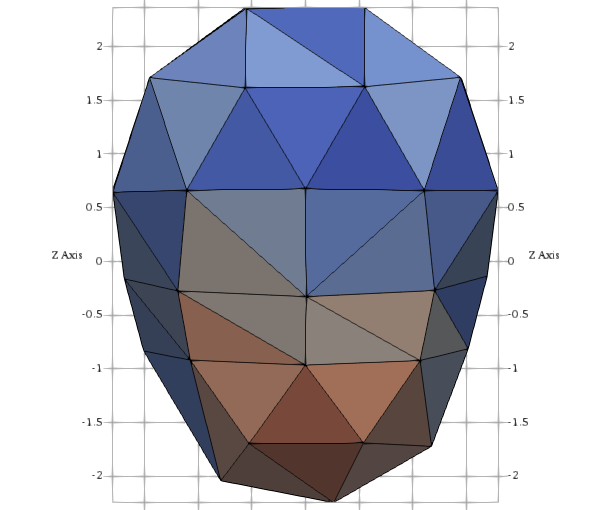}    &
			\vspace{2pt}
			\includegraphics[width=\linewidth]{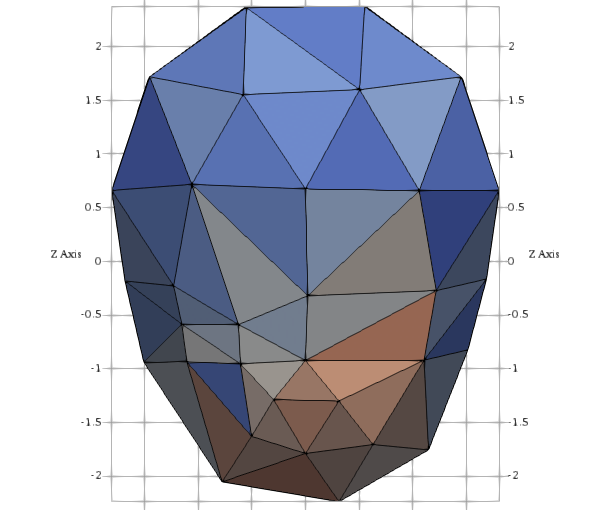}    &
			\vspace{2pt} \includegraphics[width=\linewidth]{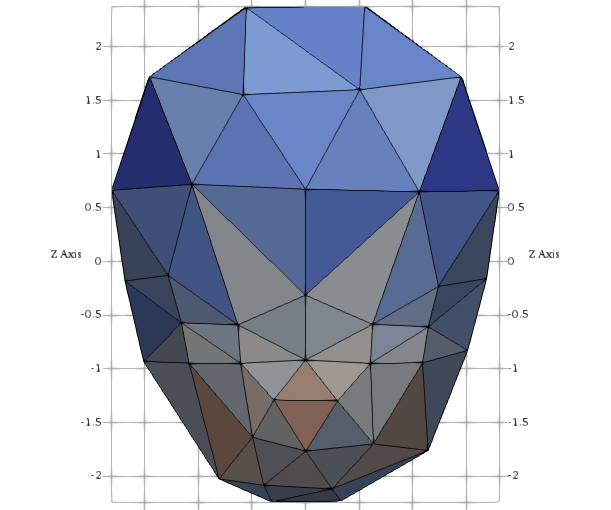}    &
			\vspace{2pt}
			\includegraphics[width=\linewidth]{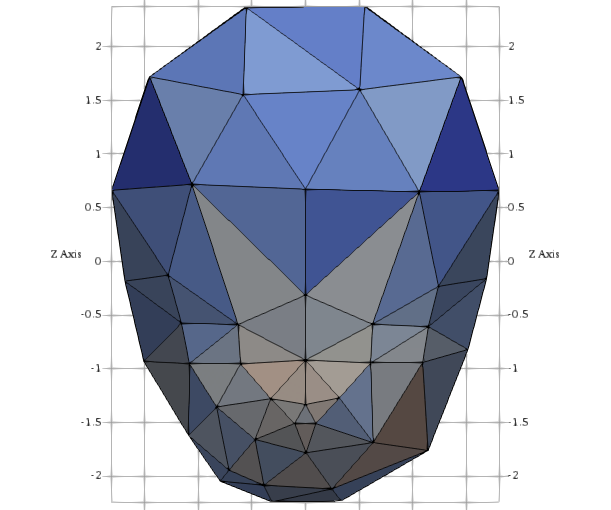}    &
			\\ \cline{1-6}
			1.0     & coarse &
			\vspace{2pt}
			\includegraphics[width=\linewidth]{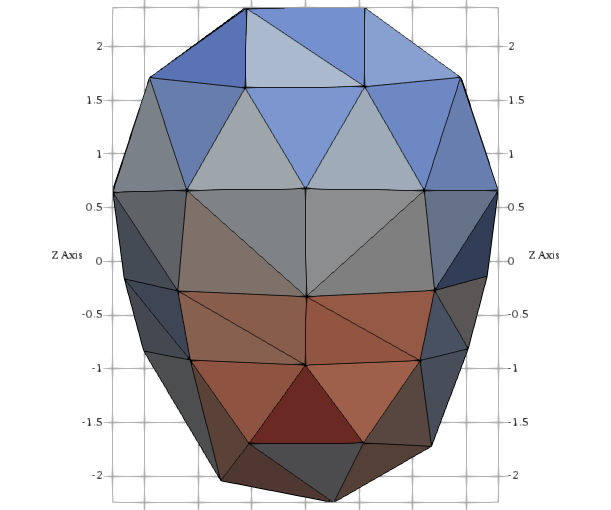}    &
			\vspace{2pt}
			\includegraphics[width=\linewidth]{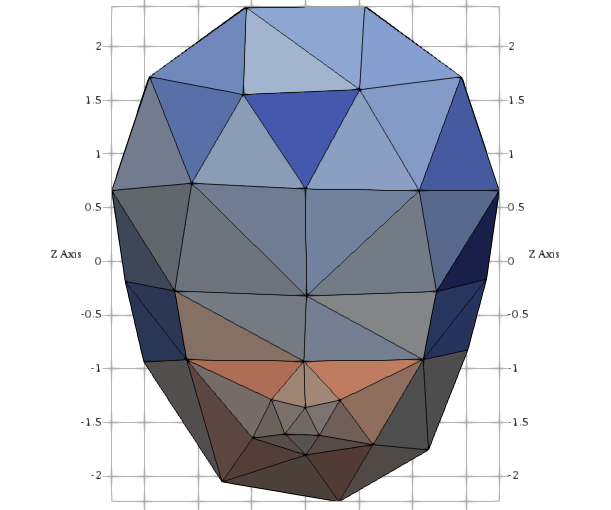}    &
			\vspace{2pt} \includegraphics[width=\linewidth]{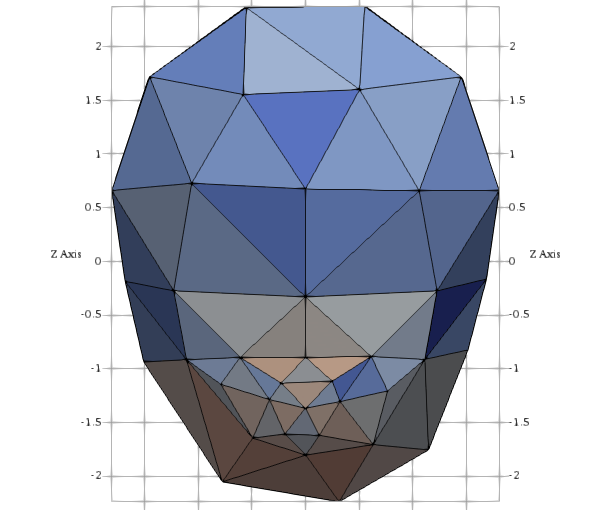}    &
			\vspace{2pt}
			\includegraphics[width=\linewidth]{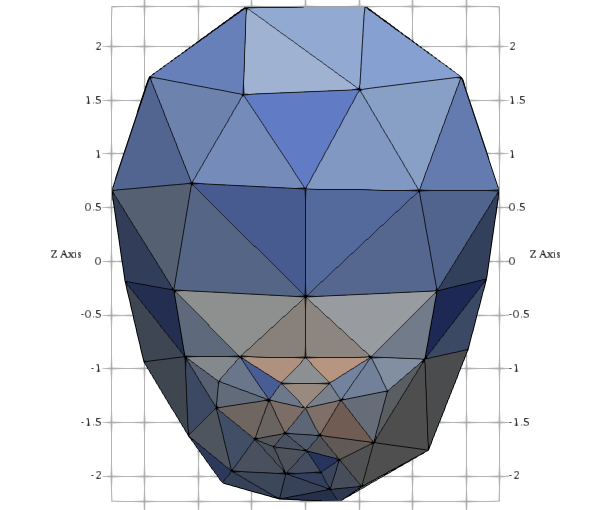}    & 
			\\ \cline{1-6}
		\end{tabular}
\end{center}
\vspace{0.25in}
\hspace*{3in}
{\Large
	\begin{minipage}[t]{3in}
		\baselineskip = .5\baselineskip
		Figure 13 \\
		Vicente Ramm, Jehanzeb H. Chaudhry, Christopher D. Cooper \\
		J.\ Comput.\ Chem.
	\end{minipage}
}
\clearpage

\begin{center}
		\begin{tabular}{|c|m{2cm}|m{2cm}|m{2cm}|m{2cm}|m{2cm}}
			\cline{1-5}
			Mesh    & \multicolumn{4}{c|}{Iteration} &  \\ \cline{2-5}
			\multicolumn{1}{|c|}{Density} & \multicolumn{1}{c|}{0}  & \multicolumn{1}{c|}{2}  & \multicolumn{1}{c|}{4}  & \multicolumn{1}{c|}{6}  & \\ \cline{1-5}
			0.5    &
			\vspace{2pt}
			\includegraphics[width=\linewidth]{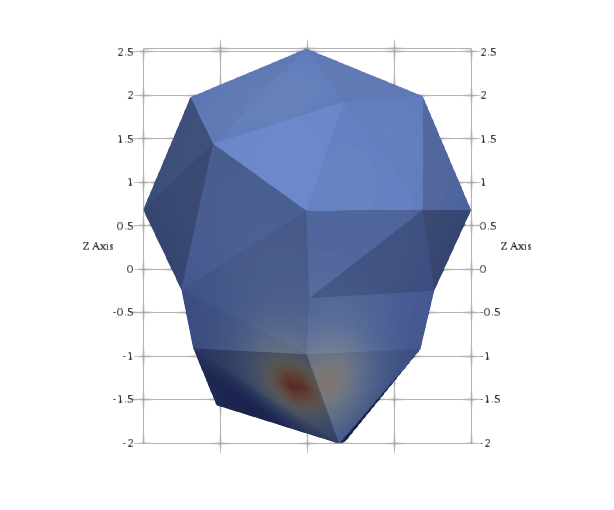}    &
			\vspace{2pt}
			\includegraphics[width=\linewidth]{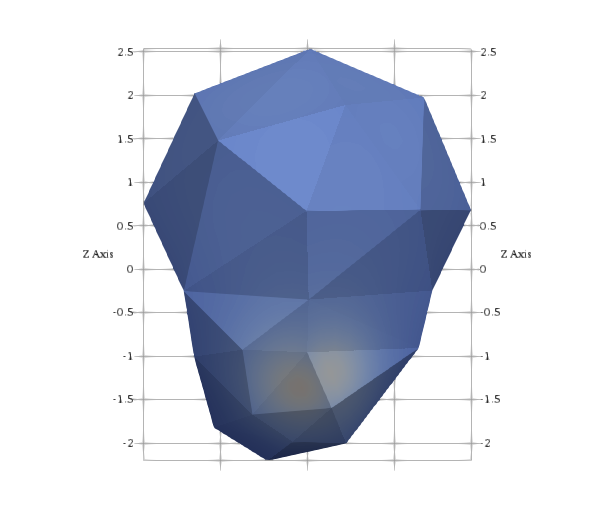}    &
			\vspace{2pt} \includegraphics[width=\linewidth]{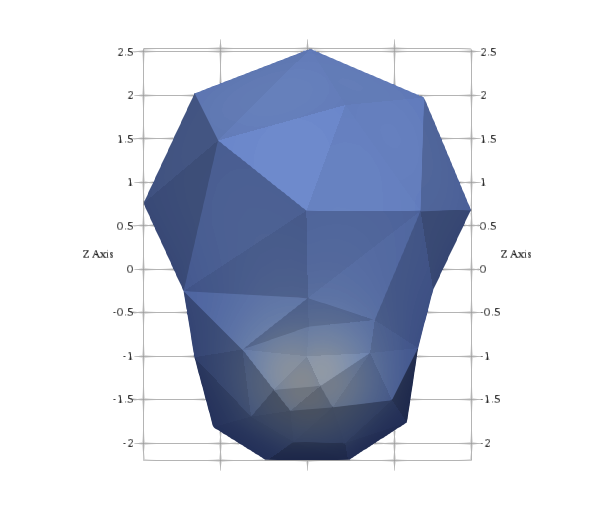}    &
			\vspace{2pt}
			\includegraphics[width=\linewidth]{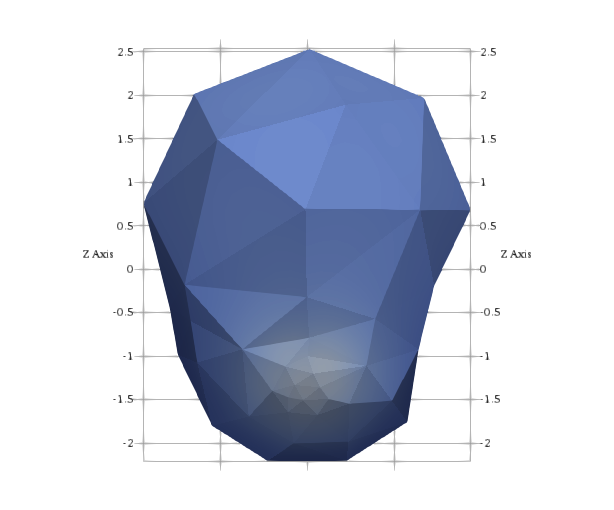}    &
			\multirow{2}{\linewidth}{
			\includegraphics[width=0.5\linewidth]{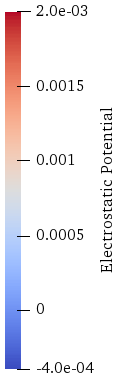}   
		}           
			\\ \cline{1-5}
			1.0        &
			\vspace{2pt}
			\includegraphics[width=\linewidth]{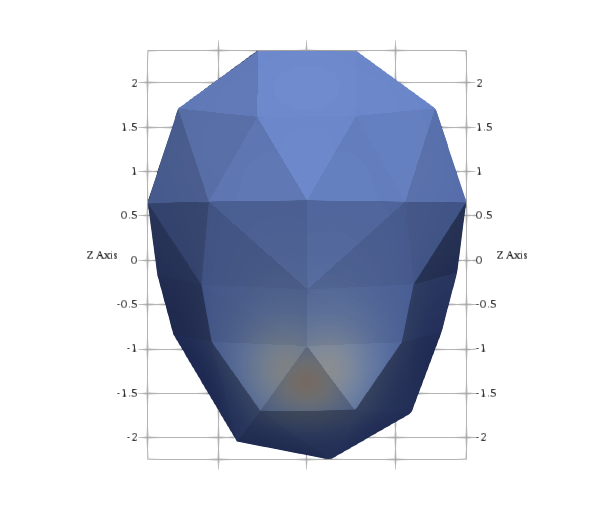}    &
			\vspace{2pt}
			\includegraphics[width=\linewidth]{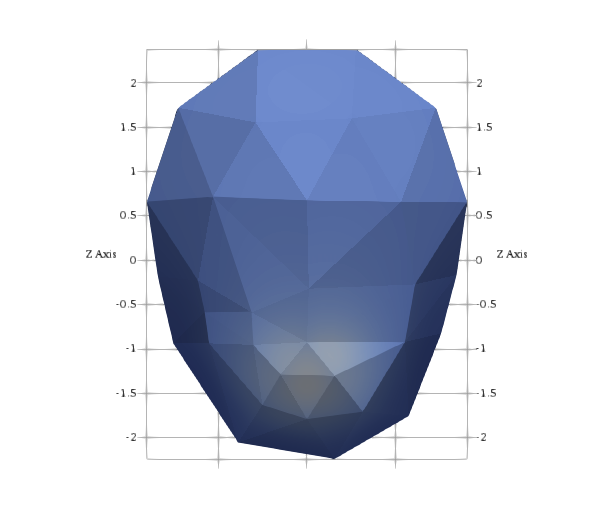}    &
			\vspace{2pt} \includegraphics[width=\linewidth]{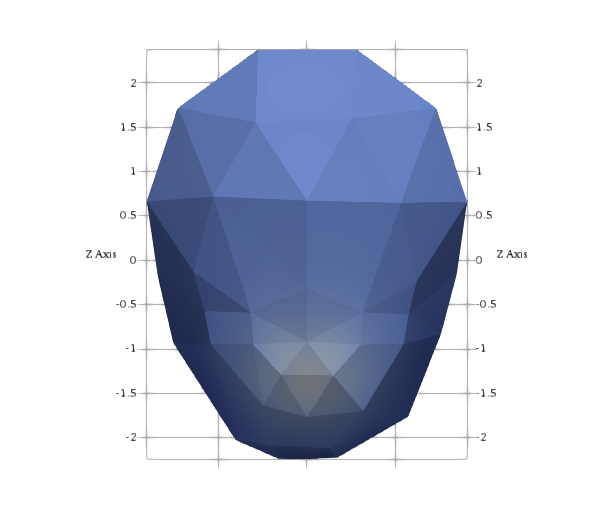}    &
			\vspace{2pt}
			\includegraphics[width=\linewidth]{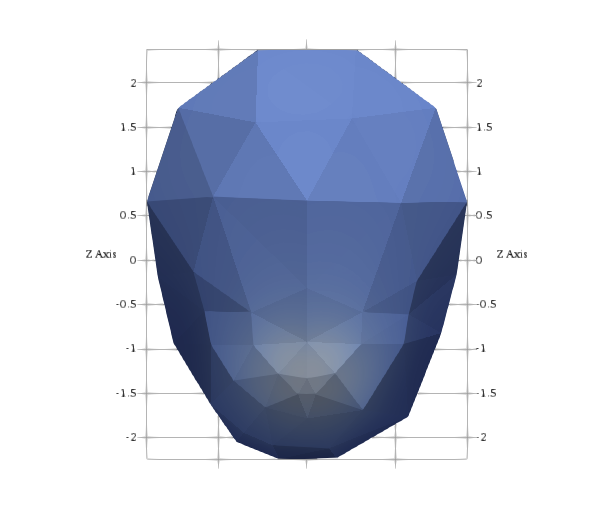}    & 
			\\ \cline{1-5}
		\end{tabular}
\end{center}
\vspace{0.25in}
\hspace*{3in}
{\Large
	\begin{minipage}[t]{3in}
		\baselineskip = .5\baselineskip
		Figure 14 \\
		Vicente Ramm, Jehanzeb H. Chaudhry, Christopher D. Cooper \\
		J.\ Comput.\ Chem.
	\end{minipage}
}
\clearpage

\begin{center}
    \includegraphics[width=0.45\textwidth]{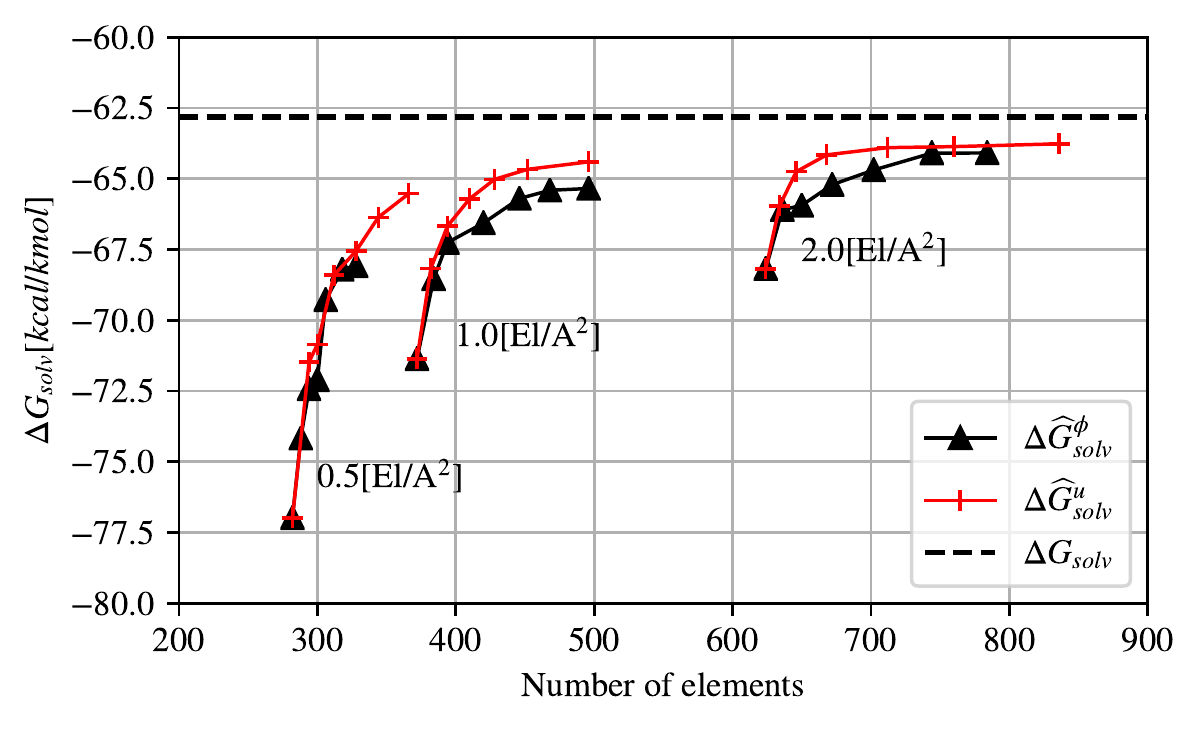} 
    \includegraphics[width=0.45\textwidth]{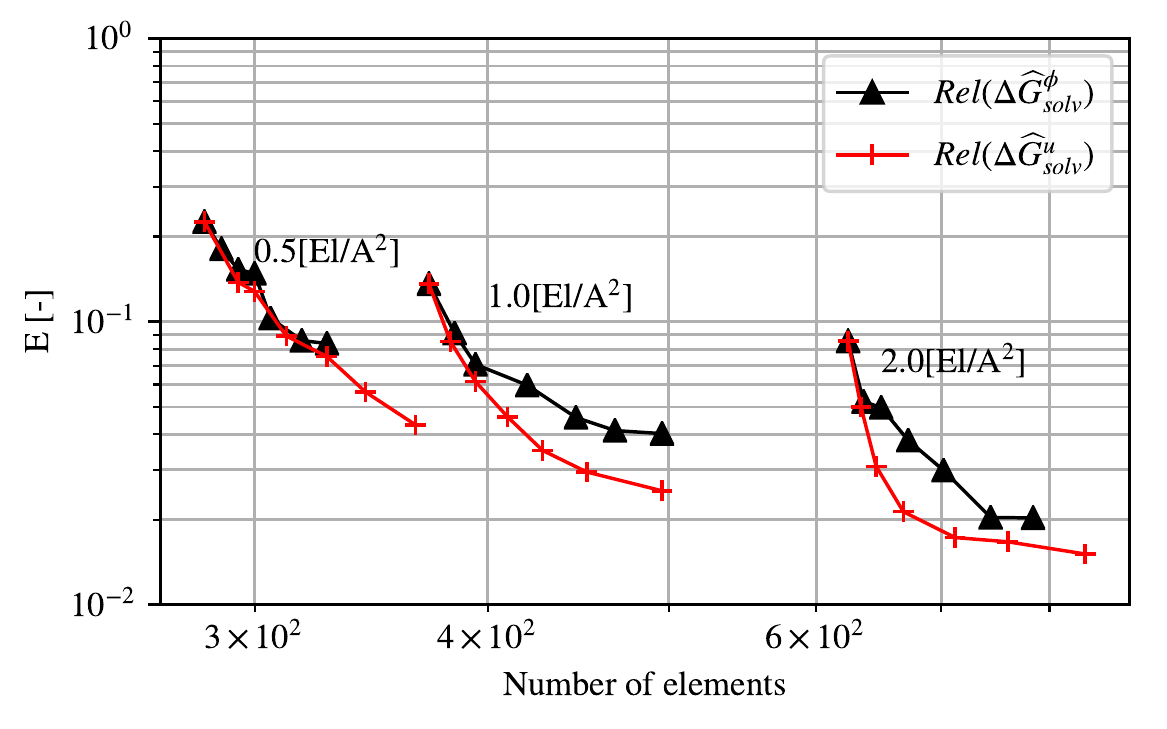} 
\end{center}
\vspace{0.25in}
\hspace*{3in}
{\Large
	\begin{minipage}[t]{3in}
		\baselineskip = .5\baselineskip
		Figure 15 \\
		Vicente Ramm, Jehanzeb H. Chaudhry, Christopher D. Cooper \\
		J.\ Comput.\ Chem.
	\end{minipage}
}
\clearpage

\begin{center}
			\includegraphics[width=0.32\textwidth]{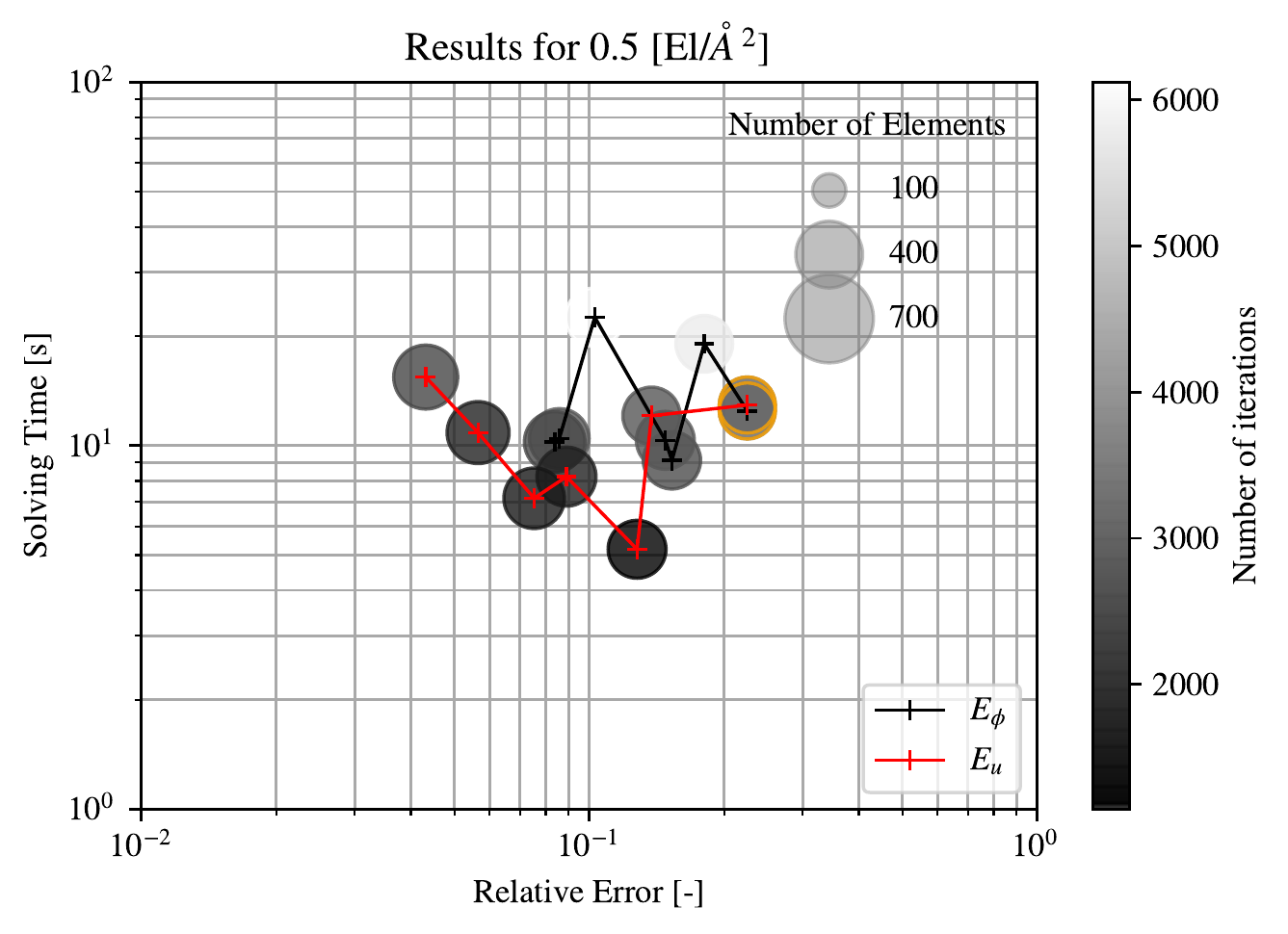} 
			\includegraphics[width=0.32\textwidth]{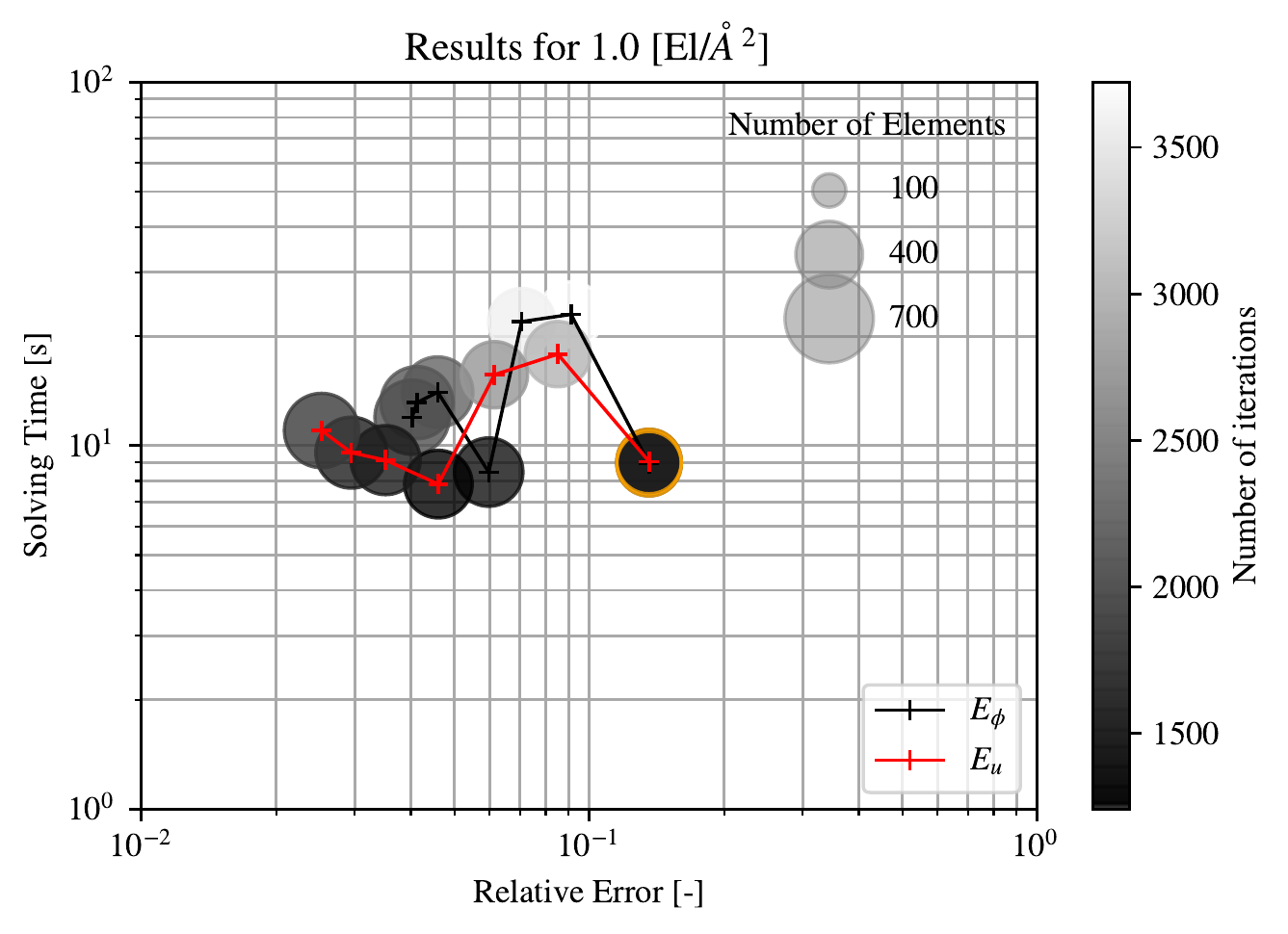} 
			\includegraphics[width=0.32\textwidth]{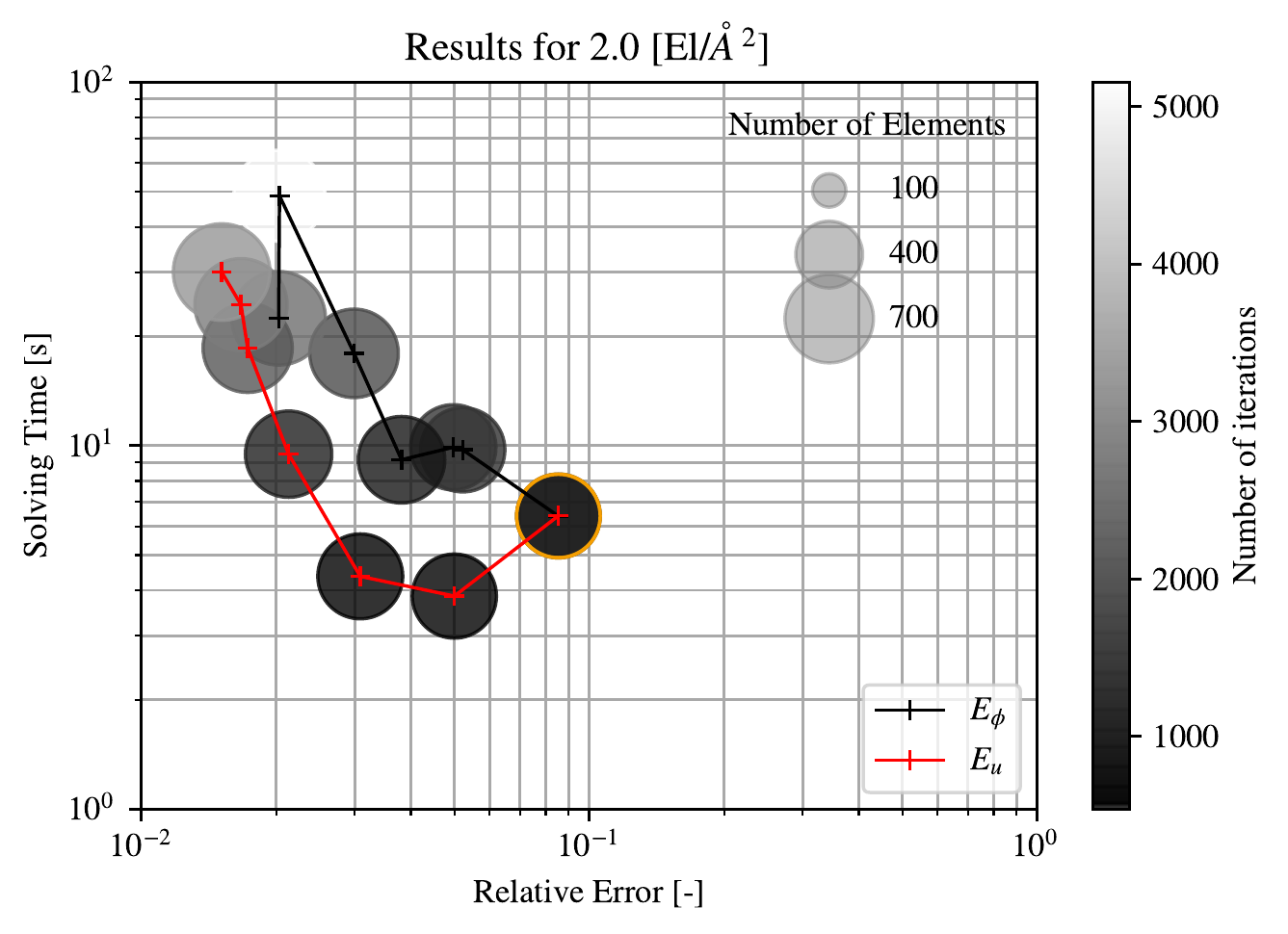} 
\end{center}
\vspace{0.25in}
\hspace*{3in}
{\Large
	\begin{minipage}[t]{3in}
		\baselineskip = .5\baselineskip
		Figure 16 \\
		Vicente Ramm, Jehanzeb H. Chaudhry, Christopher D. Cooper \\
		J.\ Comput.\ Chem.
	\end{minipage}
}

\clearpage
\begin{table}
	\centering
	\begin{tabular}{|ccrrrr|}
		\hline
		Density & $N_{\phi}$ & \multicolumn{1}{c}{$E_{\phi}$} & \multicolumn{1}{c}{$\gamma_{eff}^{\phi}$} & \multicolumn{1}{c}{$E_{u}$} & \multicolumn{1}{c|}{$\gamma_{eff}^{u}$} \\ \hline
		0.5 & 48    &  0.1461  & -0.089 & 1.5774   & -0.963 \\
		& 192       & -1.2512  & 0.764  & -0.4864  & 0.297  \\
		& 768       & -1.9297  & 1.178  & -1.2182  & 0.744  \\
		& 3072      & -2.2426  & 1.369  & -1.4882  & 0.908  \\
		& 12288     & -2.3300  & 1.422  & -1.5668  & 0.956  \\ 
		\multicolumn{3}{|c}{$\Delta \widehat{G}_{solv}=-4,9431$ kcal/mol}  & \multicolumn{3}{c|}{$\Delta G_{solv}=-6,5815$ kcal/mol}   \\ \hline
		1.0 & 86    & 0.0138   & 0.012  & 0.4806   & -2.501 \\
		& 344       & -0.2841  & -0.241 & -0.0217  & 0.113  \\
		& 1376      & -0.4393  & -0.373 & -0.1481  & 0.771  \\
		& 5504      & -0.4856  & -0.412 & -0.1815  & 0.944  \\
		& 22016     & -0.4980  & -0.423 & -0.1906  & 0.992  \\ 
		\multicolumn{3}{|c}{$\Delta \widehat{G}_{solv}=-3,5731$ kcal/mol}  & \multicolumn{3}{c|}{$\Delta G_{solv}=-3,7652$ kcal/mol}   \\ \hline
	\end{tabular}
	\caption{\label{tab:Results_on_gamma_eff}$\gamma_{eff}$ for methanol with 0.5 and 1 elements per \AA$^2$. $N_\phi$ is the number of elements in the mesh to compute $\phi$.}
\end{table}

\bibliography{bem_err,fastmethods,compbio}

\end{document}